\documentclass[11pt]{article}
\usepackage{amsfonts}
\usepackage{bm}
\usepackage{amsthm}
\usepackage{multirow}
\usepackage{authblk}
\usepackage{placeins}
\usepackage{booktabs}
\usepackage{mathtools}

\usepackage{amsmath}
\usepackage{amssymb}
\usepackage{graphicx}
\usepackage{color}

\usepackage{marginnote}
\usepackage{amstext,amsthm,amssymb,amsmath}
\usepackage{graphicx}
\usepackage{caption}
\usepackage{subcaption}

\usepackage[ruled,vlined]{algorithm2e}
\SetKwInput{AlgInput}{Input}
\SetKwInput{AlgOutput}{Output}
\SetKw{AlgReturn}{return}
\SetKw{AlgContinue}{continue}
\SetAlFnt{\footnotesize}
\SetAlCapFnt{\small}
\SetAlCapNameFnt{\small}

\renewcommand{\div}{\mathop{\rm div}\nolimits}

\hfuzz=\maxdimen
\usepackage{array}
\newcolumntype{C}[1]{>{\centering\arraybackslash}m{#1}}

\graphicspath{ {./pic/}}

\title{Multicontinuum Generalized Multiscale Finite Element Method (MC-GMsFEM). Theory and applications to upscaling of two-phase flow}
\date{}

\author[a]{Mohammed Al Kobaisi \thanks{Corresponding author}}
\author[a]{Dmitry Ammosov}
\author[b]{Yalchin Efendiev}
\author[c]{Wing Tat Leung}

\affil[a]{Chemical and Petroleum Engineering Department, Khalifa University of Science and Technology, Abu Dhabi, 127788, UAE}
\affil[b]{Department of Mathematics, Texas A\&M University, College Station, TX 77843, USA}
\affil[c]{Department of Mathematics, City University of Hong Kong, Hong Kong}

\begin{document}
\maketitle
\abstract{We develop a multicontinuum Generalized Multiscale Finite Element Method (MC-GMsFEM) for constructing coarse-scale models in heterogeneous media that simultaneously provide accurate numerical approximations and physically consistent macroscopic equations. Classical multiscale methods efficiently approximate fine-scale solutions on coarse grids using localized basis functions, but they do not offer a systematic pathway for deriving macroscopic governing equations. To overcome this limitation, we introduce a unified framework that integrates multiscale finite element constructions with multicontinuum representations.

The proposed method builds on the structure of GMsFEM and exploits a representation of multiscale basis functions that separates coarse variables and their gradients. We construct continuum-dependent basis functions using auxiliary fields defined through local problems with integral constraints, ensuring that each basis function is associated with a specific continuum. This leads to a decomposition of the coarse-scale solution into continuum variables and their gradients, establishing a direct connection between multiscale discretizations and multicontinuum homogenization.

Compared to existing multicontinuum approaches, the proposed framework provides greater flexibility in handling heterogeneous media with spatially varying numbers of continua and is naturally embedded within a standard finite element setting. This enables both systematic derivation of macroscopic equations and straightforward numerical implementation.

We apply the proposed method to the upscaling of two-phase immiscible flow in heterogeneous porous media, where multiple interacting continua, including mobile and trapped phases, arise. With the proposed approaches, we derive new
macroscopic models and show that if classical models are used, the errors can be large, especially for trapped continua. The MC-GMsFEM framework offers a robust and general approach for modeling complex multiscale systems and deriving effective macroscopic descriptions.}

\section{Introduction}

Multiscale problems arising in heterogeneous media are commonly addressed by constructing coarse-scale models that incorporate fine-scale information through multiscale basis functions. Over the past decades, a variety of multiscale methods have been developed and successfully applied \cite{chung2018constraint,chen2003mixed, hou1997multiscale, jenny2003multi, efendiev2009multiscale, efendiev2013generalized, maalqvist2014localization,chung2015mixed, chaabi2024algorithmic, peter2009multiscale, henning2013oversampling, owhadi2014polyharmonic,fu2025wavelet, xie2026robust, panasenko2018multicontinuum}. These approaches aim to approximate fine-scale solutions on coarse grids by designing basis functions that encode subgrid heterogeneity. As a result, they provide efficient coarse-grid discretizations of global problems and have been widely used both as standalone solvers \cite{efendiev2013generalized} and in conjunction with temporal discretization techniques \cite{chung2021contrast}.

Despite their success, existing multiscale methods do not generally provide a systematic framework for deriving physically consistent macroscopic differential equations. In particular, while they deliver accurate numerical approximations, the connection between the coarse-scale representation and the underlying continuum models is not clear as macroscopic degrees of freedom in multiscale methods may lack spatial continuity. To address this limitation, multicontinuum homogenization techniques have been proposed in earlier works. 
However, these approaches become difficult to apply to complex systems in which the number of interacting continua varies both spatially and dynamically. Moreover, their implementation requires the solutions of special cell problems.

The main objective of this paper is to develop a unified framework, referred to as multicontinuum Generalized Multiscale Finite Element Method (MC-GMsFEM), that constructs conforming multiscale basis functions while simultaneously enabling the derivation of macroscopic equations. The proposed framework is particularly suited for complex applications such as two-phase flow in heterogeneous media.

To motivate our approach, we briefly review the Generalized Multiscale Finite Element Method (GMsFEM). Consider the elliptic problem
\[
L(u)=f, \quad \text{with} \quad L(u)=-\nabla \cdot (\kappa(x)\nabla u),
\]
where $\kappa(x)$ is a highly heterogeneous coefficient. In GMsFEM, multiscale basis functions are constructed by solving local spectral problems in coarse neighborhoods associated with coarse-grid vertices $x_j$ \cite{galvis2010domain,efendiev2011multiscale}. The resulting basis functions take the form
\[
\phi_i^j = \zeta_i^j  \xi^j,
\]
where $\zeta_i^j$ are eigenvectors and $\xi^j$ are multiscale partition of unity functions enriched by local fine-scale information. It is well known that the multiscale basis functions admit the representation
\begin{equation}
\label{eq:eqn1}
\xi^j = \phi_0^j + Z_i^j \cdot \nabla \phi_0^j,
\end{equation}
where $\phi_0^j$ is a standard partition of unity function and $Z_i^j$ represents fine-scale corrections obtained from local problems (local bubble functions), \cite{biezemans2023non}. Consequently, the enriched basis functions can be expressed as
\begin{equation}
\label{eq:eqn2}
\phi_i^j = \zeta_i^j  \phi_0^j +\zeta_i^j Z_i^j\cdot \nabla  \phi_0^j.
\end{equation}
In many practical applications, as those considered in the paper,
$\zeta_i^j$  are slowly varying with respect to the spatial variable ($x$). 
On the other hand, $Z_i^j$ accounts for features that can be localized
(such as inclusions) and can be assumed to be smoothly varying 
with respect to the spatial variable.
In this
case, 
\begin{equation}
\label{eq:eqn3}
\phi_i^j =
N_i(x)\phi_0^j + M_i(x)\cdot \nabla \phi_0^j,
\end{equation}
where $N_i$ and $M_i$ are functions defined locally within each coarse block.
Using this representation, the coarse-scale solution can be written as
\[
u_H = \sum_{i,j} u_i^j \phi_i^j = N_i(x) U_i + M_i(x)\cdot \nabla U_i,
\]
where $U_i = \phi_0^j u_i^j$ represents the $i$-th continuum variable. This decomposition provides a direct bridge between multiscale basis functions and macroscopic continuum descriptions, enabling the derivation of coarse-scale governing equations.

In this work, we build on this structure and propose the MC-GMsFEM framework. The key idea is the construction of multiscale basis functions using continuum-dependent functions $\psi_j$, which act as characteristic indicators of different continua. These functions can be identified, for example, through local spectral problems, although their proper ordering remains a nontrivial task. In this paper, we assume that such continua can be identified, as is the case in applications such as two-phase flow.
Given the continuum functions $\psi_j$, we construct auxiliary functions $N_i$ by solving
\[
L(N_i) = r_i,
\]
subject to the constraints
\[
\int N_i \psi_j = \delta_{ij} \int \psi_j.
\]
These constraints ensure that each $N_i$ has unit average over continuum $i$ and vanishing averages over all other continua. The corresponding multiscale basis functions are then defined as
\[
\phi_i^j = N_i \phi_0^j + M_i \cdot \nabla \phi_0^j,
\]
where the correction terms $M_i \cdot \nabla \phi_0^j$ are obtained from local problems with zero boundary conditions and satisfy orthogonality conditions with respect to the continua. These corrections can be regarded as bubble functions, which guarantee that the basis functions satisfy the local problems.

The proposed framework is closely related to multicontinuum homogenization \cite{efendiev2023multicontinuum}, where the solution is approximated as
\[
u \approx N_i U_i + M_i \cdot \nabla U_i.
\]
However, there are several important distinctions. First, MC-GMsFEM provides a systematic way to construct multiscale basis functions that naturally yield this decomposition within a finite element framework. Second, it offers greater flexibility in handling heterogeneous media with spatially varying numbers of continua. Finally, the method is directly compatible with standard numerical implementations while simultaneously enabling the derivation of macroscopic equations.

A key application considered in this work is the upscaling of two-phase immiscible flow. At the macroscopic level, such systems are typically described 
(in the absence of capillary and other effects,
see e.g., \cite{helmig1997multiphase,chen2019homogenization})
 by momentum equations of the form
\[
v_i = k_r(S_i)\kappa(x)\nabla p, \quad i=1,2,
\]
coupled with mass conservation laws
\[
\partial_t S_i + \nabla \cdot v_i = 0, \quad \sum_i S_i = 1.
\]
These equations originate from fine-scale two-phase Stokes flows with complex interface dynamics in perforated domains. At the coarse scale, multiple interacting continua arise, including mobile and trapped phases (see Figure \ref{fig:two_phase_fine_results}). The MC-GMsFEM framework provides a natural mechanism to represent these continua and incorporate their effects into macroscopic models through additional continuum variables. In our formulation, we treat long trapped fluid regions as continua and show that the classical treatment via Darcy's law will give large errors. Our proposed multicontinuum approaches provide more accurate results.


The paper is organized as follows. In the next section, 
we present the general framework of MC-GMsFEM. In Section 3, we show the application of MC-GMsFEM for elliptic equations with heterogeneous coefficients and perforated domains. Finally, in Section 4, we present the application of MC-GMsFEM for two-phase flow equations, where we derive macroscopic equations and show detailed numerical comparisons.

\section{Multicontinuum GMsFEM}

Here, we present a general framework for multicontinuum GMsFEM. Throughout, we use Einstein notation interchangeably with standard notation when clarification is needed.

\begin{itemize}

\item We define multicontinuum functions $\psi_i$ (for the $i$-th continuum) by identifying continua globally.

\item 
We construct local functions $N_i$ (representing the $i$-th continuum)
 by solving local problems with average constraints
in local or extended or 
global domains by imposing the constraints on averages using appropriate
multicontinuum functions in the form 
\[
\int_K N_i \psi_j=\delta_{ij}\int_K \psi_j,
\]
 where $K$ is a coarse block. The functions $N_i$ have average $1$ in continuum $i$
and $0$ otherwise. These functions are used in constructing multicontinuum basis functions.

\item We modify $N_i$ functions by multiplying them by the standard
partition of unity functions $\phi_0^j$ and 
  construct basis functions $\phi_i^j$ in the form of
\begin{equation}
\label{eq:mcbasis}
\phi_i^j = N_i \phi_0^j + z_i^j,
\end{equation}
where $z_i^j$ vanishes on the boundary of coarse-grid block
and solves local problems with some constraints such as 
$\int_K z_i^j \psi_k=0$. Here, $\phi_0^j$
are standard finite element basis functions for each coarse node $j$.
We can show that $z_i^j$ can be written as
$z_i^j=M_i \cdot\nabla \phi_0^j$ 
 and we derive macroscopic
equations for this representation of $z_i^j$. 
For triangular macroscopic elements, this holds. For quadrilateral elements,
we need an additional constraints for $z_i^j$ when solving, where these
constraints are easily derived from linear dependence of rows of matrix
$\nabla_k \phi_0^j$. Another option is to enrich macroscopic model
by adding higher-order derivatives, e.g., mixed derivatives $\nabla^2_{12}\phi_0^j$ for 
quadrilateral elements, i.e., 
$z_i^j=M_i \cdot\nabla \phi_0^j + \widetilde{M_i}: \nabla^2_{kl}\phi_0^j$.
We will investigate such models in our future works.
In general, the multiscale solution, 
$u_H=\sum_{i,j} u_i^j \phi_i^j$,
obtained using the basis functions \eqref{eq:mcbasis} will
satisfy $\int_K u_H \psi_k/\int_K \psi_k=U_k$, 
where $U_k=\sum_j u_k^j \phi_0^j$, i.e.,
the multiscale solution will approximate the average $k$-th continuum
solution, $\int_K u\psi_k/\int_K \psi_k$, which is important for applications.

We remark that, for more accurate subgrid modeling, $z_i^j$ is
generally nonlocal \cite{chung2018constraint} and may have support 
extending over several coarse-grid blocks. This alters the 
local structure of the macroscopic model and naturally leads 
to a nonlocal macroscopic formulation.

\item We seek the solution of the global problem as
\begin{equation}
\label{eq:mcrep}
u_H = \sum_{i,j} \phi_i^j u_i^j=\sum_i (N_i U_i + M_i \cdot \nabla U_i),
\end{equation}
where $U_i=\sum_j u_i^j\phi_0^j$ and formulate macroscopic equations for $U_i$'s.
This representation is important in deriving macroscopic equations.

\end{itemize}

\section{Elliptic equation}

Next, we show the applications of MC-GMsFEM to elliptic equations.
We consider elliptic equations in heterogeneous media and in perforated domains. 
We present the derivation of
multicontinuum equations in heterogeneous domains in the paper
and the derivation of multicontinuum equations in perforated domains
in Appendix \ref{sec:ell_perfor}. We present some representative numerical results.

\subsection{The derivation of multicontinuum equations}

We consider
\begin{equation}
-\nabla \cdot (\kappa \nabla u)=f,
\end{equation}
with the weak form (assuming homogeneous Dirichlet boundary conditions)
\begin{equation}
a(u_H,v_H)=f(v_H),
\end{equation}
where
\begin{equation}
a(u_H,v_H)=\int_\Omega \kappa \nabla u_H \cdot \nabla v_H,
\qquad
f(v_H)=\int_\Omega f v_H.
\end{equation}

Using multicontinuum GMsFEM (see multiscale basis functions), we have
(see \eqref{eq:mcbasis} and \eqref{eq:mcrep})
\begin{equation}
u_H=\sum_i N_i U_i + \sum_i M_i \cdot \nabla U_i,
\ \ 
v_H=\sum_j N_j V_j + \sum_j M_j \cdot \nabla V_j.
\end{equation}
We refer to
\eqref{eq:Ni_ell} for definition of $N_i$
used in numerical simulations and \eqref{eq:Mi_ell} for the basis functions, which define $M_i$.
For simplicity, we omit $H$ index in macroscopic variables.
We have
\begin{equation}
\nabla u_H
\approx
\sum_i
\left(
\nabla N_i U_i + B_i \nabla U_i
\right),
\qquad
B_i := N_i I + \nabla M_i,
\nonumber
\end{equation}
where we ignore higher-order terms including the second order derivatives of
macroscopic variables.
Similar expression holds for $\nabla v_H$.
Substituting this into the variational formulation and
dividing the integrals over coarse regions, we have
\begin{equation}
\begin{split}
a(u_H,v_H)
=
\sum_K \sum_{i,j}
\int_K
\kappa
\Big[
(\nabla N_i \cdot \nabla N_j)U_iV_j
+
(\nabla N_i \cdot B_j)U_i\nabla V_j\\
+
(B_i \nabla U_i \cdot \nabla N_j)V_j
+
(B_i \nabla U_i)\cdot(B_j \nabla V_j)
\Big],
\nonumber
\end{split}
\end{equation}
where $K$ is a coarse grid block.
We define
\begin{equation}
\begin{split}
\alpha_{ij}:={1\over |K|}\int_K \kappa\,\nabla N_i\cdot \nabla N_j,
\quad
\beta_{ij}:={1\over |K|}\int_K \kappa\, \nabla N_i \cdot B_j,\\
\gamma_{ij}:={1\over |K|}\int_K \kappa\, B_i^T \nabla N_j,
\quad
\theta_{ij}:={1\over |K|}\int_K \kappa\, B_i^T B_j,
\nonumber
\end{split}
\end{equation}
where we avoid extra indices for some macroscopic quantities
for simplicity (e.g., $\theta_{ij}$ is a tensor).
For the forcing term,
\begin{equation}
\nonumber
f(v_H)
=
\sum_K \int_K
f
\left(
N_jV_j+M_j\cdot \nabla V_j
\right),
\nonumber
\end{equation}
with
\begin{equation}
\nonumber
F_j:={1\over |K|}\int_K fN_j,
\qquad
G_j:={1\over |K|}\int_K fM_j.
\nonumber
\end{equation}
The macroscopic variables are assumed to be smooth and can be
approximated by locally constant values within each coarse block. We use the same notation for globally defined macroscopic variables. Thus, we have (assuming Einstein summation)
\begin{equation}
\int_\Omega (\alpha_{ij}U_iV_j
+
(\beta_{ij} U_i) \cdot \nabla V_j
+
\gamma_{ij} \cdot \nabla U_i V_j
+
\theta_{ji} \nabla U_i \cdot \nabla V_j)
=
\int_\Omega\left(
F_jV_j+G_j\cdot \nabla V_j
\right).
\nonumber
\end{equation}
After integration by parts, 
\begin{equation}
\alpha_{ij}U_i
-
\nabla\cdot(\theta_{ji}\nabla U_i)
-
\nabla\cdot(\beta_{ij}U_i)
+
\gamma_{ij}\cdot \nabla U_i
=
F_j-\nabla\cdot G_j.
\nonumber
\end{equation}

\subsection{Construction of multiscale basis functions}

Let $\mathcal{T}_H$ be a coarse grid consisting of $K_d$. We denote by $x_l^H$ the coarse-grid vertices. For each $x_l^H$, we define the coarse neighborhood
\begin{equation}\label{eq:elliptic_omega_l}
\omega_l = \bigcup \{K_d \in \mathcal{T}_H: x_l^H \in \overline{K}_d\}.
\end{equation}
In addition, we introduce the general notation $\omega_l^*$, which may denote $\omega_l$, an oversampled extension of $\omega_l$, or the global domain $\Omega$.

To obtain the multiscale basis functions, we first compute the auxiliary functions $N_i^l$. Several approaches can be used for this purpose, including some spectral problems and specially designed boundary value problems. In this work, they are defined through constrained problems with coarse-block average constraints in $\omega_l^*$. More precisely, $N_i^l$ is obtained by solving the following problem
\begin{equation}
\label{eq:Ni_ell}
\begin{split}
&\int_{\omega_l^{*}} \kappa \nabla N_i^l \cdot \nabla v - \sum_{K_d \subset \omega_l^*} \sum_j \Lambda_{ij}^{ld} \int_{K_d} \psi_j v = 0,\\
&\int_{K_d} N_i^l \psi_j = \delta_{ij} \int_{K_d} \psi_j, \quad \forall K_d \subset \omega_l^*,
\end{split}
\end{equation}
where $\Lambda_{ij}^{ld}$ denote Lagrange multipliers for imposing the average constraints. 

Next, we determine the auxiliary functions $M_i^l$ corresponding to the gradient effects. In our numerical implementation, instead of explicitly resolving $M_i^l$, we introduce $z_i^l = M_i^l \cdot \nabla \phi_0^l$. The function $z_i^l$ solves the following problem
\begin{equation}
\label{eq:Mi_ell}
\begin{split}
&\int_{\omega_l} \kappa \nabla (z_i^l + N_i^l \phi_0^l) \cdot \nabla v - \sum_{K_d \subset \omega_l} \sum_j \Theta_{ij}^{ld} \int_{K_d} \psi_j v = 0,\\
&\int_{K_d} z_i^l \psi_j = 0, \quad \forall K_d \subset \omega_l,\\
&z_i^l|_{\partial K_d} = 0, \quad \forall K_d \subset \omega_l,
\end{split}
\end{equation}
where $\Theta_{ij}^{ld}$ are Lagrange multipliers associated with the zero average constraints. In our numerical results, these constraints were especially important for the perforated case. We refer to \eqref{eq:mcbasis} for discussions on the representation of $z_i^l$ via $M_i^l$.

Finally, the multiscale basis functions are computed as
\begin{equation}\label{eq:phi_expansion_with_z}
\phi_i^l = N_i^l \phi_0^l + M_i^l \cdot \nabla \phi_0^l = N_i^l \phi_0^l + z_i^l.
\end{equation}

In the perforated case, the auxiliary functions are determined by solving constrained problems in the perforated computational domains. The boundary conditions on the perforations must also be taken into account. The constrained problems used in the numerical results are given by \eqref{eq:Ni_ell_perf} and \eqref{eq:Mi_ell_perf} in Appendix~\ref{sec:ell_perfor}.

\subsection{Numerical results}

In this subsection, we present numerical experiments to test the proposed MC-GMsFEM. The results are obtained using the auxiliary functions $N_i$ computed in the coarse neighborhoods $\omega_l$. The global approach provides similar results and omitted for brevity.

We consider two-dimensional elliptic problems with two continua in $\Omega = (0,1) \times (0,1)$. Both heterogeneous and perforated cases are studied. We write $\Omega = D \cup \mathcal{B} = D_1 \cup D_2 \cup \mathcal{B}$, where $D_1$ and $D_2$ correspond to the first and second continuum regions, respectively, and $\mathcal{B}$ denotes the set of perforations. For the heterogeneous case, we set $\mathcal{B} = \emptyset$.

In all tests, we set $f = 1$ and impose zero Dirichlet boundary conditions on $\partial D$, which includes perforation boundaries for the perforated case. The coefficient $\kappa$ is chosen according to the media case as follows:
\begin{itemize}
\item Heterogeneous case:
\begin{equation*}
\kappa = \begin{cases}
1, & x \in D_1, \\
10^4, & x \in D_2,
\end{cases}
\end{equation*}
\item Perforated case: $\kappa = 1$.
\end{itemize}

\begin{figure}[hbt!]
\centering
\setlength{\tabcolsep}{1pt}
\renewcommand{\arraystretch}{0.70}

\newcommand{\panel}[1]{\includegraphics[width=\linewidth]{#1}}
\newcommand{\rowlabel}[1]{\rotatebox[origin=c]{90}{\scriptsize #1}}

\begin{tabular}{@{} C{0.020\textwidth} C{0.235\textwidth} C{0.235\textwidth} C{0.235\textwidth} C{0.235\textwidth} @{}}
 & {\scriptsize Heterogeneous 1} & {\scriptsize Heterogeneous 2} & {\scriptsize Perforated 1} & {\scriptsize Perforated 2} \\

\rowlabel{Microstructure}
 & \panel{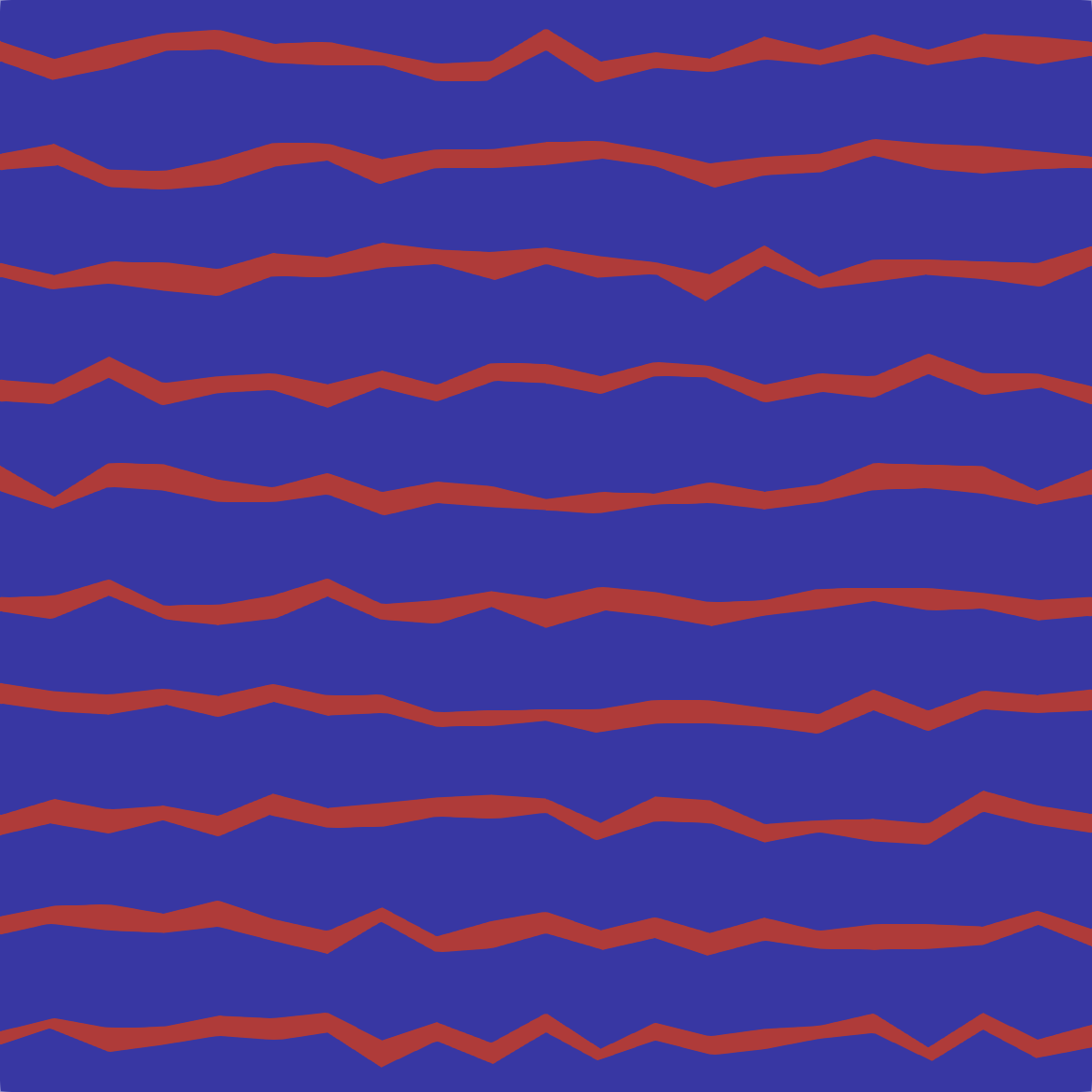}
 & \panel{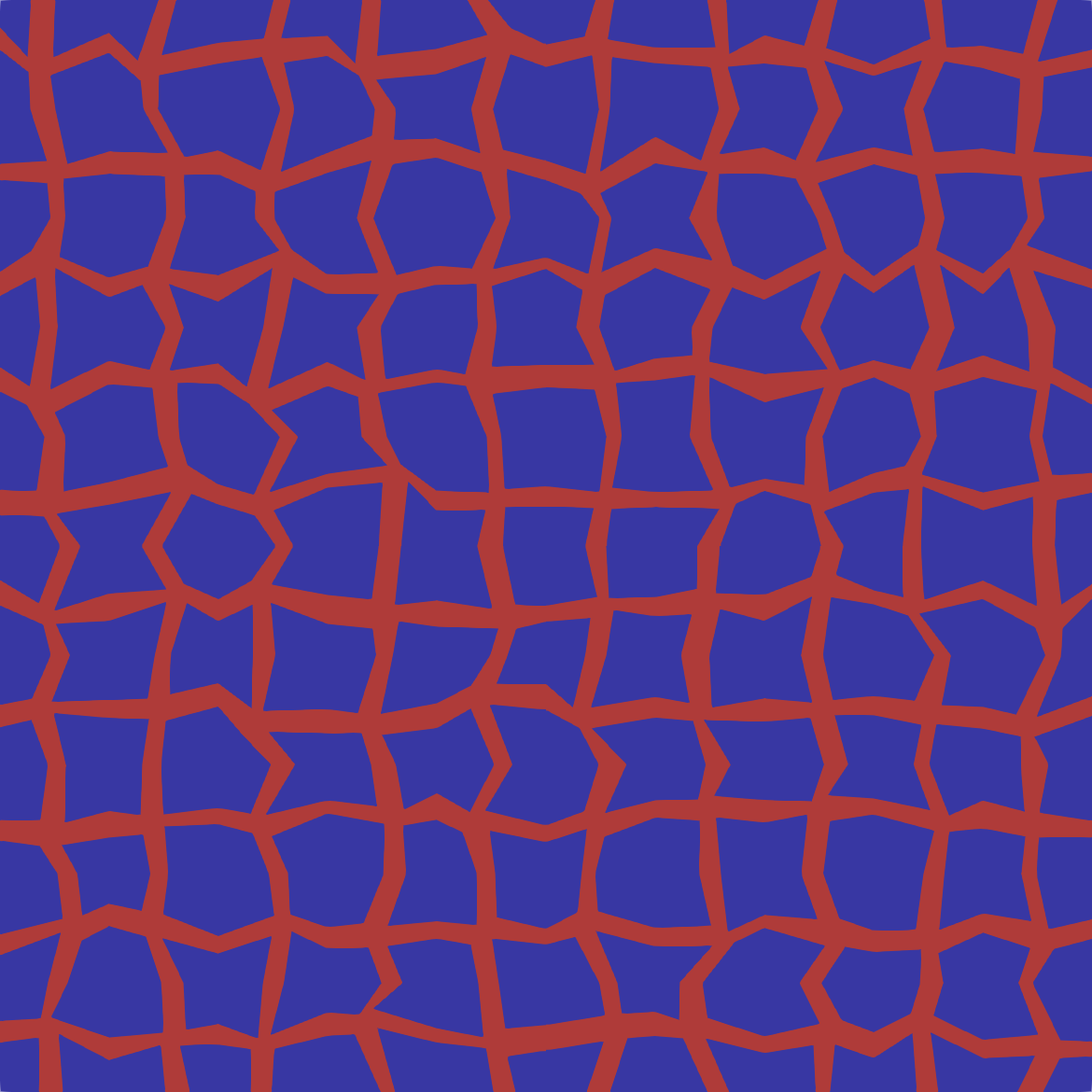}
 & \panel{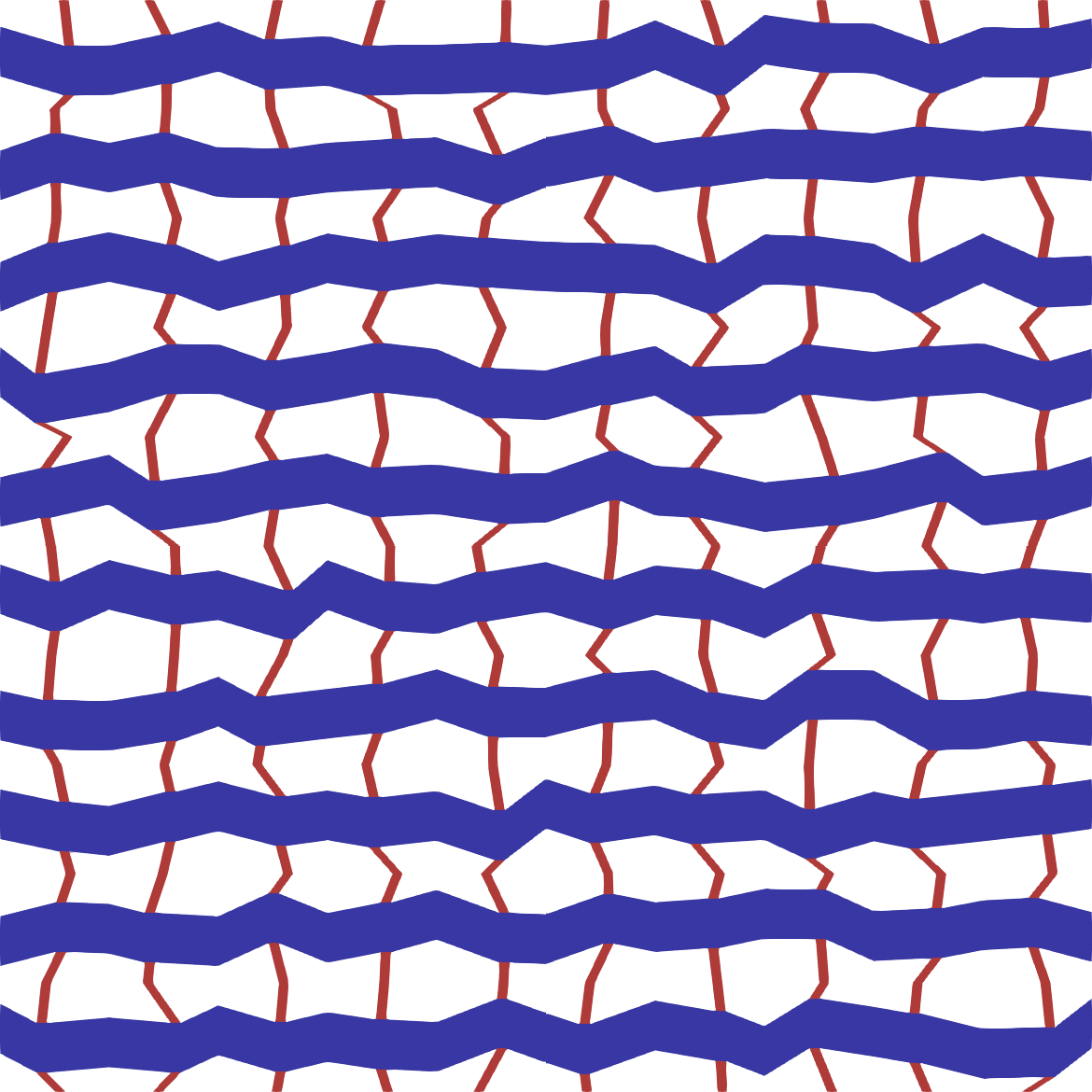}
 & \panel{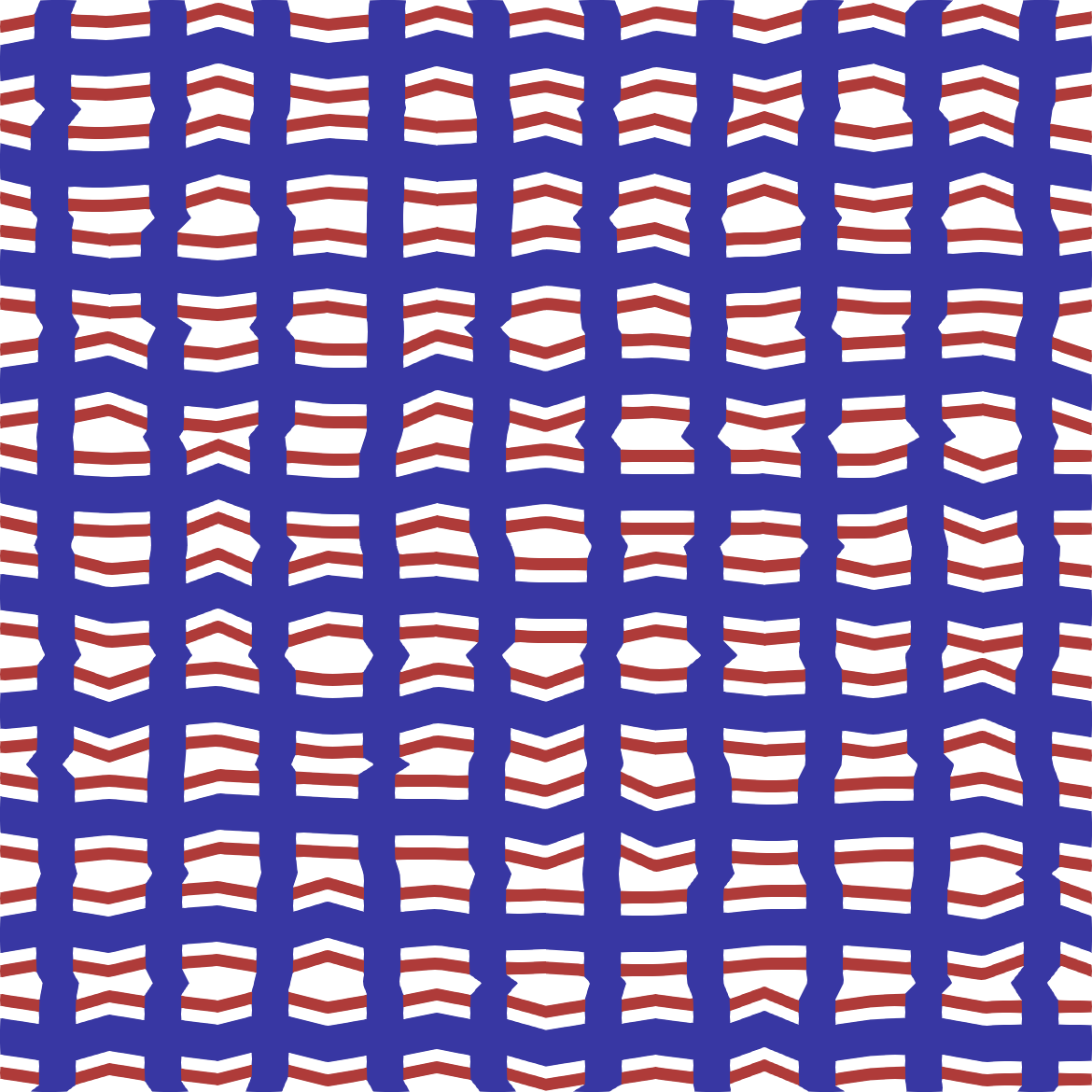} \\

\rowlabel{Reference solution}
 & \panel{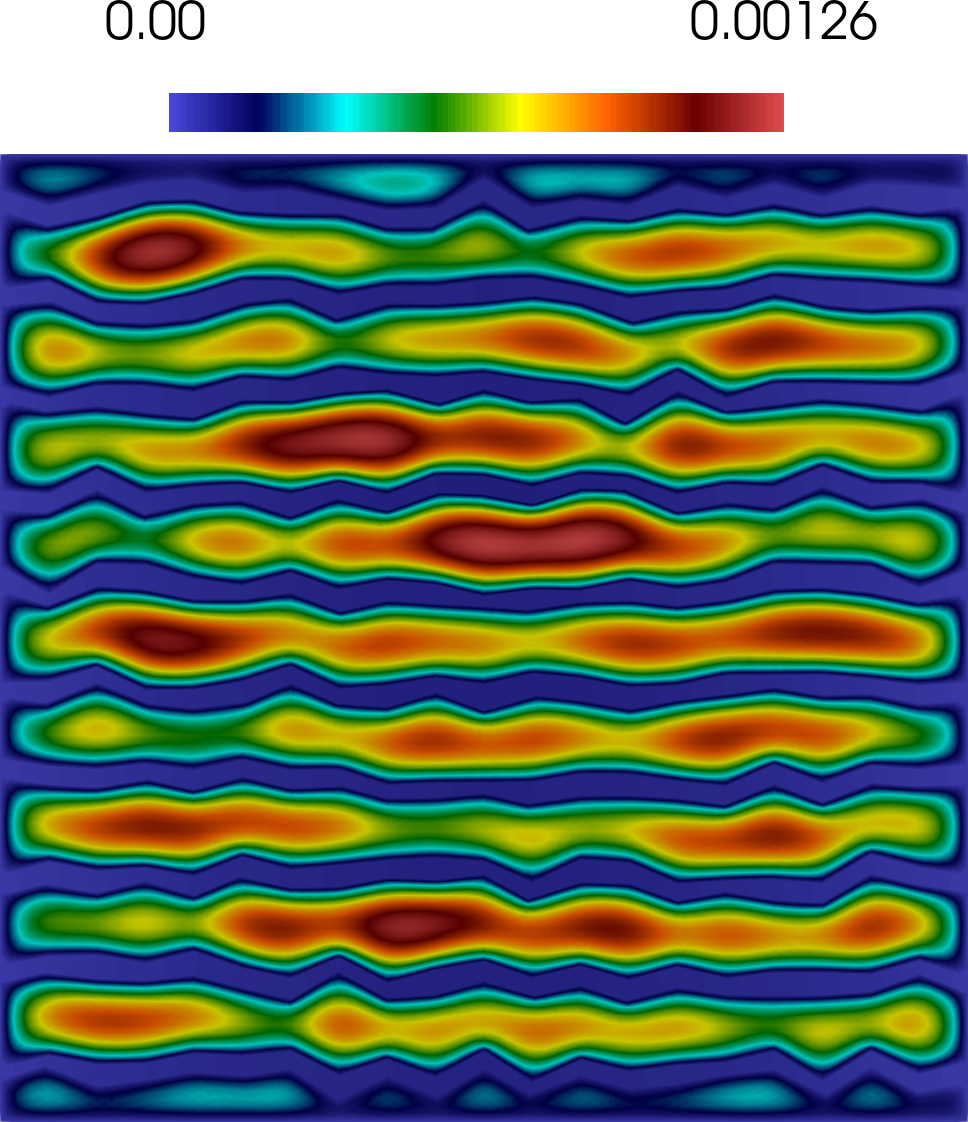}
 & \panel{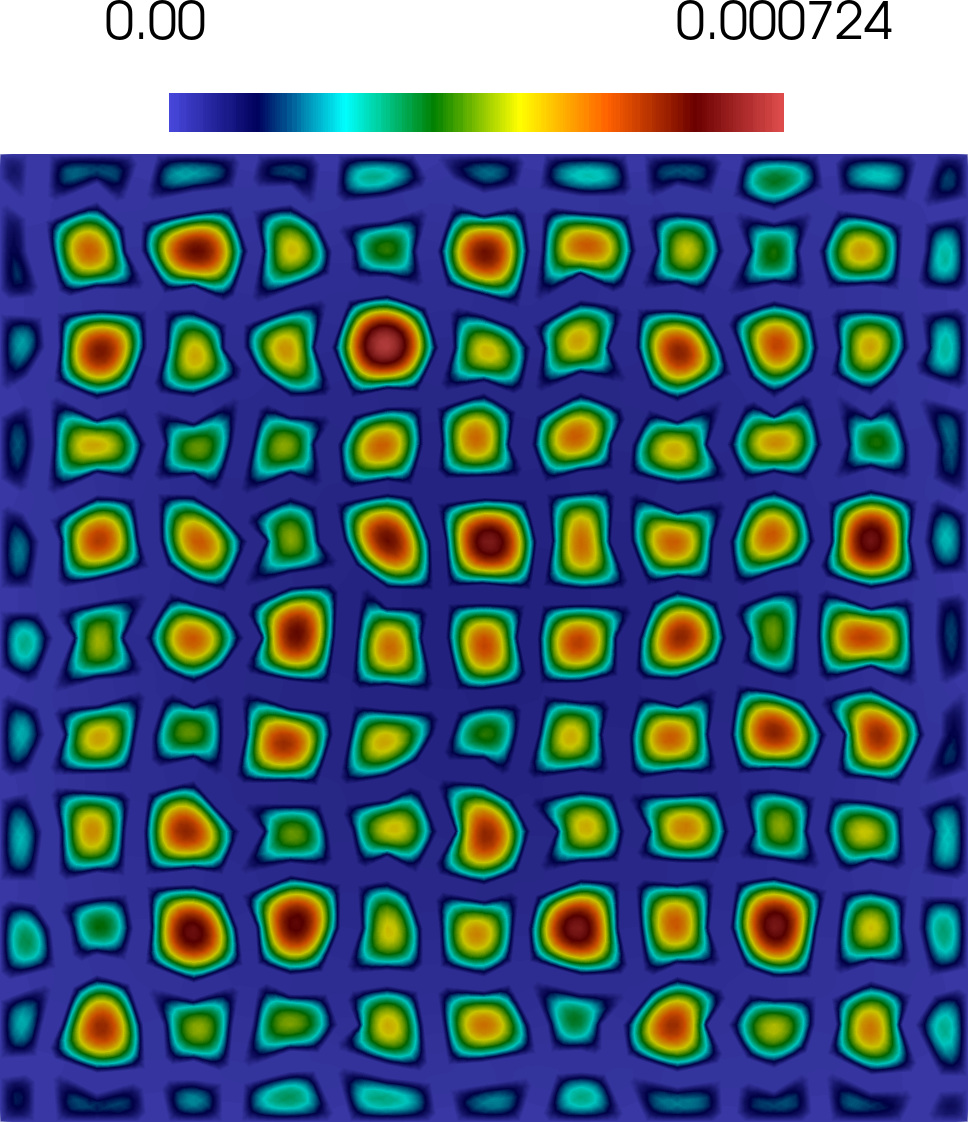}
 & \panel{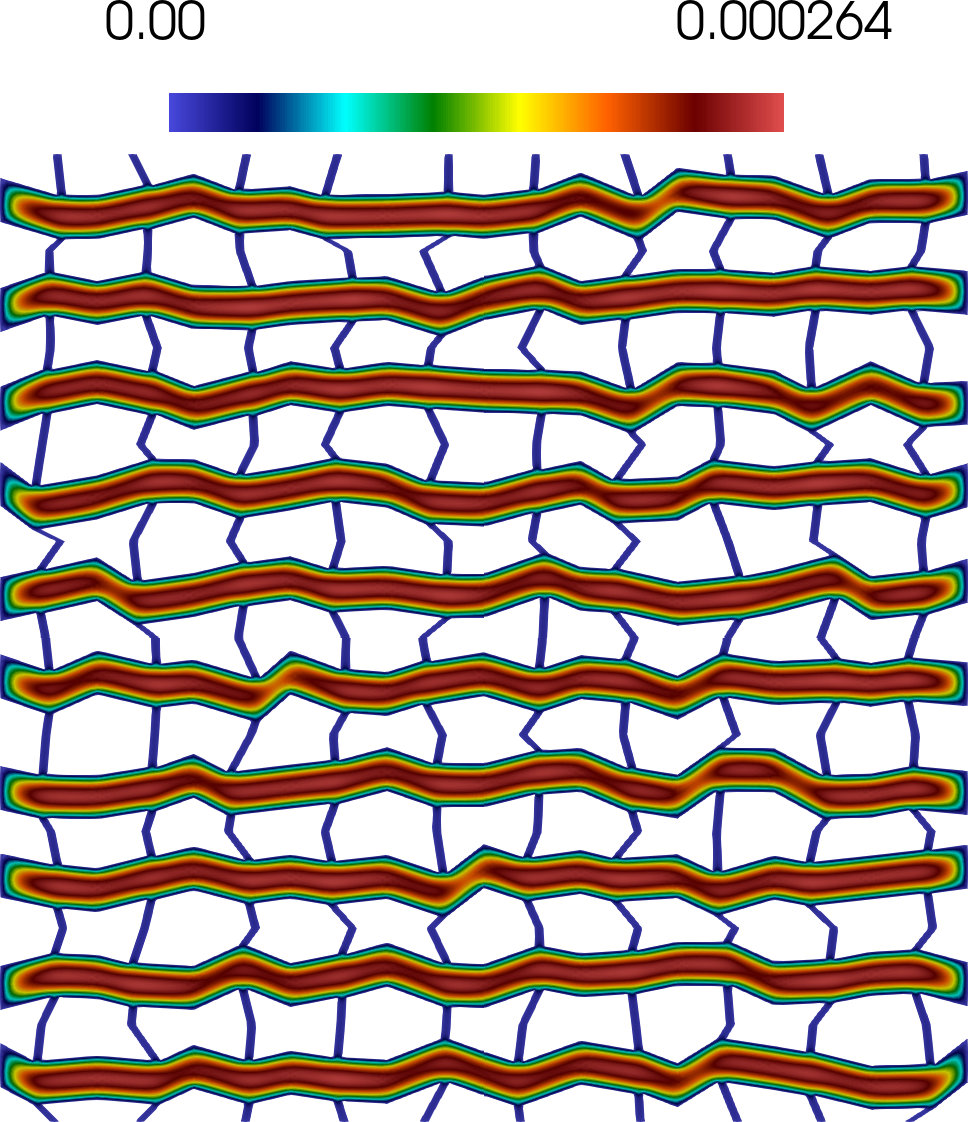}
 & \panel{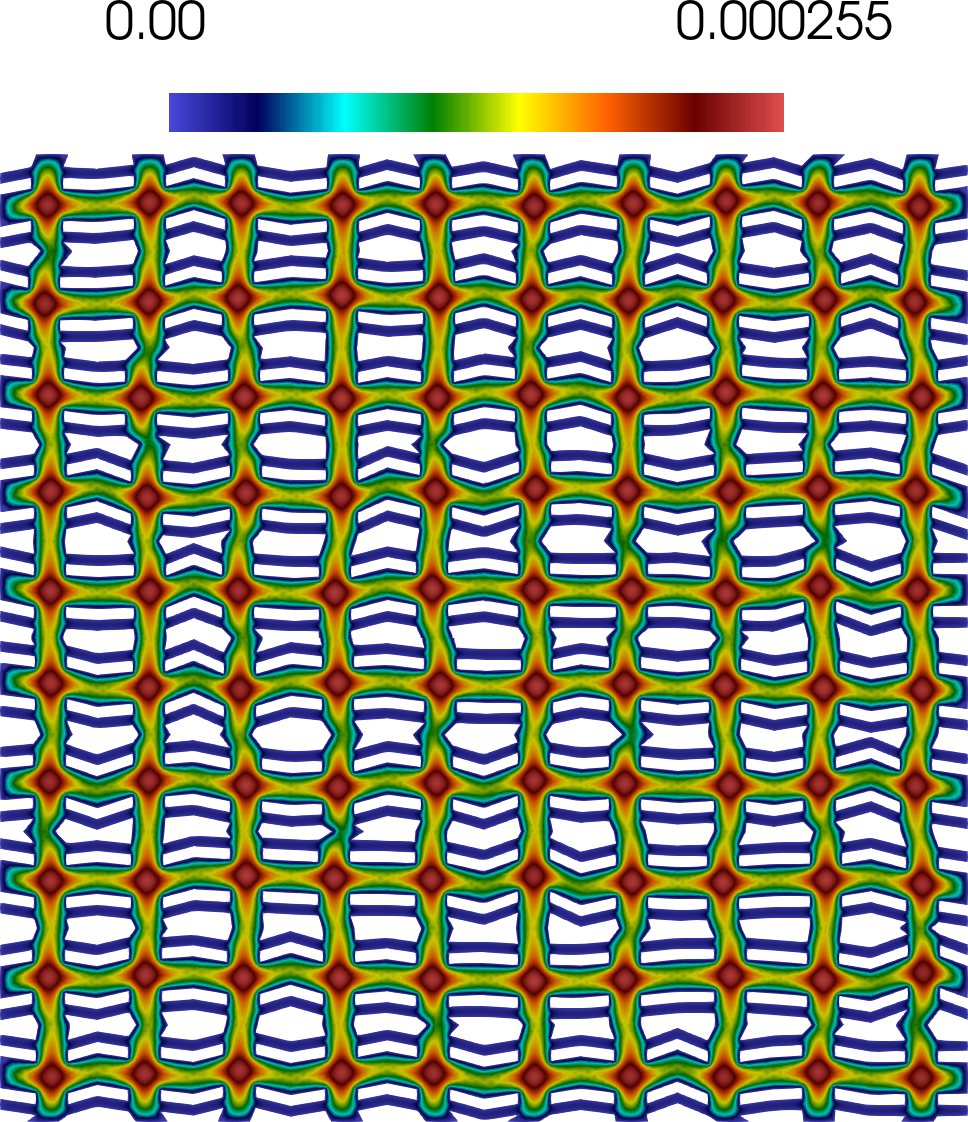} \\

\rowlabel{Multiscale solution}
 & \panel{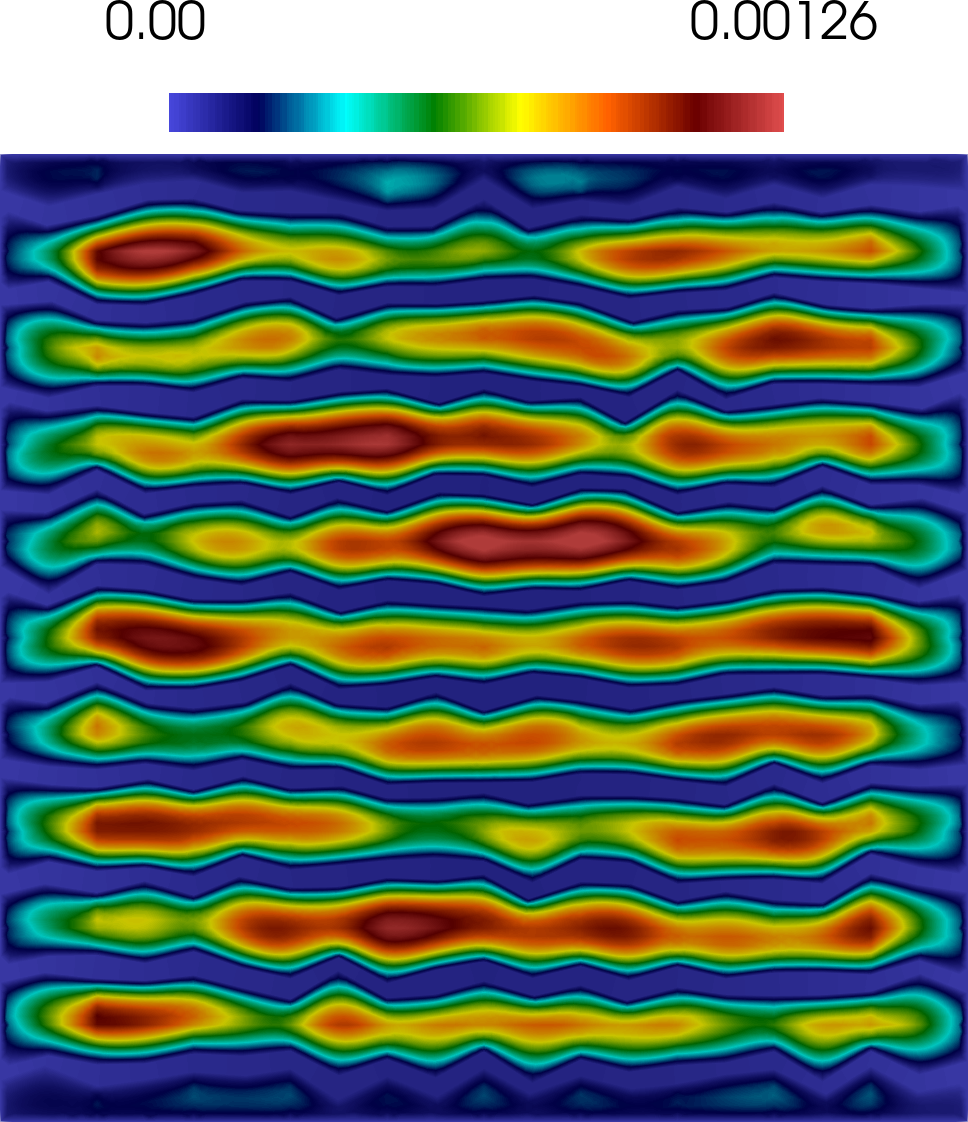}
 & \panel{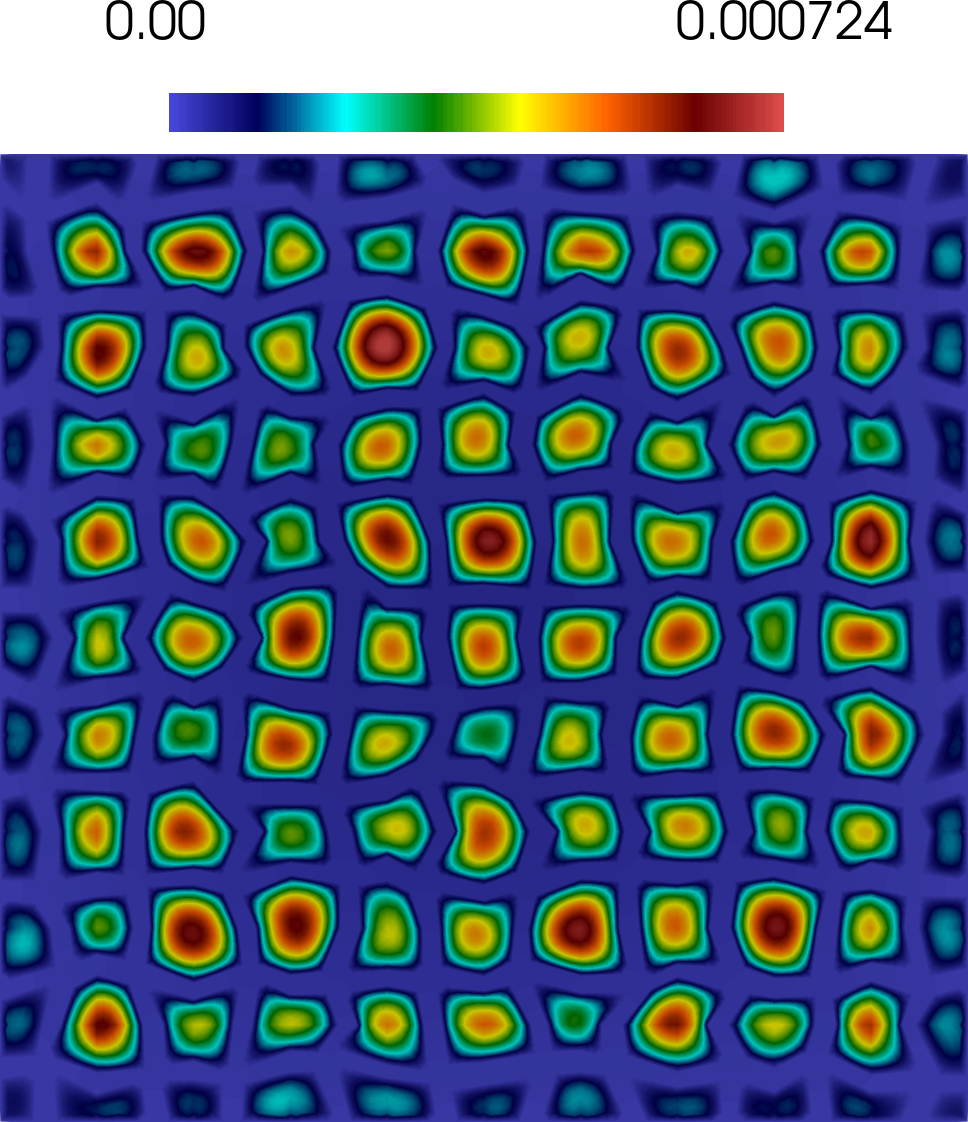}
 & \panel{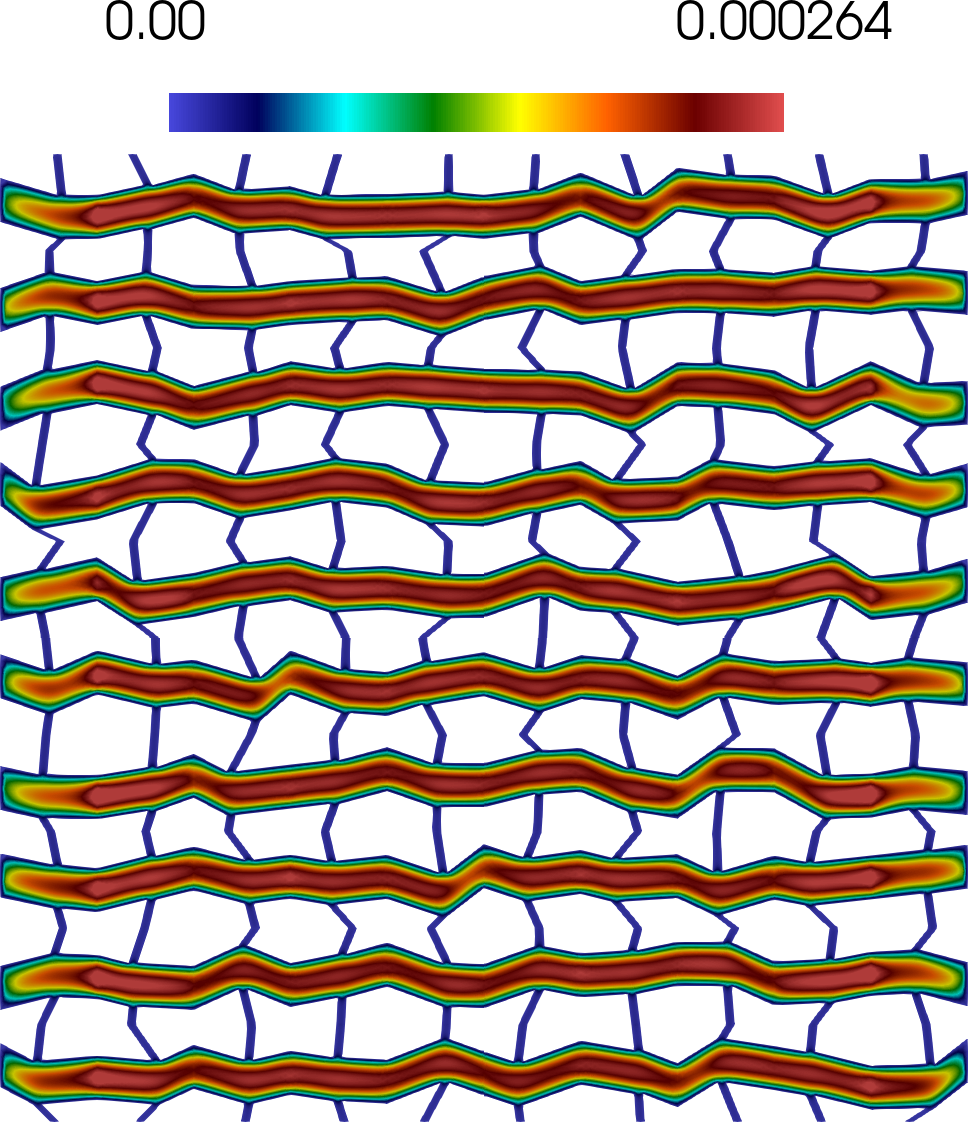}
 & \panel{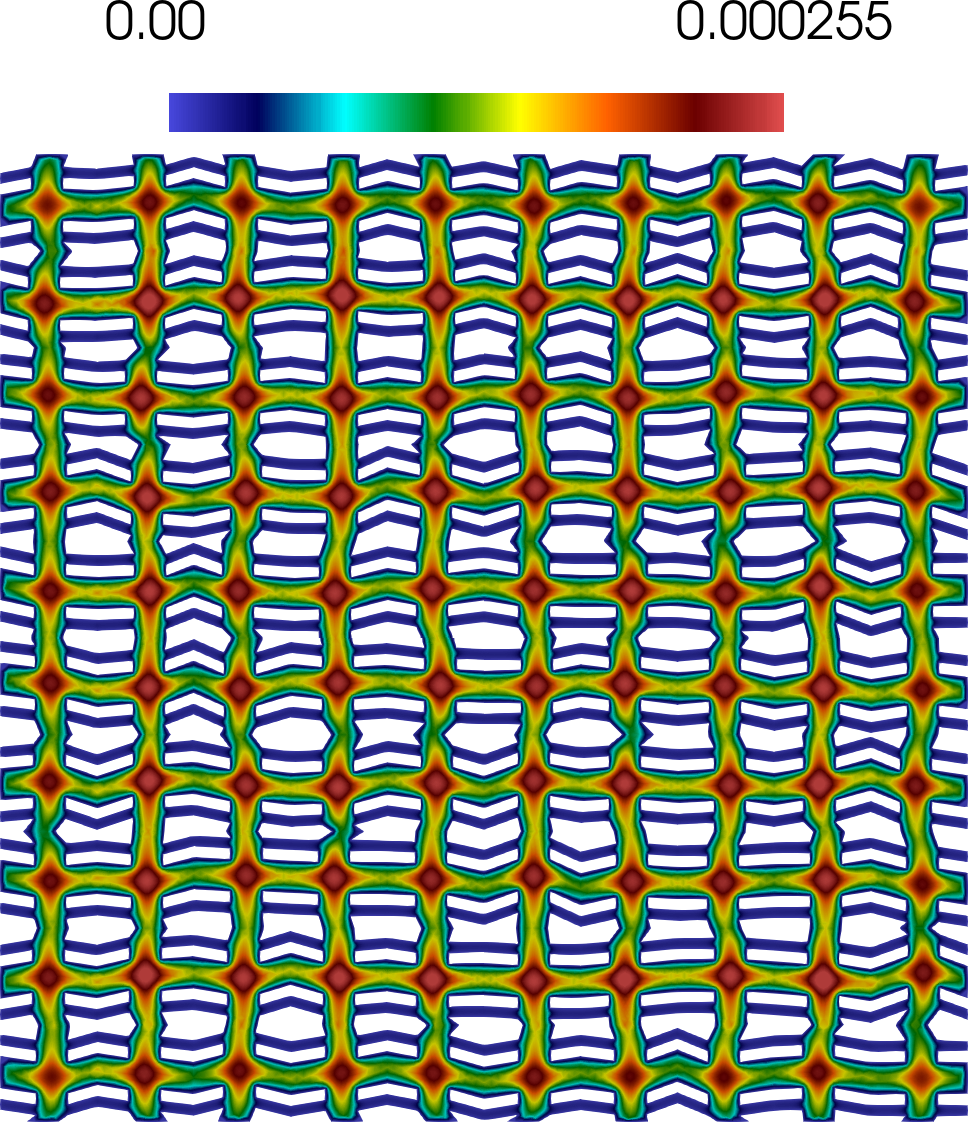}
\end{tabular}

\caption{
Numerical results for the elliptic problems. The rows show the microstructure, reference solution, and downscaled MC-GMsFEM solution. The four columns correspond to the heterogeneous and perforated microstructures.
}
\label{fig:elliptic_fine_grid_results}
\end{figure}

We examine two heterogeneous microstructures and two perforated microstructures, as shown in the first row of Figure~\ref{fig:elliptic_fine_grid_results}. The blue and red regions correspond to $D_1$ and $D_2$, respectively. The first heterogeneous microstructure is discretized using an unstructured fine grid with 29148 triangular cells and 14789 vertices, while the second heterogeneous microstructure is discretized using an unstructured fine grid with 32013 triangular cells and 16237 vertices. For the first perforated microstructure, we use an unstructured fine grid with 169558 triangular cells and 91503 vertices. The fine grid for the second perforated microstructure consists of 130152 triangular cells and 72700 vertices.

For both heterogeneous and perforated cases, we use a coarse grid consisting of $10 \times 10$ rectangular blocks, with coarse mesh size $H=1/10$. For each continuum $i$, the relative error of the MC-GMsFEM upscaled solution is computed as
\begin{equation}\label{eq:errors_upscaled}
e_{2}^{(i)} = \frac{
\sqrt{
\sum_K
\left|
\frac{1}{|K|}\int_K U_i
-
\frac{1}{|K \cap D_i|}
\int_{K \cap D_i} u^{\mathrm{ref}}
\right|^2
}
}{
\sqrt{
\sum_K
\left|
\frac{1}{|K \cap D_i|}
\int_{K \cap D_i} u^{\mathrm{ref}}
\right|^2
}
}
\times 100\%,
\end{equation}
where $u^{\mathrm{ref}}$ is the reference fine-grid solution, $U_i = \sum_j u_{i}^{j} \varphi_0^j$, and $u_{i}^{j}$ are the coarse-grid degrees of freedom of the MC-GMsFEM solution.

The MC-GMsFEM can also reconstruct the fine-grid solution from the computed coarse-grid solution. Therefore, for the elliptic problems considered here, we also compute the following downscaling errors
\begin{equation*}
\begin{gathered}
e_{2}^{\mathrm{down}} =
\frac{\sqrt{\int_D (u^{\mathrm{ref}} - u^{\mathrm{mc}})^2}}
{\sqrt{\int_D (u^{\mathrm{ref}})^2}}
\times 100 \%,\\
e_{a}^{\mathrm{down}} =
\frac{\sqrt{a(u^{\mathrm{ref}} - u^{\mathrm{mc}}, u^{\mathrm{ref}} - u^{\mathrm{mc}})}}
{\sqrt{a(u^{\mathrm{ref}}, u^{\mathrm{ref}})}}
\times 100 \%,
\end{gathered}
\end{equation*}
where $u^{\mathrm{mc}} = u_{i}^{j} \phi_i^j$ is the downscaled MC-GMsFEM solution and $a(u,v) = \int_D \kappa \nabla u \cdot \nabla v$.

For the fine-grid approximation, we use the finite element method with standard piecewise linear basis functions. The numerical implementation is based on the FEniCS computational package~\cite{logg2012automated}. The fine grids are generated using Gmsh \cite{geuzaine2009gmsh}, while visualization of the numerical results is performed using ParaView \cite{ayachit2015paraview}.

In the second and third rows of Figure~\ref{fig:elliptic_fine_grid_results}, we present the fine-scale distributions of the reference and downscaled multiscale solutions for the different microstructures. The visualizations of the upscaled solutions are omitted for brevity. One can see that the reference and downscaled solutions are very similar, indicating the high accuracy of the proposed MC-GMsFEM. Table~\ref{tabs:elliptic_errors} presents the corresponding relative errors. Both the upscaled and downscaled errors are small for all the considered microstructures.

\begin{table}[hbt!]
\caption{Relative errors of the MC-GMsFEM solutions for the elliptic problems.}
\label{tabs:elliptic_errors}
\centering
\begin{tabular}{lllll}
\hline
Microstructure & $e_{2}^{(1)}$ & $e_{2}^{(2)}$ & $e_{2}^{\mathrm{down}}$ & $e_{a}^{\mathrm{down}}$ \\
\hline
Heterogeneous 1 & 7.67 \% & 6.23 \% & 7.81 \% &  16.31 \% \\
Heterogeneous 2 & 6.11 \% & 15.24 \% & 6.90 \% & 19.55 \% \\ \addlinespace
Perforated 1 & 5.19 \% & 7.31 \% & 7.73 \% & 11.42 \% \\
Perforated 2 & 3.50 \% & 2.95 \% & 6.21 \% & 10.47 \% \\
\hline
\end{tabular}
\end{table}


\FloatBarrier
\section{Two-phase flow simulations}

Let $\Omega \subset \mathbb{R}^2$ can be represented as $\Omega = \Omega_f \cup \Omega_s$, where $\Omega_f$ and $\Omega_s$ are fluid and solid subregions, respectively. We consider the following coupled Stokes flow and transport equations in $\Omega_f$
\begin{equation}\label{eq:two_phase_fine_model}
\begin{split}
-\div \sigma &= f,\\
\div v &= 0,\\
\frac{\partial c}{\partial t} + v \cdot \nabla c &= 0,
\end{split}
\end{equation}
where $v(x, t)$ is the velocity, $p(x, t)$ is the pressure, $c(x, t)$ is the concentration, $f(x)$ is a body force, and
the Cauchy stress $\sigma(v, p, c)$ defined as
\[
\sigma(v, p, c)=2\mu(c)\varepsilon(v)-pI, \qquad \varepsilon(v)=\frac{1}{2}\left(\nabla v+(\nabla v)^T\right)
\]
with $\mu(c)$ denoting the dynamic viscosity, defined by
\begin{equation}\label{eq:mu_and_psi}
\begin{gathered}
\mu(c) = \sum_k \mu_k \psi_k(c), \quad
\psi_k(c) =
\begin{cases}
  1, & c \in I_k,\\
  0, & \text{otherwise},
\end{cases} \\
I_k = \{c: c_k \leq c < c_{k+1}\}.
\end{gathered}
\end{equation}
Here, $I_k$ denotes the interval of concentration values associated with the $k$-th fluid phase, $\psi_k(c)$ is its characteristic function, and $\mu_k$ is the dynamic viscosity of the $k$-th phase.

We complement \eqref{eq:two_phase_fine_model} with the boundary conditions
\begin{equation}\label{eq:two_phase_fine_bcs}
\begin{split}
v = v_{\mathrm{in}}, \quad c = c_{\mathrm{in}} \quad &\text{on } \Gamma_{l}, \\ 
v = v_{\mathrm{out}}, \quad &\text{on } \Gamma_{r}, \\ 
v \cdot n = 0, \quad \tau \cdot (\sigma n) + \beta (v \cdot \tau) = 0 \quad &\text{on } \Gamma_{s},
\end{split}
\end{equation}
where $v_{\mathrm{in}}(x, t)$ and $v_{\mathrm{out}}(x, t)$ are the inlet and outlet velocities, $c_{\mathrm{in}}(x, t)$ is the inlet concentration,
$\beta(c)$ is a friction coefficient, $n$ is the unit outward normal vector,
and $\tau$ is a unit tangent vector. 

We denote by $\Gamma_l$ and $\Gamma_r$ the left and right parts of the exterior boundary of $\Omega_f$, respectively, and define
\[
\Gamma_s = \partial \Omega_f \setminus (\Gamma_l \cup \Gamma_r),
\]
which denotes the slip boundary, including the boundaries of the perforations.

In addition, we complement the transport equation with the initial condition
\begin{equation}\label{eq:two_phase_fine_ics}
c(x, t_0) = c_0(x) \quad \text{in } \Omega_f.
\end{equation}

\subsection{Single pressure continuum. The relation to our numerical results and existing two-phase flow model}

We first consider a case of a single-pressure continuum, which is commonly
used. In this case, we assume a smooth pressure field (denote $p_0$)
 and keep the term
$\nabla p_0$ (instead of $\nabla_x p$ to indicate the single continuum)
 in the equation and only expand the velocity field and 
concentration 
via multiscale basis functions. As before, we 
 use for basis functions
for the velocity field
\[
\phi_i^{v, j} = N_i^v \phi_0^j + M_i^v \cdot \nabla\phi_0^j,
\]
where $N_i^v$ is a globally (or locally) defined solution with constraints
and $M_i^v$ is zero on the boundary of coarse regions
(see \eqref{eq:two_phase_N_i} for definition of $N_i^v$
used in numerical simulations and \eqref{eq:two_phase_phi_i} for the basis functions, which define $M_i^v$).
We will first derive macroscopic equations
for Stokes equations. The macroscopic equations for 
the transport equations are derived in \cite{chen2019homogenization} and will use the same formulation. 
Using multicontinuum GMsFEM (see multiscale basis functions), we have
(see \eqref{eq:mcbasis} and \eqref{eq:mcrep})
\[
v=\sum_{i,j} \phi_i^{v, j} v_i^j= \sum_i (N_i^v V_i + M_i^v\cdot \nabla V_i),
\]
where $V_i=\sum_j v_i^j \phi_0^j$.

Let $\Omega_f$ be the fluid region and $\Gamma_s$ be the slip boundary. 
The Stokes variational form is
\begin{equation}
\int_{\Omega_f}
2\mu\,\varepsilon(v):\varepsilon(w)
+
\int_{\Gamma_s}
\beta\, v_\tau \cdot w_\tau
+
\int_{\Omega_f}
\nabla p_0(x)\cdot w
=
\int_{\Omega_f}
f\cdot w
\end{equation}
with
\begin{equation}
\nabla\cdot v = 0,
\nonumber
\end{equation}
where $\epsilon(v)$ is a symmetric gradient.
We leave $\nabla p_0$ unchanged.

We assume the velocity is represented by enriched multicontinuum basis functions
\begin{equation}
\label{eq:vel_exp}
v(x)
=
N_i^v(x)V_i(x)
+
M_i^v(x)\cdot \nabla V_i(x),
\end{equation}
and test function
\begin{equation}
w(x)
=
N_j^v(x)W_j(x)
+
M_j^v(x)\cdot \nabla W_j(x),
\nonumber
\end{equation}
where 
$N_i^v(x)$ - velocity continuum partition/basis functions,
$M_i^v(x)$ - velocity enrichment/corrector basis,
$V_i(x)$ - continuum velocity variables,
 $W_j(x)$ - continuum test functions.
Then, we have (as before)
\begin{equation}
\nabla v =
\nabla(N_i^v V_i + M_i^v \cdot \nabla V_i) \approx
(\nabla N_i)V_i
+
B_i^v \nabla V_i
\nonumber
\end{equation}
$B_i^v = N_i^v I + \nabla M_i^v$ (with indices $B_i^{v, klnr}=N_i^{v, kl}\delta_{nr}+\nabla_r M_i^{v, kln}$)
and similar expansion for $w$.
Therefore,
\begin{equation}
\varepsilon(v)
\approx
\varepsilon(N_i^v)V_i
+
\varepsilon(B_i^v)\nabla V_i
\end{equation}

Next, we substitute the expansion into the variational form
\begin{align}
\nonumber
&
\int_{\Omega_f}
2\mu
\,
\varepsilon\left(
N_i^v V_i + M_i^v \cdot \nabla V_i
\right)
:
\varepsilon\left(
N_j^v W_j + M_j^v \cdot \nabla W_j
\right)
\\
\nonumber
&+
\int_{\Gamma_s}
\beta
\left(
N_i^v V_i + M_i^v \cdot \nabla V_i
\right)_\tau
\cdot
\left(
N_j^v W_j + M_j^v \cdot \nabla W_j
\right)_\tau
\\
\nonumber
&+
\int_{\Omega_f}
\nabla p_0(x) \cdot
\left(
N_j^v W_j + M_j^v \cdot \nabla W_j
\right)
=
\int_{\Omega_f}
f\cdot
\left(
N_j^v W_j + M_j^v \cdot \nabla W_j
\right).
\end{align}
As before, we divide the integral over coarse blocks, where macroscopic 
variables can be approximated by their mid-point values and denote
\begin{equation*}
\alpha_{ij}
=
\frac{1}{|K|} (\int_{K\cap\Omega_f}
2\mu\,\varepsilon(N_i^v):\varepsilon(N_j^v) + \int_{\Gamma_s\cap \overline{K}}
\beta (N_i^v)_\tau \cdot (N_j^v)_\tau),
\end{equation*}

\begin{equation*}
\gamma_{ij}
=
\frac{1}{|K|} (\int_{K\cap\Omega_f}
2\mu\,\varepsilon(B_i^v):\varepsilon(N_j^v)+\int_{\Gamma_s\cap \overline{K}}
\beta (M_i^v)_\tau \cdot (N_j^v)_\tau),
\end{equation*}

\begin{equation*}
\theta_{ij}
=
\frac{1}{|K|} (\int_{K\cap\Omega_f}
2\mu\,\varepsilon(B_i^v):\varepsilon(B_j^v)
+
\int_{\Gamma_s\cap \overline{K}}
\beta (M_i^v)_\tau \cdot (M_j^v)_\tau),
\end{equation*}

\begin{equation*}
\rho_j = \frac{1}{|K|} \int_{K\cap\Omega_f}
(N_j^v)^T,\,\ \ 
\eta_j
=
\frac{1}{|K|} \int_{K\cap\Omega_f}
f \cdot N_j^v,\,
\ \ 
\zeta_j
=
\frac{1}{|K|} \int_{K\cap\Omega_f}
 f \cdot M_j^v.
\end{equation*}

Assuming that  $\int_{\Omega_f}
\nabla p_0(x) \cdot (
N_j^v W_j + M_j^v \cdot \nabla W_j)\approx  \int_{\Omega_f}\nabla p_0 \cdot N_j^v W_j$,
 we get the following macroscopic variational form
\begin{equation}
\label{eq:weak2}
\begin{split}
\int_{\Omega}
\alpha_{ij}V_iW_j
+
\int_{\Omega}
\gamma_{ij}\nabla V_i\,W_j
+
\int_{\Omega}
\gamma_{ji}V_i\,\nabla W_j
+\\
\int_{\Omega}
\theta_{ij}\nabla V_i\cdot \nabla W_j
+
\int_{\Omega}
\rho_j \nabla  p_0(x)W_j
=
\int_{\Omega}
\eta_jW_j
+
\int_{\Omega}
\zeta_j\cdot \nabla W_j.
\end{split}
\end{equation}
This equation in the strong form has the following form.
\begin{equation}
\label{eq:stokes1}
\alpha_{ij}V_i
-
\nabla\cdot(\theta_{ij}\nabla V_i)
+
\gamma_{ij}\nabla V_i
-
\nabla \cdot (\gamma_{ji} V_i)
+
\rho_j \nabla p_0
=
\eta_j-\nabla\cdot \zeta_j
\end{equation}

Because we use a single continuum for pressure, the incompressibility
equation can be derived by integrating $\div(v)=0$ over each coarse block. 
In the next section, we derive the incompressibility equation for multicontinuum.
If we take in \eqref{eq:prestest}, $M_j^p=0$ and $N_j^p=1$ ($j=1$ one continuum for pressure), we get instead of
equation \eqref{eq:div33}, the following equation
\begin{equation}
\label{eq:div1}
D_{ji}\cdot V_i + E_{ji}:\nabla V_i=0.
\end{equation}
In this (single) equation, the first term accounts the transfer term between the continua and the second term is incompressibility for the phases.  

 We note that if the effect
of $M$ is negligible (see our numerical results), the resulting model, 
\eqref{eq:stokes1} and \eqref{eq:div1}, becomes
a standard two-phase flow equations. In this case, we have (in the absence of 
external forces)
\begin{equation}
V_i=-\mathcal{K}(c) \nabla p_0
\end{equation}
 This is a standard two-phase flow
equations in the absence of gravity.

The mass conservation equation can be obtained 
from the 
concentration equation
using (see \cite{al2025dynamic}). To demonstrate 
this, we
note that in our numerical examples and, in practice, one chooses
the continua based on concentration. More precisely, we assume two
distinct values for $c$ which identifies each phase. 
We denote 
the continua region by $\psi_k$, where $\psi_k$ is the characteristic
function for the region, where the concentration is constant.
Multiplying the transport equation by 
$\psi_k(c)$, we obtain
\begin{equation}
\label{eq:multicontinuum_concentration_model}
\begin{split}
0=\int (c)_t\psi_k   + \int v \cdot \nabla c \psi_k =
\int (c\psi_k)_t -  \int c(\psi_k)_t   +\\
 \int_{\partial_{\psi_k^{ext}}} v\cdot n c -
 \int c v \cdot \nabla \psi_k=
(C_k)_t  + \\
\int_{\partial_{\psi_k^{ext}}} v\cdot n c - \int c(-v \cdot \nabla \psi_k) -
 \int c v \cdot \nabla \psi_k=
(C_k)_t  + \int_{\partial K} v\cdot n c\psi_k.
\end{split}
\end{equation}
Next, we use the following expansions
\begin{equation}
c= N_j^c C_j,\ \ 
v= N_j^v V_j + M_j^{v} \cdot \nabla V_j
\end{equation}
to compute 
\begin{equation}
\begin{split}
\int_{\partial K} v\cdot n c\psi_k=\int_{\partial K} (N_j^v V_j + M_j^{v} \cdot \nabla V_j ) \cdot n N_m^c C_m \psi_k \approx \\
\int_{\partial K}  \left\{ N_j^v N_k^c (V_j\cdot n) C_k + ( M_j^{v} \cdot n N_k^c) \nabla V_j C_k  \right\} =\\
\int_{K} \left\{ \div(\xi_{jk} V_j C_k) +\nabla_n (\sigma_{jk}^n \nabla V_j C_k) \right\},
\end{split}
\end{equation}
where $\xi_{jk}$ and $\sigma_{jk}^n$ are appropriately defined.
We can assume
$\sigma$ is negligible if convective effects are dominant. 
Note that if we take 
$N_k^c=\psi_k / \int \psi_k$, then $(C_k)_t + \div(\xi_{jk} V_j C_k^*)=0$,
where $C_k^* = C_k / \int \psi_k$ is a constant concentration value in continuum $k$.


\subsection{Multicontinuum pressure}

Next, we consider multicontinuum pressure case and 
associate 
a pressure field with each phase via multiscale basis
functions.
We expand the pressure as
\begin{equation}
p(x)=N_i^p(x)\,P_i(x)+M_i^p(x)\cdot \nabla P_i(x),
\nonumber
\end{equation}
where $N_i^p(x)$ - pressure basis functions,
$M_i^p(x)$ - pressure enrichment functions, and
$P_i(x)$-  coarse pressure unknowns.

First, we will derive the multicontinuum mass conservation.
We use the pressure test functions
\begin{equation}
\label{eq:prestest}
q(x)=N_j^p(x)\,Q_j(x)+M_j^p(x)\cdot \nabla Q_j(x)
\end{equation}
and impose
\begin{equation}
\int_{\Omega_f} q\,(\nabla\cdot v)=0
\end{equation}
with velocity expansion, as before
\begin{equation}
v(x)=N^v_i(x)V_i(x)+M^v_i(x)\cdot \nabla V_i(x).
\nonumber
\end{equation}
Thus,
\begin{equation}
\int_{\Omega_f}
\left(N_j^pQ_j+M_j^p\cdot \nabla Q_j\right)
\nabla\cdot
\left(N_i^vV_i+M_i^v\cdot \nabla V_i\right)
=0.
\nonumber
\end{equation}
We have (ignoring higher order terms)
\begin{equation}
\nabla\cdot v \approx 
(\nabla \cdot N_i^v) V_i + N_i^v: (\nabla V_i) + (\nabla\cdot M_i^v):\nabla V_i, 
\nonumber
\end{equation}
where $\overline{B}_i^v=N_i^v + \nabla\cdot M_i^v$. More precisely,
we write the index form (assuming Einstein notation)
\begin{equation}
\nabla_k v_k = 
(\nabla_k  N_i^{v,lk}) V_i^l + N_i^{v,lk}\nabla_k V_i^l + (\nabla_k M_i^{v,rlk})\nabla_r V_i^l,
\nonumber 
\end{equation}
where $\overline{B}_i^{v,rl}=N_i^{v,lk}\delta_{kr} + \nabla_k M_i^{v,rlk}$.
Here, the superindex refers to the vector and subindex refers
to the continuum.
The weak mass conservation equation  is 
\begin{align}
\sum_K \int_{K\cap\Omega_f}
\left(N_j^pQ_j+M_j^{p,n} \nabla_n Q_j\right)
\Big[
(\nabla_k  N_i^{v,lk})V_i^l
+
\overline{B}_i^{v,rl} \nabla_r V_i^l
\Big]
=0.
\nonumber
\end{align}

We have the following coupling tensors
\begin{equation}
\label{eq:macro1}
\begin{split}
D_{ji}^l
=
{1\over |K|}\int_{K\cap \Omega_f}
N_j^p\,(\nabla_k  N_i^{v,lk}),
\quad
E_{ji}^{kl}
=
{1\over |K|}\int_{K\cap \Omega_f}
N_j^p \overline{B}_i^{v,kl}\\
F_{ji}^{nl}
=
{1\over |K|}\int_{K\cap \Omega_f}
M_j^{p,n} (\nabla_k  N_i^{v,lk})\quad 
G_{ji}^{nkl}
=
{1\over |K|}\int_{K\cap \Omega_f}
M_j^{p,n} \overline{B}_i^{v,kl}.
\nonumber
\end{split}
\end{equation}
The  weak formulation of the mass conservation is 
\begin{align}
\int_{\Omega}
(
D_{ji}^lV_i^lQ_j + E_{ji}^{kl}\nabla_k V_i^l Q_j
+ F_{ji}^{nl}V_i^l\nabla_n Q_j
+ G_{ji}^{nkl}\nabla_k V_i^l\nabla_n Q_j
)
=0.
\nonumber
\end{align}
In the strong form,
\begin{equation}
\label{eq:div3}
D_{ji}^lV_i^l
+
E_{ji}^{kl}\nabla_k V_i^l
-
\nabla_n (F_{ji}^{nl} V_i^l)
-
\nabla_n(G_{ji}^{nkl}\nabla_k V_i^l)
=0.
\end{equation}
Or
\begin{equation}
\label{eq:div33}
D_{ji}\cdot V_i
+
E_{ji}:\nabla V_i
-
\nabla\cdot  (F_{ji} \cdot V_i)
-
\nabla\cdot (G_{ji}:\nabla V_i)
=0.
\end{equation}

Note that if the terms $\overline{B}_i^v$ are negligible, the mass conservation equation
becomes
\begin{equation}
\label{eq:div2}
D_{ji} \cdot V_i
-
\nabla \cdot (F_{ji}\cdot V_i)
\approx 0.
\end{equation}

Next, we briefly discuss the derivation of momentum equations. 
The main difference will be in the terms involving the pressure,
which will be the same as \eqref{eq:div3} and added to momentum equations.
In particular, assuming the notation
\begin{equation}
w=N_j^vW_j+M_j^v\cdot \nabla W_j,
\end{equation}
the extra term in the momentum equation will have the form 
\begin{equation}
-\int_{\Omega_f}
\left(N_i^pP_i+M_i^p\cdot \nabla P_i\right)
\nabla\cdot
\left(N_j^v W_j+M_j^v\cdot \nabla W_j\right).
\end{equation}
The resulting weak form of the equation can be obtained by adding the extra terms in \eqref{eq:weak2} and will be 
\begin{equation}
\begin{split}
\int_{\Omega}
\alpha_{ij}V_iW_j
+
\int_{\Omega}
\gamma_{ij}\nabla V_i\,W_j
+
\int_{\Omega}
\gamma_{ji}V_i\,\nabla W_j
+
\int_{\Omega}
\theta_{ij}\nabla V_i\cdot \nabla W_j
\\
-
\int_{\Omega}
(D_{ij}P_i W_j
+
E_{ij}P_i \nabla W_j
+
F_{ij}\nabla P_i W_j
+
G_{ij}\nabla P_i\cdot \nabla W_j)\\
=
\int_{\Omega}
\eta_jW_j
+
\int_{\Omega}
\zeta_j\cdot \nabla W_j.
\end{split}
\end{equation}

This equation in the strong form has the following form.
\begin{equation}
\label{eq:stokes2}
\begin{split}
\alpha_{ij}V_i
-
\nabla\cdot(\theta_{ij}\nabla V_i)
+
\gamma_{ij}\nabla V_i
-
\nabla \cdot (\gamma_{ji} V_i)
\\
-D_{ij}P_i
+
\nabla\cdot(E_{ij}P_i) 
-
F_{ij}\nabla P_i
+
\nabla\cdot(G_{ij}\nabla P_i)
=
\eta_j-\nabla\cdot \zeta_j.
\end{split}
\end{equation}

The equation \eqref{eq:stokes2} coupled with mass conservation equations
\eqref{eq:div2} represents macroscopic model for two-phase flow,
when coefficients depend on concentration. The macroscopic equations
for concentrations are the same as in the case of single pressure continuum.

We can use multicontinuum expansion for Cahn-Hilliard coupled with Stokes
also. We will present the results elsewhere. 

\subsection{Construction of multiscale basis functions}

To construct multiscale basis functions for the velocity, we first need to compute the auxiliary functions. Several choices are possible for their construction. In our numerical examples, the velocity auxiliary functions are obtained by solving constrained Stokes problems with average constraints on the coarse-grid edges. For the pressure field, no auxiliary functions are constructed, since we adopt the single-continuum pressure formulation and use standard piecewise constant basis functions on the coarse grid.

Let $\Omega = (0, L_1) \times (0, L_2)$ be a slab-type domain, and let $\Omega_f = \Omega \setminus \Omega_s$. We introduce a coarse grid $\mathcal{T}_H$ with $n_c^H \times 1$ coarse blocks. The coarse blocks are
defined as $K_r = (x_{r-1}, x_r) \times (0, L_2)$, $r=1, \dots, n_c^H$, with the size $H=L_1/n_c^H$.

We denote by \(\mathcal{E}_H\) the set of all vertical coarse-grid edges, including
both boundary edges and edges separating neighboring coarse blocks. For each edge \(E_l \in \mathcal{E}_H\), we then define the associated coarse neighborhood by
\[
    \omega_l
    =
    \bigcup
    \left\{
    K_j \in \mathcal{T}_H:
    E_l \subset \partial K_j
    \right\}.
\]
In addition, we introduce a general notation $\omega_l^{*}$, which may denote either an oversampled extension of $\omega_l$ or the global domain $\Omega$. To account for perforations, we define $\omega_l^f = \omega_l \cap \Omega_f$ and $\omega_l^{*f} = \omega_l^{*} \cap \Omega_f$.

For macroscopic model derivation, we kept simplified indices for $N$ and $M$
to show their smooth variations with respect to the spatial variables, i.e.,
the indices representing the spatial locations of macroscopic grid blocks
are  not used in macroscopic model derivations.
In this section, we will use indices for discrete spatial locations.
The auxiliary functions $N_i^{v,l}$ are constructed by solving the Stokes flow problem in $\omega_l^{*f}$ subject to coarse-edge average constraints and compatibility conditions. More precisely, $N_i^{v,l}$ consists of vector-valued functions $N_{i}^{v,ls}$, $s=1,2$, corresponding to the velocity directions. Each $N_{i}^{v,ls}$ can be written as $N_{i}^{v,ls} = (N_{i}^{v,ls1}, N_{i}^{v,ls2})$. Thus, we obtain the following problems
\begin{equation}\label{eq:two_phase_N_i}
\begin{split}
&\int_{\omega_l^{*f}} 2\mu\,\varepsilon(N_{i}^{v, ls}):\varepsilon(w) 
+ \int_{\Gamma_s \cap \overline{\omega}_l^{*f}} \beta (N_{i}^{v, ls})_\tau \cdot w_\tau
- \int_{\omega_l^{*f}} \varrho_{i}^{ls} \nabla \cdot w \\
&- \sum_{E_r \in \mathcal{E}_H (\omega_l^{*})} \sum_{j} \Lambda_{ij}^{slr} \int_{E_r \cap \overline{\omega}_l^{*f}} \psi_j w
= 0,\\
&\int_{\omega_l^{*f}} \div (N_{i}^{v, ls}) q - \sum_{K_d \subset \omega_l^{*}} \Theta_{i}^{lds} \int_{K_d \cap \omega_l^{*f}} q = 0,\\
&\int_{E_r \cap \overline{\omega}_l^{*f}} N_{i}^{v, ls} \psi_j = e_s \delta_{ij} \int_{E_r \cap \overline{\omega}_l^{*f}} \psi_j, \qquad \forall E_r \in \mathcal{E}_H(\omega_l^*),\\
&\sum_{E \in \mathcal{E}_H(K_d)} \int_{E \cap \overline{\omega}_l^{*f}} N_{i}^{v, ls} \cdot \tilde{n}^d - \Theta_{i}^{lds} |K_d \cap \omega_l^{*f}| = 0, \qquad \forall K_d \subset \omega_l^*,
\end{split}
\end{equation}
where $e_s$ is the $s$-th column of the identity matrix $I_2$, $\tilde{n}^{d}$ is the unit outward normal vector on $\partial K_d \cap \overline{\omega}_l^{*f}$, and $\Lambda_{ij}^{slr}$ and $\Theta_{i}^{lds}$ denote the Lagrange multipliers associated with the coarse-edge average constraints and compatibility conditions, respectively. Here, $\mathcal{E}_H(D)$ denotes the restriction of $\mathcal{E}_H$ to $D \subseteq \Omega$, i.e., $\mathcal{E}_H(D) = \{ E \in \mathcal{E}_H: E \subset \overline{D} \}$.

We can represent the multiscale basis functions as $\phi_{i}^{v,ls} = N_{i}^{v,ls}\phi_0^l + M_{i}^{v,ls}\cdot \nabla \phi_0^l$. Instead of explicitly computing $M_{i}^{v,ls}$, we directly determine $\phi_{i}^{v,ls}$ by solving
\begin{equation}\label{eq:two_phase_phi_i}
\begin{split}
&\int_{\omega_l^{f}} 2\mu\,\varepsilon(\phi_{i}^{v, ls}):\varepsilon(w) 
+ \int_{\Gamma_s \cap \overline{\omega}_l^{f}} \beta (\phi_{i}^{v, ls})_\tau \cdot w_\tau
- \int_{\omega_l^{f}} \tilde{\varrho}_{i}^{ls} \nabla \cdot w = 0,\\
&\int_{\omega_l^{f}} \div (\phi_{i}^{v, ls}) q = \sum_{K_d \subset \omega_l} g_{is}^{d} \int_{K_d \cap \omega_l^{f}} q,\\
&\phi_{i}^{v, ls}|_{E \cap \overline{\omega}_l^f} = (N_{i}^{v, ls} \phi_0^l)|_{E \cap \overline{\omega}_l^f}, \quad \forall E \in \mathcal{E}_H (\omega_l),
\end{split}
\end{equation}
where $g_{is}^{d}$ is determined by the compatibility condition $g_{is}^{d} = \frac{1}{|K_d \cap \omega_l^f|} \int_{E_l \cap \overline{\omega}_l^f} N_{i}^{v, ls} \cdot \tilde{n}^d$. We note that $\phi_{i}^{v,ls}$ can be expressed via $N_i^{ls}$ and $M_i^{ls}$
(see \eqref{eq:mcbasis}) in our case since 
our macroscopic domain is one dimensional slab. In this case, 
using global $N_{i}^{v, s}$, $M_{i}^{v,ls}$ is rescaling of $\phi_{i}^{v, ls}$ by $\nabla_1\phi_0^l$, i.e., $M_{i}^{v,ls}=(\phi_{i}^{v, ls}-N_{i}^{v, s}\phi_0^l)/\nabla_1\phi_0^l$,
where $M_{i}^{v,ls}$ has only one direction.

One can obtain an approximation of $\phi_{i}^{v,ls}$ by solving \eqref{eq:two_phase_N_i} in the coarse neighborhood $\omega_l^f$, with zero Dirichlet boundary conditions imposed on $\bigcup_{\widetilde{E} \in \mathcal{E}_H (\omega_l) \setminus \{E_l\}} (\widetilde{E} \cap \overline{\omega}_l^f)$. In this case, the average constraints are imposed only on $E_l$. We refer to this approximation as the local approach in our numerical results. 

\subsection{Numerical results}

To test the proposed MC-GMsFEM approach, we consider two-dimensional model problems with two fluids in perforated media. We define fluid phases by
\begin{equation}\label{eq:two_phase_continua_definition}
I_1 = \{c: 0 \leq c < 0.5\}, \qquad
I_2 = \{c: 0.5 \leq c \leq 1 \},
\end{equation}
and set the viscosities $\mu_1 = 100$ and $\mu_2 = 1$. Using $I_k$ and $\mu_k$, we determine $\psi_k(c)$ and $\mu(c)$ \eqref{eq:mu_and_psi}.

\begin{figure}[hbt!]
\centering
\includegraphics[width=0.48\textwidth]{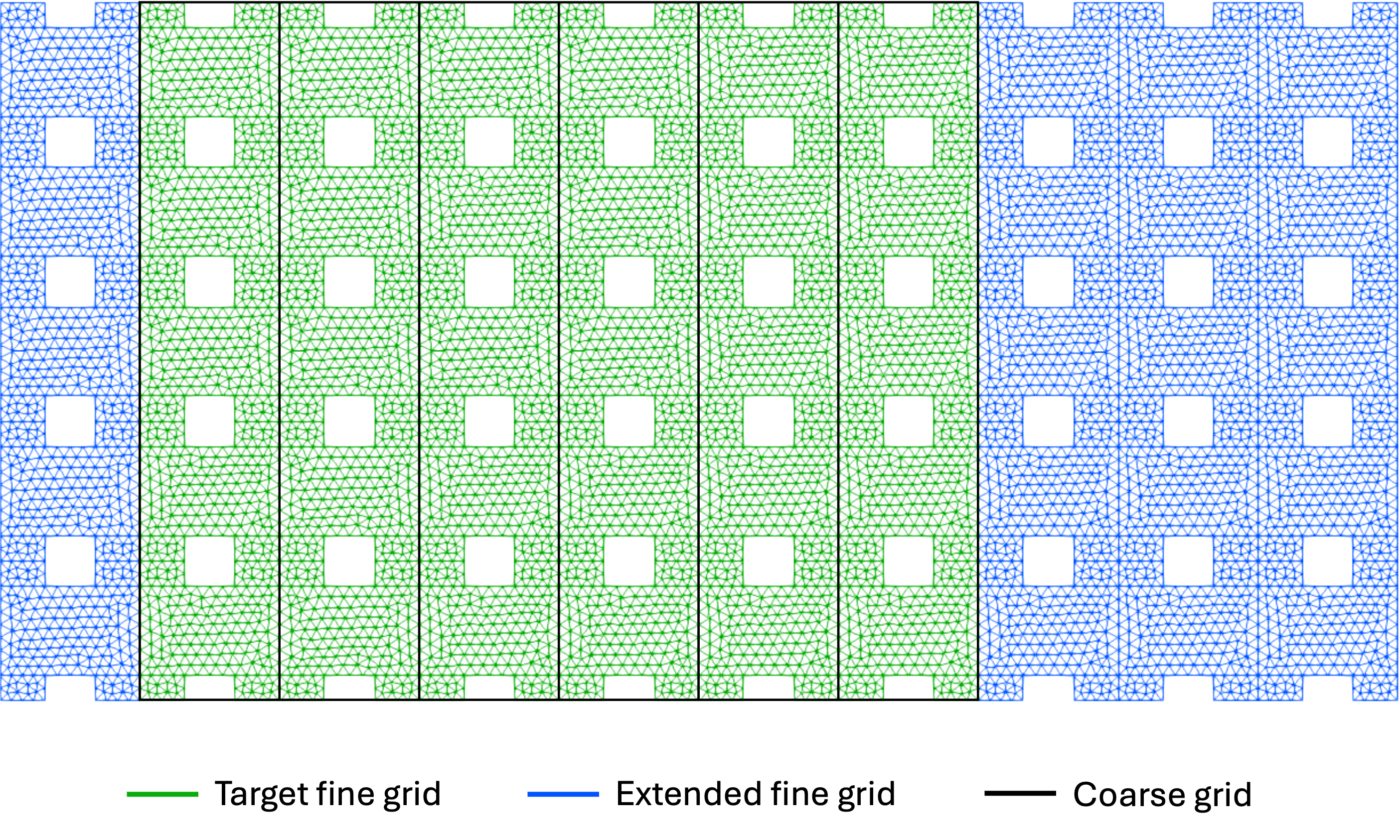}
\hfill
\includegraphics[width=0.48\textwidth]{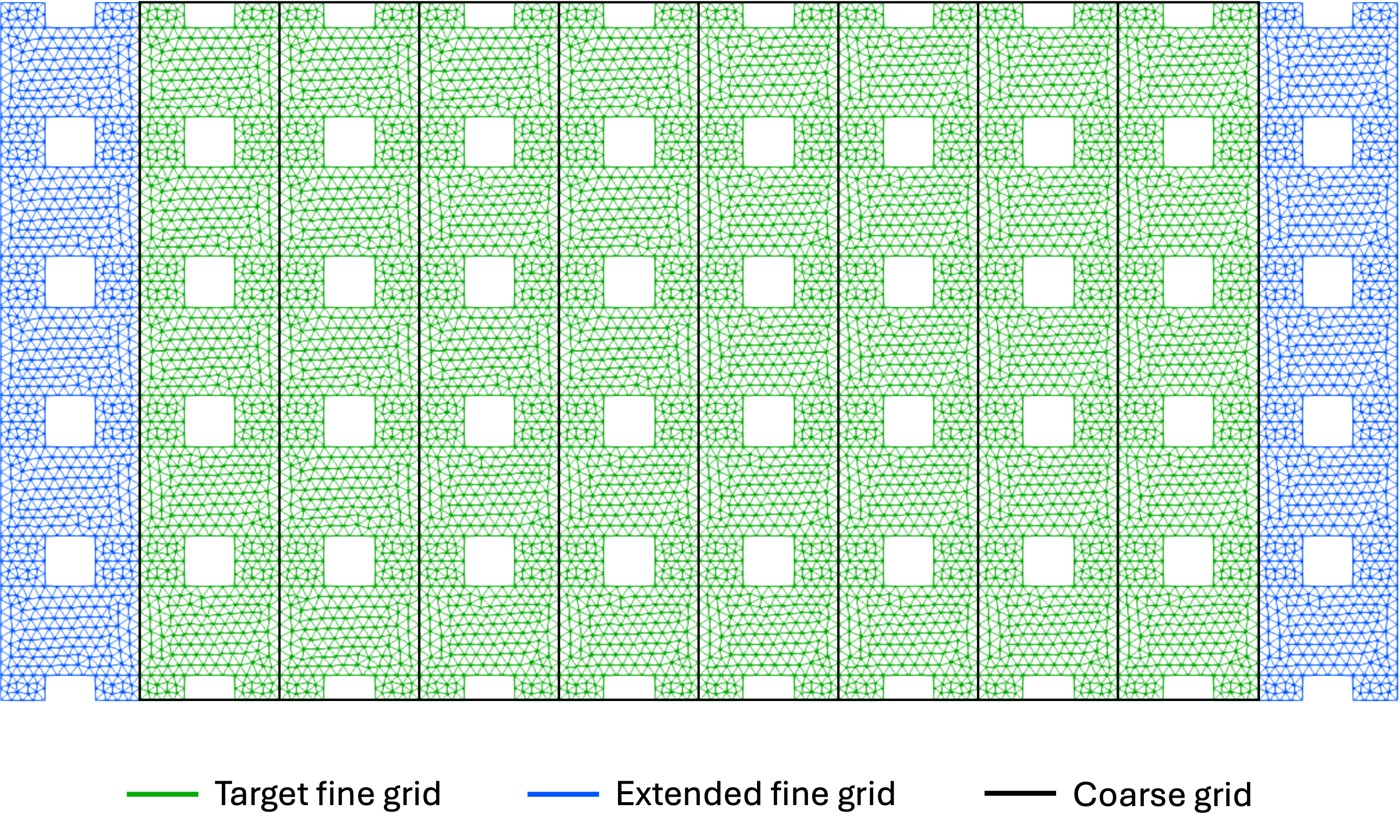}
\caption{Computational grids for Target domains 1 and 2 (from left to right).}
\label{fig:two_phase_grids}
\end{figure}

To generate boundary and initial data for the target model problems, we solve auxiliary problems in a domain larger than the target domain $\Omega$. More precisely, we introduce
\[
\Omega^{\mathrm{ext}} = (0, 1) \times (0, 0.5),
\qquad
\Omega_f^{\mathrm{ext}}
=
\Omega^{\mathrm{ext}}
\setminus
\Omega_s^{\mathrm{ext}},
\]
such that $\Omega \subset \Omega^{\mathrm{ext}}$ and
$\Omega_f \subset \Omega_f^{\mathrm{ext}}$. In particular, we consider two target domains $\Omega$ with their corresponding $\Omega_f$:
\begin{itemize}
\item Target domain 1: $\Omega = (0.1, 0.7) \times (0, 0.5)$;
\item Target domain 2: $\Omega = (0.1, 0.9) \times (0, 0.5)$.
\end{itemize}
See Figure \ref{fig:two_phase_grids} for an illustration.

The auxiliary problems solve \eqref{eq:two_phase_fine_model} for $(v_{\mathrm{ext}},p_{\mathrm{ext}},c_{\mathrm{ext}})$ in $\Omega_f^{\mathrm{ext}}$ with $f=(0, 0)$, the initial condition $c_{\mathrm{ext}}|_{t=0}=0.333$, and the boundary conditions
\begin{equation}\label{eq:two_phase_aux_bcs}
\begin{split}
v_{\mathrm{ext}} = (1,0), \quad c_{\mathrm{ext}} = 1
\quad &\text{on } \Gamma_{l}^{\mathrm{ext}}, \\ 
\sigma n = 0
\quad &\text{on } \Gamma_{r}^{\mathrm{ext}}, \\ 
v_{\mathrm{ext}}\cdot n = 0, \quad
\tau \cdot (\sigma n)
+
\beta (v_{\mathrm{ext}}\cdot \tau)
= 0
\quad &\text{on } \Gamma_{s}^{\mathrm{ext}},
\end{split}
\end{equation}
where $\Gamma_{l}^{\mathrm{ext}}$,
$\Gamma_{r}^{\mathrm{ext}}$, and
$\Gamma_{s}^{\mathrm{ext}}$ denote the left, right, and slip
boundaries of $\Omega_f^{\mathrm{ext}}$, respectively.

The friction coefficient $\beta(c)$ is defined as
\begin{equation}\label{eq:two_phase_beta}
\beta(c) = \frac{\mu(c)}{l_s},
\end{equation}
where $l_s = 0.1 \, d_p$ and $d_p = 0.036$ is the side length of the rectangular perforations.

The fine grid of $\Omega_f^{\mathrm{ext}}$ consists of 16174 triangular cells and 8583 vertices (see Figure \ref{fig:two_phase_grids}), generated using Gmsh \cite{geuzaine2009gmsh}. Note that the fine grids for $\Omega_f$ are subsets of that for $\Omega_f^{\mathrm{ext}}$. The Stokes problem is discretized using the finite element method with Taylor--Hood elements, implemented with the FEniCS package \cite{logg2012automated}. The Navier slip and no-penetration conditions are imposed weakly using Nitsche's method \cite{gjerde2022nitsche}.

The transport equation is solved using a Lagrangian--Eulerian particle-mesh method based on the LeoPart library \cite{maljaars2021leopart}.
In particular, we randomly seed 25 particles per cell. The particle data are projected onto a piecewise constant function space using an $l^2$ projection. To impose the inlet boundary condition, we fix the particle property to 1 in the cells adjacent to $\Gamma_{l}^{\mathrm{ext}}$.

We solve the Stokes flow and transport problems in a sequential manner. First, we solve the flow problem using the concentration from the previous time step. Then, we solve the transport problem using the obtained velocity. We employ the Forward Euler method for time discretization. The time step is chosen adaptively to keep $\text{CFL} = 0.5$, i.e., $\Delta t_k = \text{CFL} \, h_{\text{min}} / \max_{x \in \Omega_f^{\mathrm{ext}}} |v(x, t_k)|$, where $h_{\text{min}}$ is the minimum fine-grid cell diameter.

The duration of the simulations depends on the development of the fingers and consists of two stages. The pre-simulation stage lasts until the condition 
\begin{equation}\label{eq:two_phase_main_simulation_stage_criteria}
\int_{\Gamma_*} \psi_2(c_{\mathrm{ext}}^k) > 0, \quad \Gamma_* = \{ x \in \Omega_f^{\mathrm{ext}} : x_1 = 0.7 \}
\end{equation}
is satisfied, which indicates that the fingers of the second fluid have developed. To ensure that the fingers are sufficiently developed, we prolong the pre-simulation by ten additional time steps. After that, we perform the main simulation stage for 451 time steps, with the first time step used to define the initial condition. We test the proposed MC-GMsFEM for this main stage.

After solving the auxiliary problems, we use their solutions to define boundary conditions \eqref{eq:two_phase_fine_bcs} and initial condition \eqref{eq:two_phase_fine_ics} of the target problem
\begin{equation*}
\begin{gathered}
c_{\mathrm{in}} = c_{\mathrm{ext}}|_{\Gamma_{l}}, \quad v_{\mathrm{in}} = v_{\mathrm{ext}}|_{\Gamma_{l}}, \quad
v_{\mathrm{out}} = v_{\mathrm{ext}}|_{\Gamma_{r}},\\
c_0 = c_{\mathrm{ext}}(t_0)|_{\Omega_f},
\end{gathered}
\end{equation*}
where $t_0$ corresponds to the beginning of the main stage. We also set $f=(0, 0)$ for the target problem.


The considered coupled flow and transport problem is known to be unstable, and each run can produce different finger patterns. We therefore perform several independent runs of the auxiliary problem. For each run, the resulting auxiliary solution is used to define the boundary and initial data for the target problem through the restrictions described above. Thus, the target problem has the same general formulation in all tests, while the different test cases correspond to different finger realizations and, consequently, to different induced boundary and initial data. 

\begin{figure}[p]
\centering

\begin{subfigure}[b]{\textwidth}
\centering
\includegraphics[height=0.15\textheight]{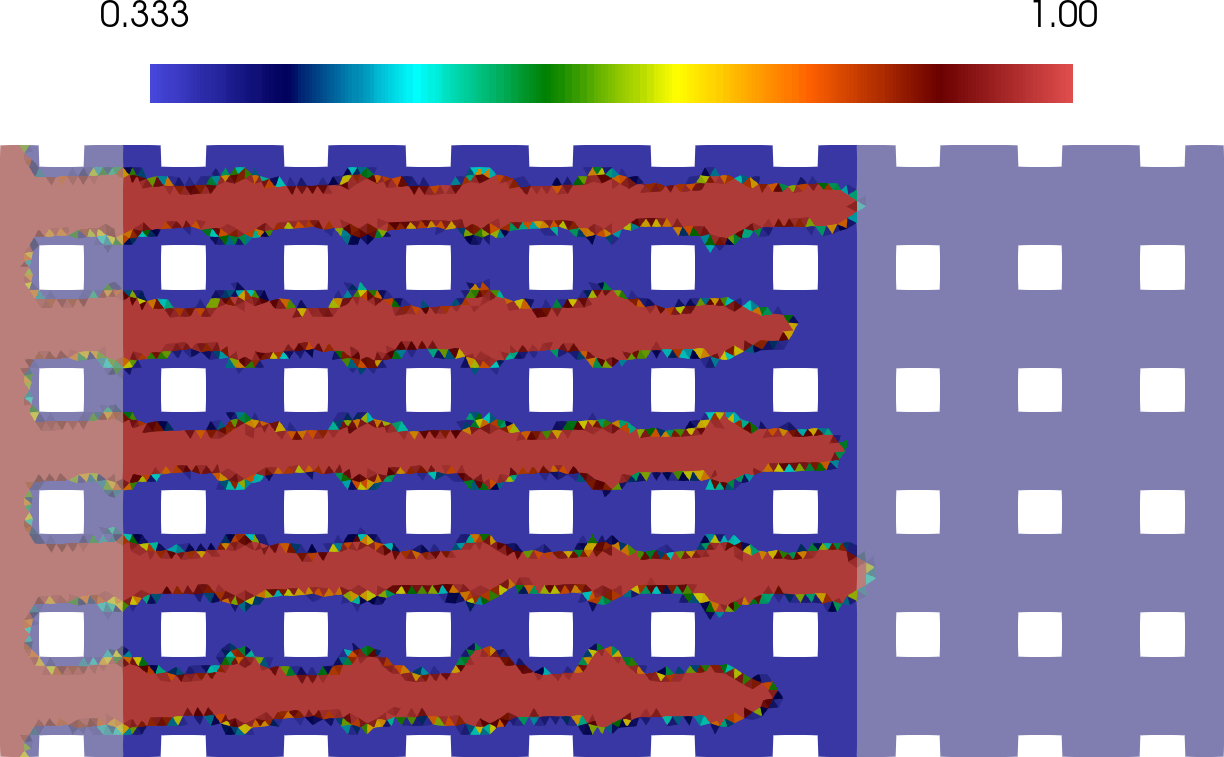} \hfill
\includegraphics[height=0.15\textheight]{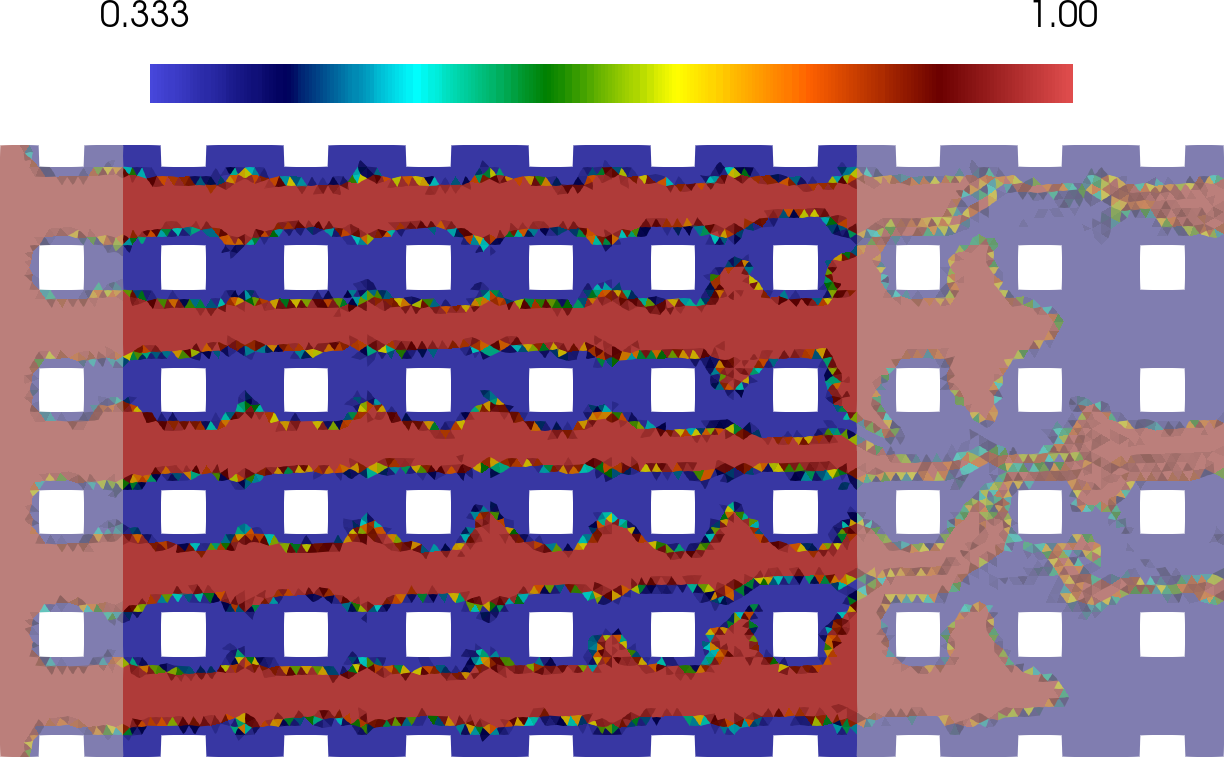}
\caption*{Test A}
\end{subfigure}

\begin{subfigure}[b]{\textwidth}
\centering
\includegraphics[height=0.15\textheight]{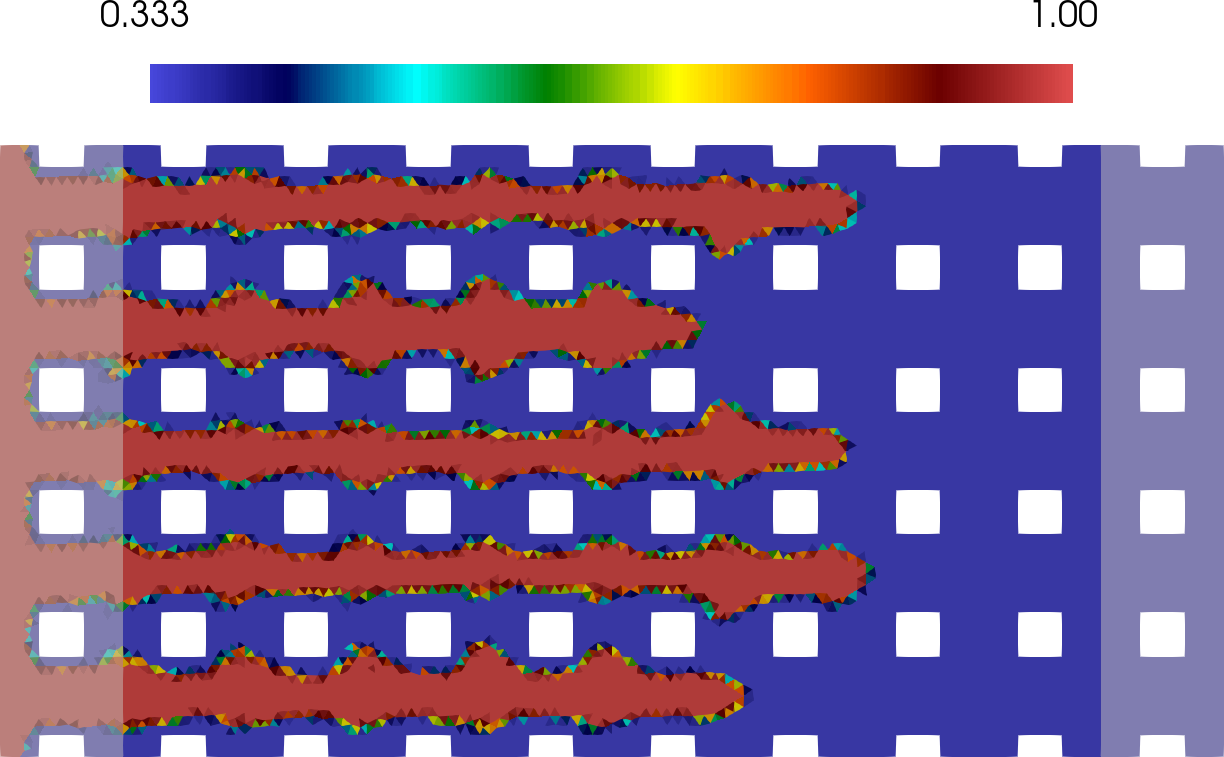} \hfill
\includegraphics[height=0.15\textheight]{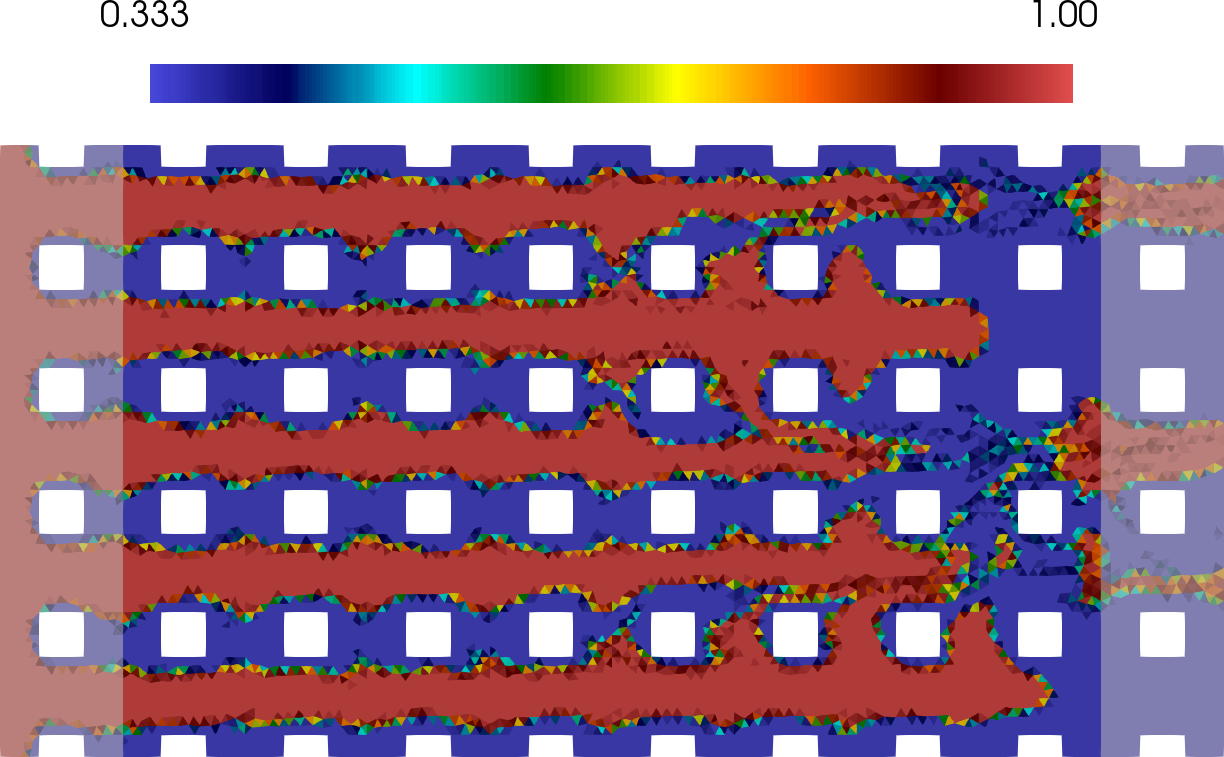}
\caption*{Test B}
\end{subfigure}

\begin{subfigure}[b]{\textwidth}
\centering
\includegraphics[height=0.15\textheight]{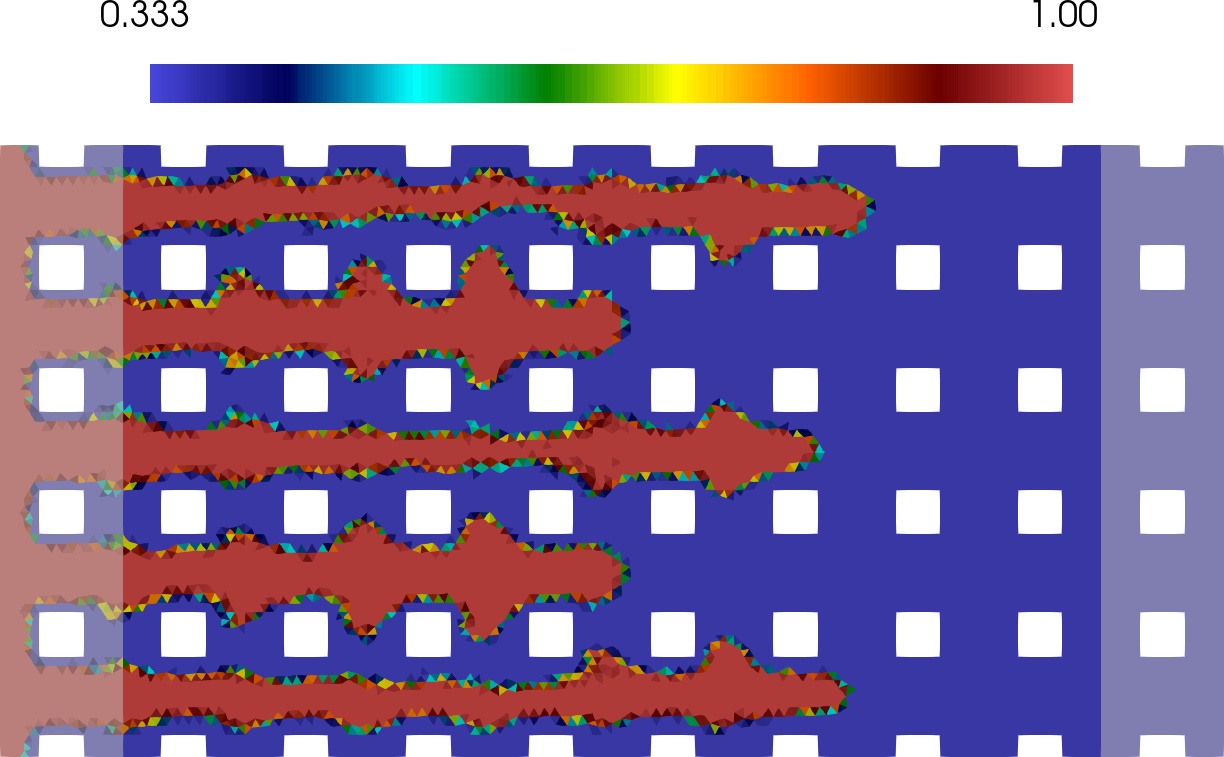} \hfill
\includegraphics[height=0.15\textheight]{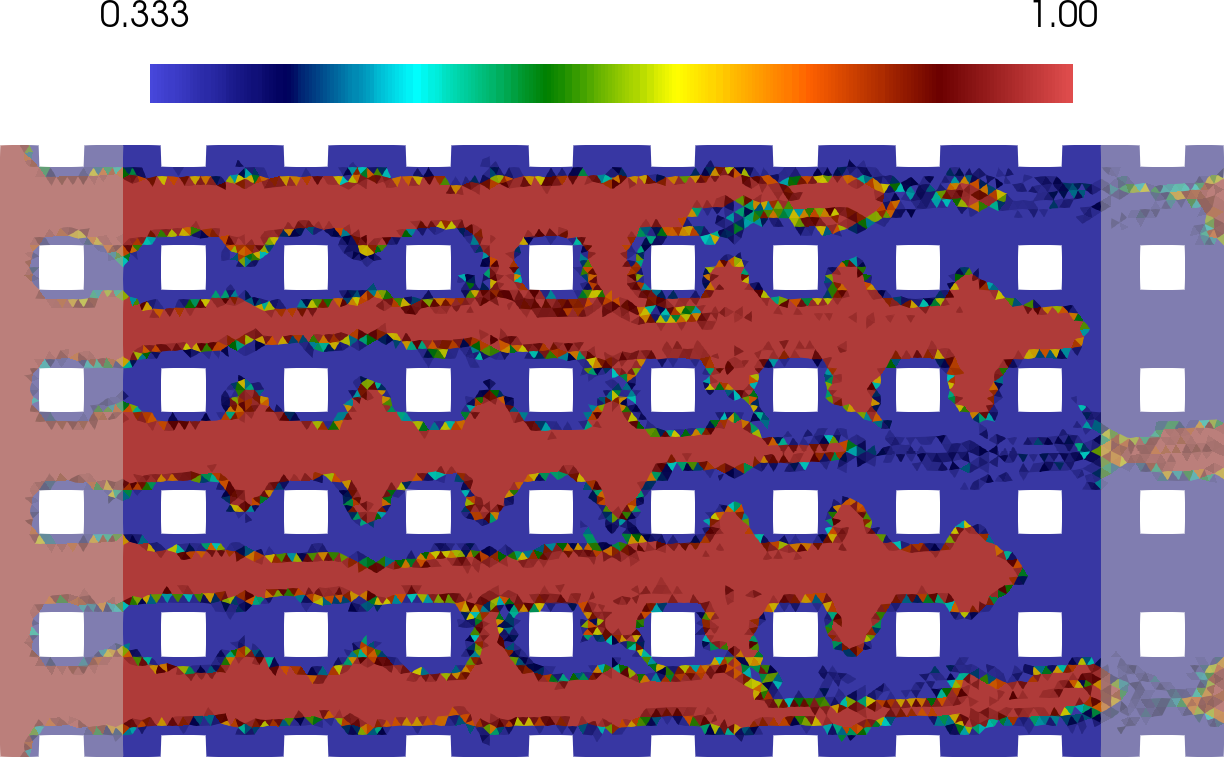}
\caption*{Test C}
\end{subfigure}

\begin{subfigure}[b]{\textwidth}
\centering
\includegraphics[height=0.15\textheight]{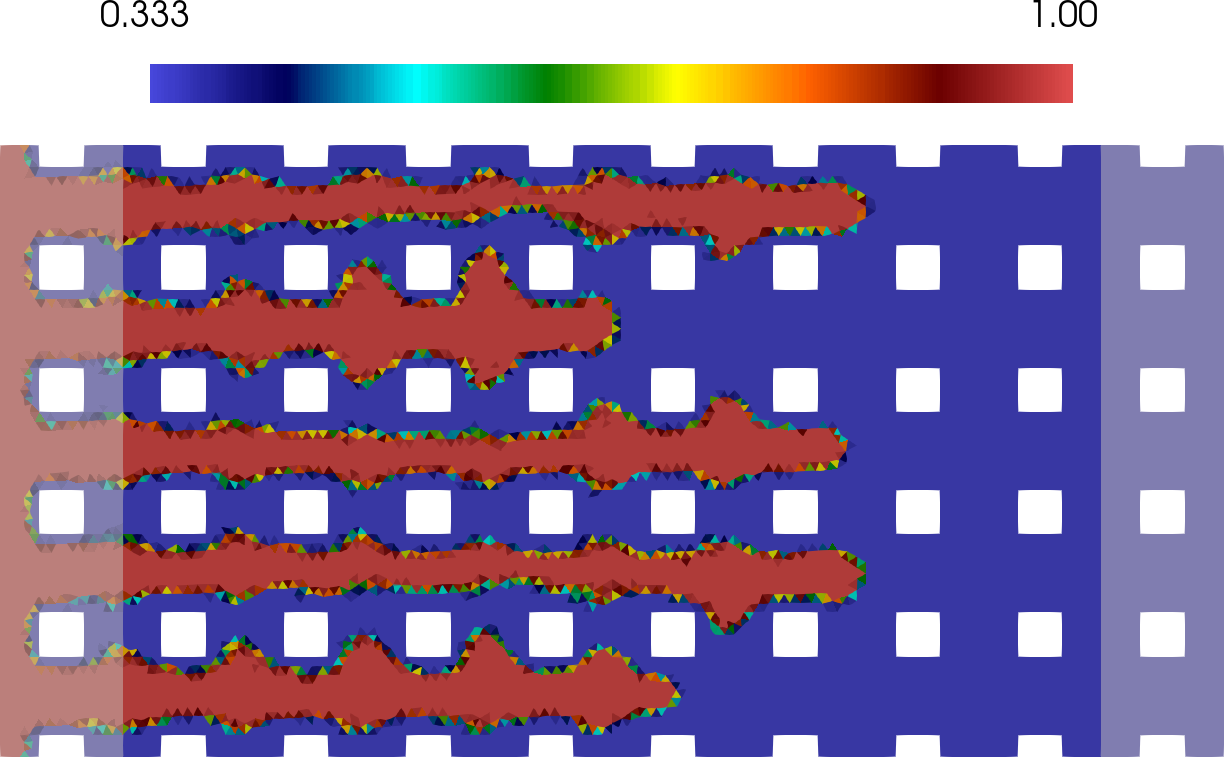} \hfill
\includegraphics[height=0.15\textheight]{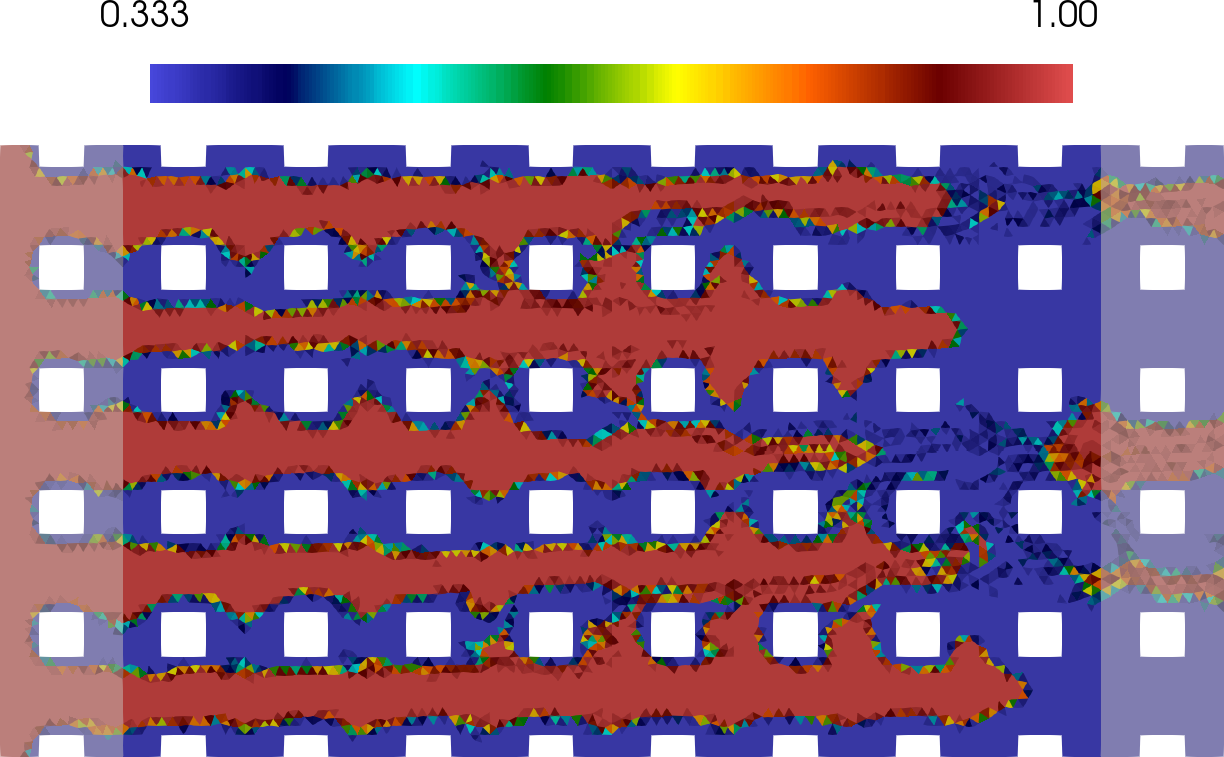}
\caption*{Test D}
\end{subfigure}

\begin{subfigure}[b]{\textwidth}
\centering
\includegraphics[height=0.15\textheight]{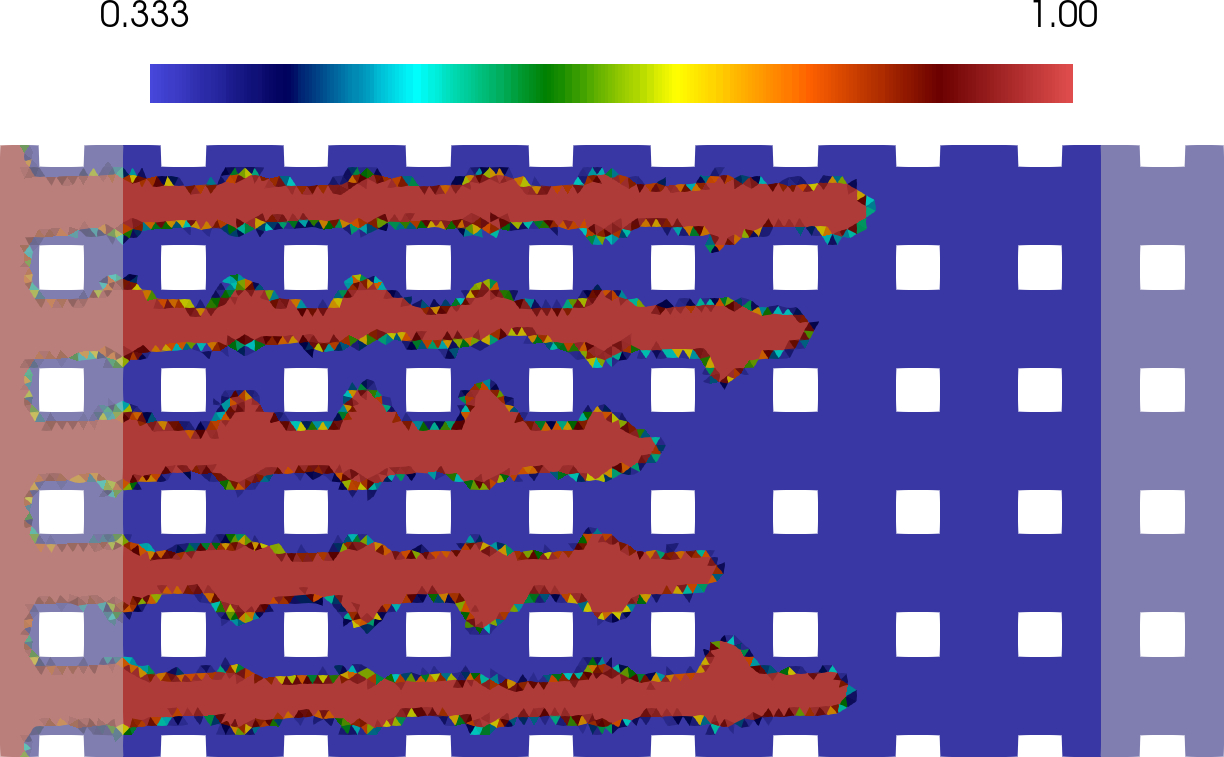} \hfill
\includegraphics[height=0.15\textheight]{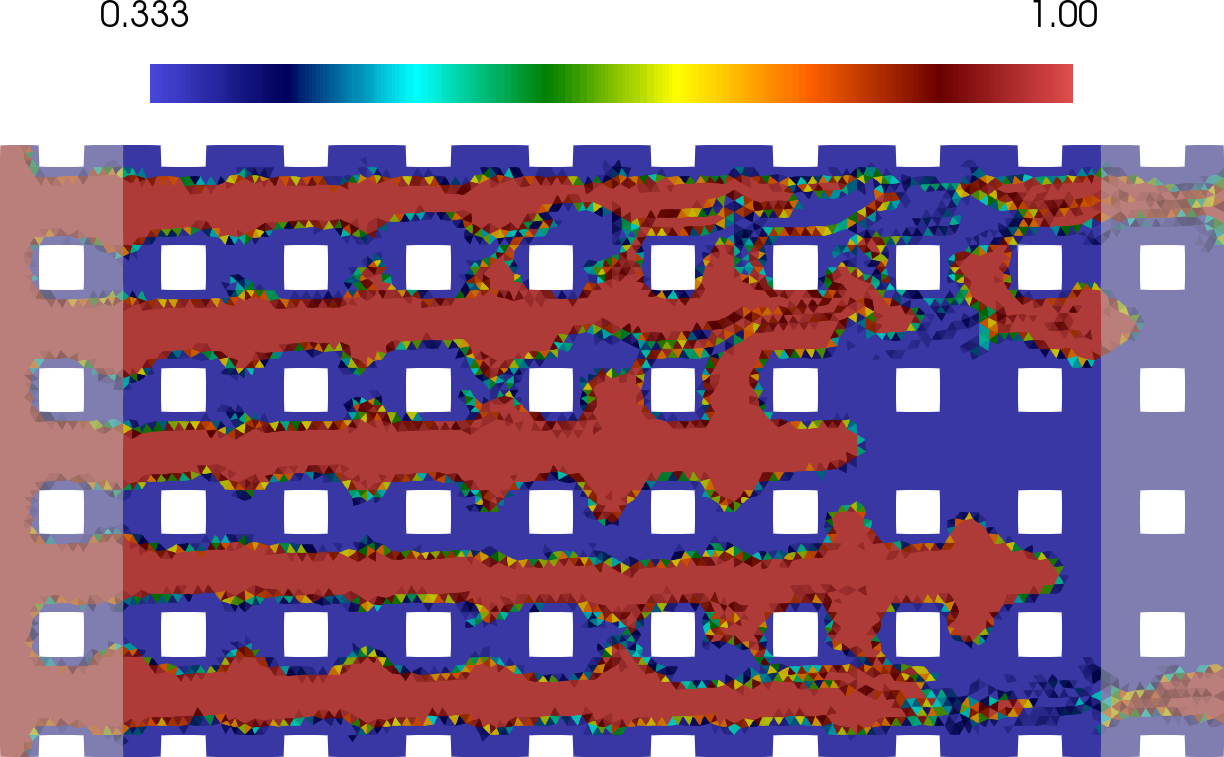}
\caption*{Test E}
\end{subfigure}

\caption{Initial and final fine-scale concentration distributions of the main simulation stage.}
\label{fig:two_phase_fine_results}
\end{figure}

In these numerical experiments, we consider both local and global MC-GMsFEM approaches. In the global approach, $N_i^{v}$ is computed by solving a global constraint problem, followed by the construction of the corresponding $\phi_i^{v,j}$. In the local approach, $\phi_i^{v,j}$ is approximated by solving a local problem with an average constraint on the target coarse edge. We tested both approaches for various finger realizations, and both provided good results. Throughout this section, we adopt the single-continuum pressure formulation of the macroscopic problem. In the following tests, we present representative numerical results to demonstrate the accuracy of the proposed method. Figure \ref{fig:two_phase_fine_results} presents the concentration distributions during the main simulation stage obtained from several runs. We use semi-transparent shading to indicate the regions outside the target domains. Thus, the unshaded regions correspond to the target domains considered in the corresponding tests.

From the concentration distributions, one can see that the more viscous fluid phase forms trapped regions, or ganglia. Such regions are often characterized by different properties than the connected regions of the same phase and may therefore exhibit different effective dynamics. One advantage of the proposed MC-GMsFEM is that it can account for spatially and temporally changing continua. This allows us to separate ganglia from the connected part of the more viscous phase and treat them as distinct continua. Therefore, in addition to the velocities of the two main continua/phases, the MC-GMsFEM provides separate velocities for the ganglia continua, capturing their motion independently from the connected phase.

We use a simple algorithm based on breadth-first search (BFS) to detect ganglia regions. The piecewise-constant indicator field $\psi_1$ is represented as a graph: grid cells are treated as graph nodes, and adjacency relations between cells define graph edges. We define ganglia as isolated connected regions of $I_1$ that do not touch the left and right boundaries and whose horizontal extent is more than two coarse-grid blocks. In addition, we connect two grid cells that share a common edge-neighbor cell; in a triangular mesh, such cells meet at a grid vertex. The detected ganglia regions are then removed from $\psi_1$ and assigned to new fields $\psi_k$, $k > 2$, representing ganglia continua. Pseudocode for the algorithm is provided in Appendix \ref{sec:ganglia_detection_algorithm}.

We emphasize that, our main goal is to derive accurate macroscopic two-phase flow models, and, thus, in the present numerical experiments, the indicator functions $\psi_k$, or equivalently the concentration-based decomposition into fluid phases and ganglia regions, are assumed to be available at each time step of the main simulation. This assumption is consistent with the general MC-GMsFEM framework introduced earlier, where the continua functions $\psi_k$ provide the structural information used to construct the auxiliary functions, multiscale basis functions, coarse-scale variables, and the corresponding macroscopic model. In the two-phase flow setting considered here, these continua are naturally determined by the concentration field and by the connected components of the corresponding phases. Thus, in our numerical tests, we compute $\psi_k = \psi_k(c_{\mathrm{ext}})|_{\Omega_f}$ from the fine-grid concentration field at each time step and use these functions to construct the multicontinuum basis functions and the corresponding macroscopic model. The present numerical experiments assess the accuracy of the proposed MC-GMsFEM coarse-scale approximation and the derived macroscopic equations for spatially and temporally varying multicontinuum configurations. The development of fully coarse-scale procedures will be studied in future work.

The coarse-scale simulations are performed in different target domains according to the test case. For Target domain 1, we use a uniform coarse grid consisting of $6 \times 1$ rectangular cells. For Target domain 2, we use a uniform coarse grid consisting of $8 \times 1$ rectangular cells. In both cases, the coarse-grid size (width of the coarse-grid blocks) is $H=0.1$. See Figure \ref{fig:two_phase_grids} for illustration. Time discretization is performed using the Forward Euler method with the same time steps as in the fine-grid simulations.

Note that in these numerical results, we mainly focus on the ability of MC-GMsFEM to upscale the two-phase flow problem under consideration and simultaneously derive the corresponding macroscopic equations. Therefore, we evaluate errors only for the upscaled solutions.

We consider several types of reference solutions. First, the reference velocities are obtained by averaging the fine-grid velocity fields, either with or without ganglia separation. This allows us to assess the effect of treating ganglia as distinct continua. Second, the reference concentrations are computed either by directly averaging the fine-grid concentration field or by solving the macroscopic transport model using the reference velocities. These choices lead to different comparison settings, which are specified in the corresponding tests. 

More precisely, we use the following reference solutions:
\begin{itemize}
\item $V_i^{\mathrm{ref,ng}}$ is the reference multicontinuum velocity without ganglia separation:
\begin{equation*}
V_i^{\mathrm{ref,ng}}|_{E} = \frac{\int_E v_1 \psi_i}{\int_E \psi_i},
\end{equation*}
where $E$ is a coarse-grid edge and $\psi_i$ are the regular indicator functions without ganglia separation. Since our macroscopic model is one-dimensional, we compare velocities in $x_1$ direction.

\item $V_i^{\mathrm{ref,wg}}$ is the reference multicontinuum velocity with ganglia separation:
\begin{equation*}
V_i^{\mathrm{ref,wg}}|_{E} = \frac{\int_E v_1 \hat{\psi}_i}{\int_E \hat{\psi}_i},
\end{equation*}
where $\hat{\psi}_i$ are indicator functions in which ganglia are treated as separate continua.

\item $C_i^{\mathrm{ref}}$ is the reference multicontinuum concentration:
\begin{equation*}
C_i^{\mathrm{ref}}|_{K} = \int_K c \psi_i,
\end{equation*}
where $K$ is a coarse-grid block.

\item $C_i^{\mathrm{mc}}(V^{\mathrm{ref,ng}})$ is the multicontinuum concentration obtained by solving the macroscopic transport model using $V_i^{\mathrm{ref,ng}}$.

\item $C_i^{\mathrm{mc}}(V^{\mathrm{ref,wg}})$ is the multicontinuum concentration obtained by solving the macroscopic transport model using $V_i^{\mathrm{ref,wg}}$.
\end{itemize}

In addition, to track the dynamics of the multicontinuum concentrations, we define
\begin{equation*}
C_{\mathrm{avg}, i}^{\mathrm{ref}}(t_k) 
= 
\frac{1}{|\Omega_{r}|} 
\int_{\Omega_{r}} C_i^{\mathrm{ref}}(t_k),
\end{equation*}
where $t_k$ denotes the time level and $\Omega_r$ is a subdomain of $\Omega$ that differs for each target domain. In particular, $\Omega_r = (0.3, 0.7) \times (0, 0.5)$ for Target domain 1, while $\Omega_r = (0.4, 0.9) \times (0, 0.5)$ for Target domain 2. We define $C_{\mathrm{avg}, i}^{\mathrm{mc}}(V^{\mathrm{ref,ng}}; t_k)$ and $C_{\mathrm{avg}, i}^{\mathrm{mc}}(V^{\mathrm{ref,wg}}; t_k)$, corresponding to $C_i^{\mathrm{mc}}(V^{\mathrm{ref,ng}})$ and $C_i^{\mathrm{mc}}(V^{\mathrm{ref,wg}})$, respectively, in the same manner.

The coarse-grid multicontinuum solutions are denoted by $V_i^{\mathrm{mc}}$ and $C_i^{\mathrm{mc}}(V^{\mathrm{mc}})$. Note that we use the indicator functions with ganglia separation, $\hat{\psi}_i$, to compute $V_i^{\mathrm{mc}}$. For $C_i^{\mathrm{mc}}(V^{\mathrm{mc}})$, however, we employ the regular indicator functions $\psi_i$ and use the velocities $V_1^{\mathrm{mc}}$ and $V_2^{\mathrm{mc}}$ corresponding to the two main continua/phases.

We use the following norms when comparing the solutions
\begin{equation*}
\begin{gathered}
\|V_i\|_2 = \sqrt{\sum_E (V_i|_E)^2}, \quad
\|C_i\|_2 = \sqrt{\sum_K (C_i|_K)^2}, \\
\|C_{\mathrm{avg}, i}\|_2 = \sqrt{\sum_{t_k} (C_{\mathrm{avg}, i}(t_k))^2},
\end{gathered}
\end{equation*}
where $V_i$, $C_i$, and $C_{\mathrm{avg},i}$ denote the quantities being compared, either reference or coarse-grid multicontinuum solutions.

\FloatBarrier
\subsubsection{Test A}

In this test, we consider Target domain 1, $\Omega = (0.1, 0.7) \times (0, 0.5)$. The fine-grid concentrations are presented in the first row of Figure \ref{fig:two_phase_fine_results}. Due to the size of the target domain and the criterion for starting the main simulation stage \eqref{eq:two_phase_main_simulation_stage_criteria}, the two fluid continua are already present on each coarse edge at the initial condition. Nevertheless, the continua evolve over time, and ganglia arise with spatially varying distributions. We use the global MC-GMsFEM approach.

\begin{figure}[hbt!]
\centering
\includegraphics[width=0.48\textwidth]{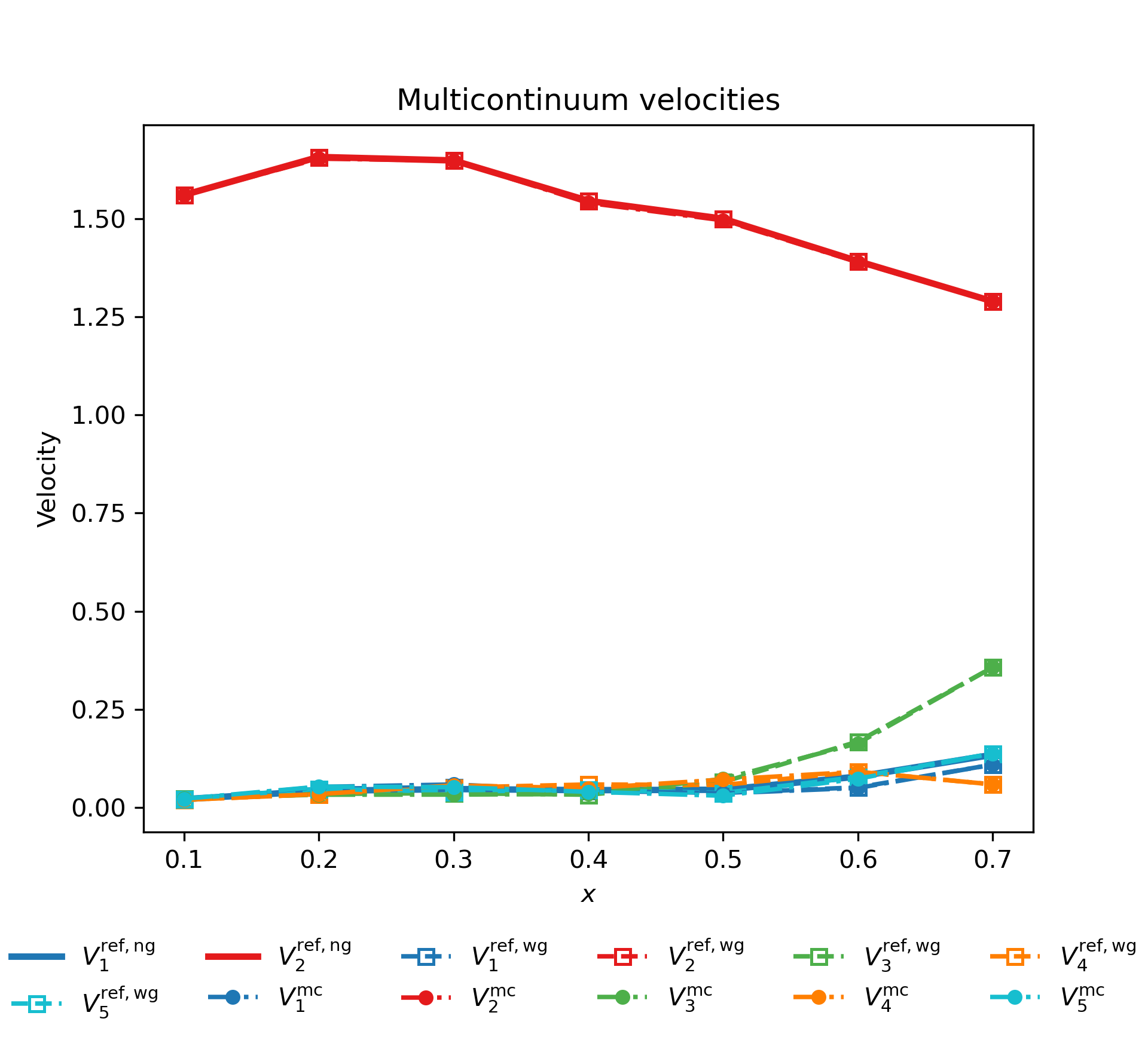}
\includegraphics[width=0.48\textwidth]{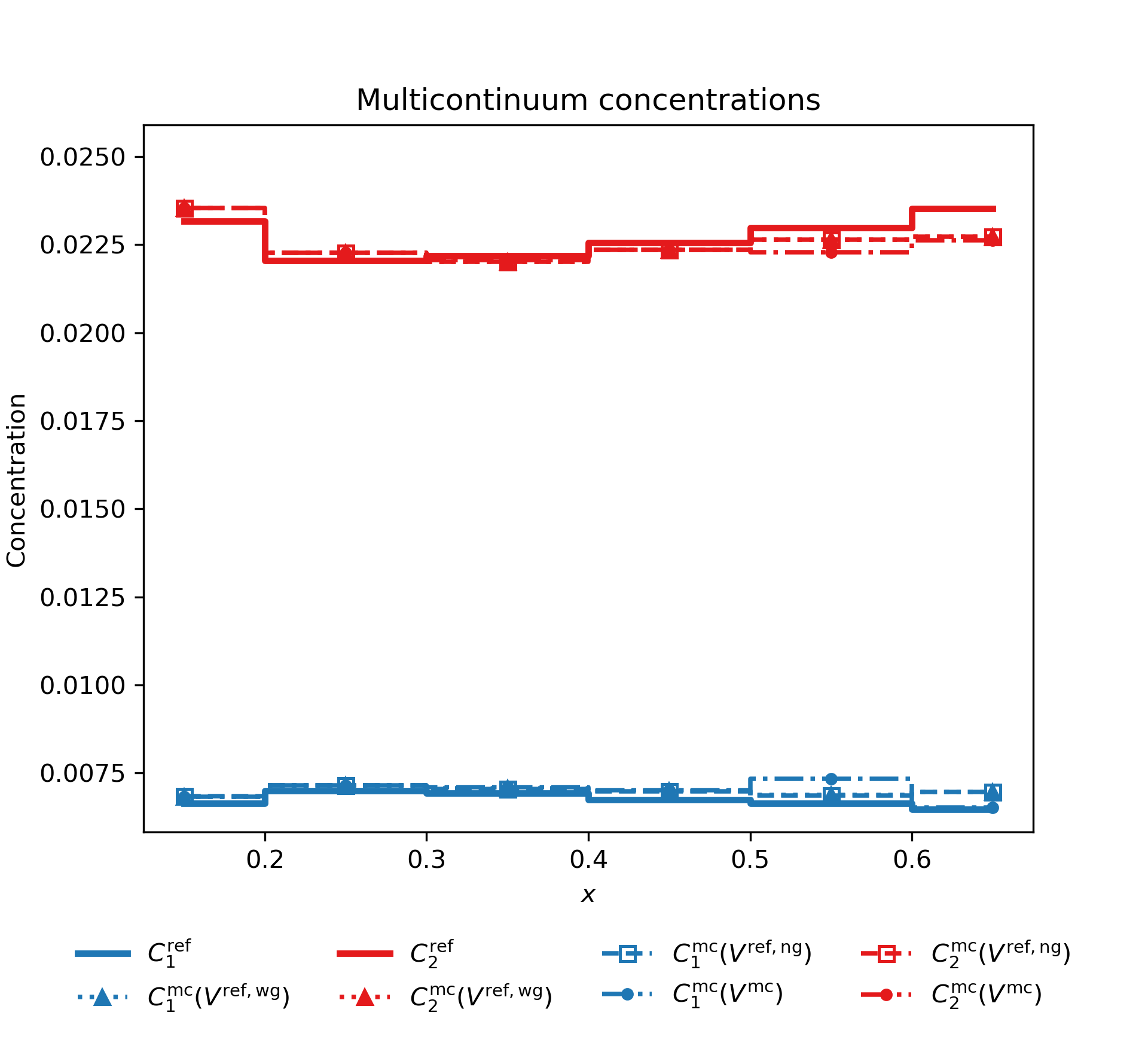}
\includegraphics[width=0.48\textwidth]{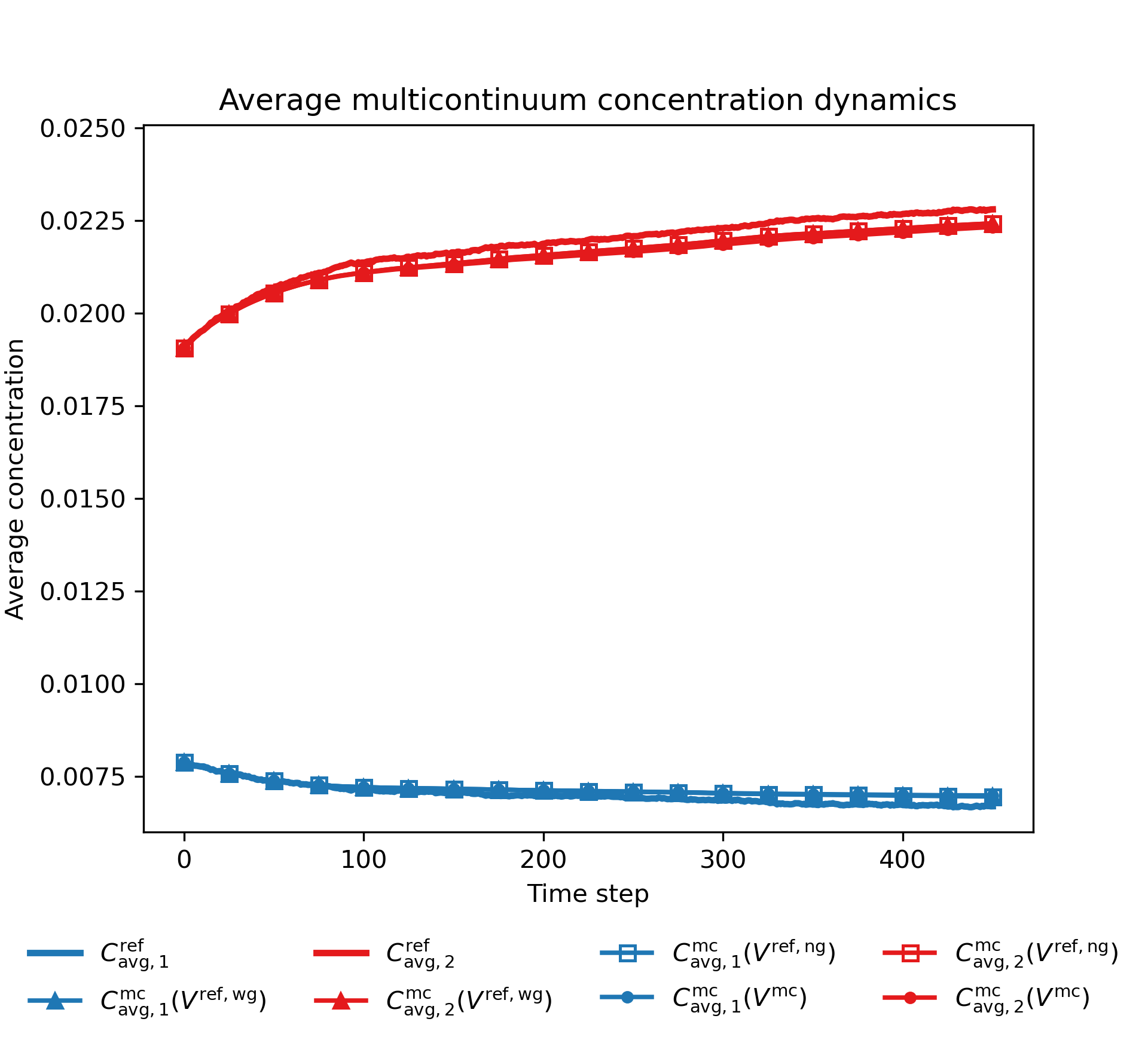}
\caption{Multicontinuum velocities, concentrations, and average concentration dynamics for Test A.}
\label{fig:two_phase_test_a_plots}
\end{figure}

Figure \ref{fig:two_phase_test_a_plots} presents the multicontinuum velocities and concentrations at the final time, as well as the evolution of the average multicontinuum concentrations over time. One can see that the corresponding curves are close to each other, indicating the high accuracy of the proposed MC-GMsFEM approach. The second-continuum velocity is significantly higher than those of the first continuum and the ganglia continua. Similarly, the concentration values of the second continuum are larger due to the definition of the continua \eqref{eq:two_phase_continua_definition}. For the average multicontinuum concentration dynamics, the second continuum slowly grows over time, while the first one decreases. This behavior is expected, since we track the right half of the domain.

Table \ref{tabs:two_phase_test_a_V_errors} presents the relative errors of the multicontinuum velocities at the final time. The top table shows the errors for the velocities of the two main continua, while the bottom table shows the errors for the ganglia velocities. In the first row of the top table, we compute the errors between the reference solutions with and without ganglia separation. We see that the errors are reasonable. The second-continuum velocity error is zero because the ganglia are excluded only from the first continuum. The second row computes the errors between the multicontinuum velocities and the reference velocities without ganglia separation. The obtained errors mostly correspond to those in the first row. The third row shows the errors between the multicontinuum velocities and the reference velocities with ganglia separation. This type of error represents the MC-GMsFEM approximation error without the influence of ganglia treatment. We can see that the obtained errors are small. From the bottom table, we see that the ganglia velocities are also accurately captured.
\begin{table}[hbt!]
\caption{Relative errors of multicontinuum velocities at the final time for Test A.}
\label{tabs:two_phase_test_a_V_errors}
\centering
\begin{tabular}{lll}
\hline
 Velocity error $e_V^{(i)}$                                                                                   & $e_V^{(1)}$   & $e_V^{(2)}$   \\
\hline
 $\|V_i^{\mathrm{ref, wg}} - V_i^{\mathrm{ref, ng}}\|_2 / \|V_i^{\mathrm{ref, ng}}\|_2 \times 100\%$ & 21.49 \%      & 0.00 \%       \\
 $\|V_i^{\mathrm{mc}} - V_i^{\mathrm{ref, ng}}\|_2 / \|V_i^{\mathrm{ref, ng}}\|_2 \times 100\%$      & 23.45 \%      & 0.21 \%       \\ \addlinespace
 $\|V_i^{\mathrm{mc}} - V_i^{\mathrm{ref, wg}}\|_2 / \|V_i^{\mathrm{ref, wg}}\|_2 \times 100\%$      & 5.42 \%       & 0.21 \%       \\
\hline
\end{tabular}

\vspace{0.4em}

\begin{tabular}{llll}
\hline
 Velocity error for ganglia $e_V^{(i)}$                                                                                   & $e_V^{(3)}$   & $e_V^{(4)}$   & $e_V^{(5)}$   \\
\hline
 $\|V_i^{\mathrm{mc}} - V_i^{\mathrm{ref, wg}}\|_2 / \|V_i^{\mathrm{ref, wg}}\|_2 \times 100\%$      & 2.52 \%       & 11.85 \%      & 8.86 \%       \\
\hline
\end{tabular}
\end{table}

In Table \ref{tabs:two_phase_test_a_C_errors}, we present the multicontinuum concentration errors at the final time. The first three rows compare three types of multicontinuum concentrations with the reference multicontinuum concentrations: those obtained using the reference velocities without ganglia separation, those obtained using the reference velocities with ganglia separation, and those obtained using the multicontinuum velocities. The next two rows compare the concentrations obtained using the reference velocities with ganglia separation and the multicontinuum velocities with those obtained using the reference velocities without ganglia separation. Finally, the last row compares the concentrations obtained using the multicontinuum velocities with those obtained using the reference velocities with ganglia separation. This error represents the MC-GMsFEM approximation error without the influence of numerical diffusion and ganglia separation. Overall, all errors are small. The effect of ganglia separation on the concentration errors is minor, and the second-continuum errors are smaller than those of the first continuum.
\begin{table}[hbt!]
\caption{Relative errors of multicontinuum concentrations at the final time for Test A.}
\label{tabs:two_phase_test_a_C_errors}
\centering
\begin{tabular}{lll}
\hline
 Concentration error $e_C^{(i)}$                                                                                                                                      & $e_C^{(1)}$   & $e_C^{(2)}$   \\
\hline
 $\|C_i^{\mathrm{mc}}(V^{\mathrm{ref, ng}}) - C_i^{\mathrm{ref}}\|_2 / \|C_i^{\mathrm{ref}}\|_2 \times 100\%$                                           & 3.90 \%       & 1.81 \%       \\
 $\|C_i^{\mathrm{mc}}(V^{\mathrm{ref, wg}}) - C_i^{\mathrm{ref}}\|_2 / \|C_i^{\mathrm{ref}}\|_2 \times 100\%$                                           & 3.97 \%       & 1.81 \%       \\
 $\|C_i^{\mathrm{mc}}(V^{\mathrm{mc}}) - C_i^{\mathrm{ref}}\|_2 / \|C_i^{\mathrm{ref}}\|_2 \times 100\%$                                                & 4.80 \%       & 2.22 \%       \\ \addlinespace
 $\|C_i^{\mathrm{mc}}(V^{\mathrm{ref, wg}}) - C_i^{\mathrm{mc}}(V^{\mathrm{ref, ng}})\|_2 / \|C_i^{\mathrm{mc}}(V^{\mathrm{ref, ng}})\|_2 \times 100\%$ & 0.19 \%       & 0.00 \%       \\
 $\|C_i^{\mathrm{mc}}(V^{\mathrm{mc}}) - C_i^{\mathrm{mc}}(V^{\mathrm{ref, ng}})\|_2 / \|C_i^{\mathrm{mc}}(V^{\mathrm{ref, ng}})\|_2 \times 100\%$      & 3.80 \%       & 0.67 \%       \\ \addlinespace
 $\|C_i^{\mathrm{mc}}(V^{\mathrm{mc}}) - C_i^{\mathrm{mc}}(V^{\mathrm{ref, wg}})\|_2 / \|C_i^{\mathrm{mc}}(V^{\mathrm{ref, wg}})\|_2 \times 100\%$      & 3.68 \%       & 0.67 \%       \\
\hline
\end{tabular}
\end{table}

Table \ref{tabs:two_phase_test_a_C_tracer_errors} presents the relative errors of the average multicontinuum concentration dynamics. It has the same comparison structure as Table \ref{tabs:two_phase_test_a_C_errors}. One can see that all the errors are minor. In the last row, which represents the MC-GMsFEM approximation error, all the errors are less than 1\%.
\begin{table}[hbt!]
\caption{Relative errors of average multicontinuum concentration dynamics for Test A.}
\label{tabs:two_phase_test_a_C_tracer_errors}
\centering
\begin{tabular}{lll}
\hline
 Average concentration dynamics error $e_{C, \mathrm{avg}}^{(i)}$                                                                                                                                                                      & $e_{C, \mathrm{avg}}^{(1)}$   & $e_{C, \mathrm{avg}}^{(2)}$   \\
\hline
 $\|C_{\mathrm{avg}, i}^{\mathrm{mc}}(V^{\mathrm{ref, ng}}) - C_{\mathrm{avg}, i}^{\mathrm{ref}}\|_2 / \|C_{\mathrm{avg}, i}^{\mathrm{ref}}\|_2 \times 100\%$                                           & 2.45 \%                       & 1.42 \%                       \\
 $\|C_{\mathrm{avg}, i}^{\mathrm{mc}}(V^{\mathrm{ref, wg}}) - C_{\mathrm{avg}, i}^{\mathrm{ref}}\|_2 / \|C_{\mathrm{avg}, i}^{\mathrm{ref}}\|_2 \times 100\%$                                           & 2.49 \%                       & 1.42 \%                       \\
 $\|C_{\mathrm{avg}, i}^{\mathrm{mc}}(V^{\mathrm{mc}}) - C_{\mathrm{avg}, i}^{\mathrm{ref}}\|_2 / \|C_{\mathrm{avg}, i}^{\mathrm{ref}}\|_2 \times 100\%$                                                & 2.59 \%                       & 1.77 \%                       \\ \addlinespace
 $\|C_{\mathrm{avg}, i}^{\mathrm{mc}}(V^{\mathrm{ref, wg}}) - C_{\mathrm{avg}, i}^{\mathrm{mc}}(V^{\mathrm{ref, ng}})\|_2 / \|C_{\mathrm{avg}, i}^{\mathrm{mc}}(V^{\mathrm{ref, ng}})\|_2 \times 100\%$ & 0.04 \%                       & 0.00 \%                       \\
 $\|C_{\mathrm{avg}, i}^{\mathrm{mc}}(V^{\mathrm{mc}}) - C_{\mathrm{avg}, i}^{\mathrm{mc}}(V^{\mathrm{ref, ng}})\|_2 / \|C_{\mathrm{avg}, i}^{\mathrm{mc}}(V^{\mathrm{ref, ng}})\|_2 \times 100\%$      & 0.16 \%                       & 0.37 \%                       \\ \addlinespace
 $\|C_{\mathrm{avg}, i}^{\mathrm{mc}}(V^{\mathrm{mc}}) - C_{\mathrm{avg}, i}^{\mathrm{mc}}(V^{\mathrm{ref, wg}})\|_2 / \|C_{\mathrm{avg}, i}^{\mathrm{mc}}(V^{\mathrm{ref, wg}})\|_2 \times 100\%$      & 0.12 \%                       & 0.37 \%                       \\
\hline
\end{tabular}
\end{table}

\FloatBarrier
\subsubsection{Test B}

This test and all subsequent tests present numerical results for Target domain 2, $\Omega = (0.1, 0.9) \times (0, 0.5)$. The fine-grid concentrations are depicted in the second row of Figure \ref{fig:two_phase_fine_results}. This test represents a more complicated case, where the second continuum is not present on all coarse edges at the initial condition and propagates over time. We use the global MC-GMsFEM approach.

\begin{figure}[hbt!]
\centering
\includegraphics[width=0.48\textwidth]{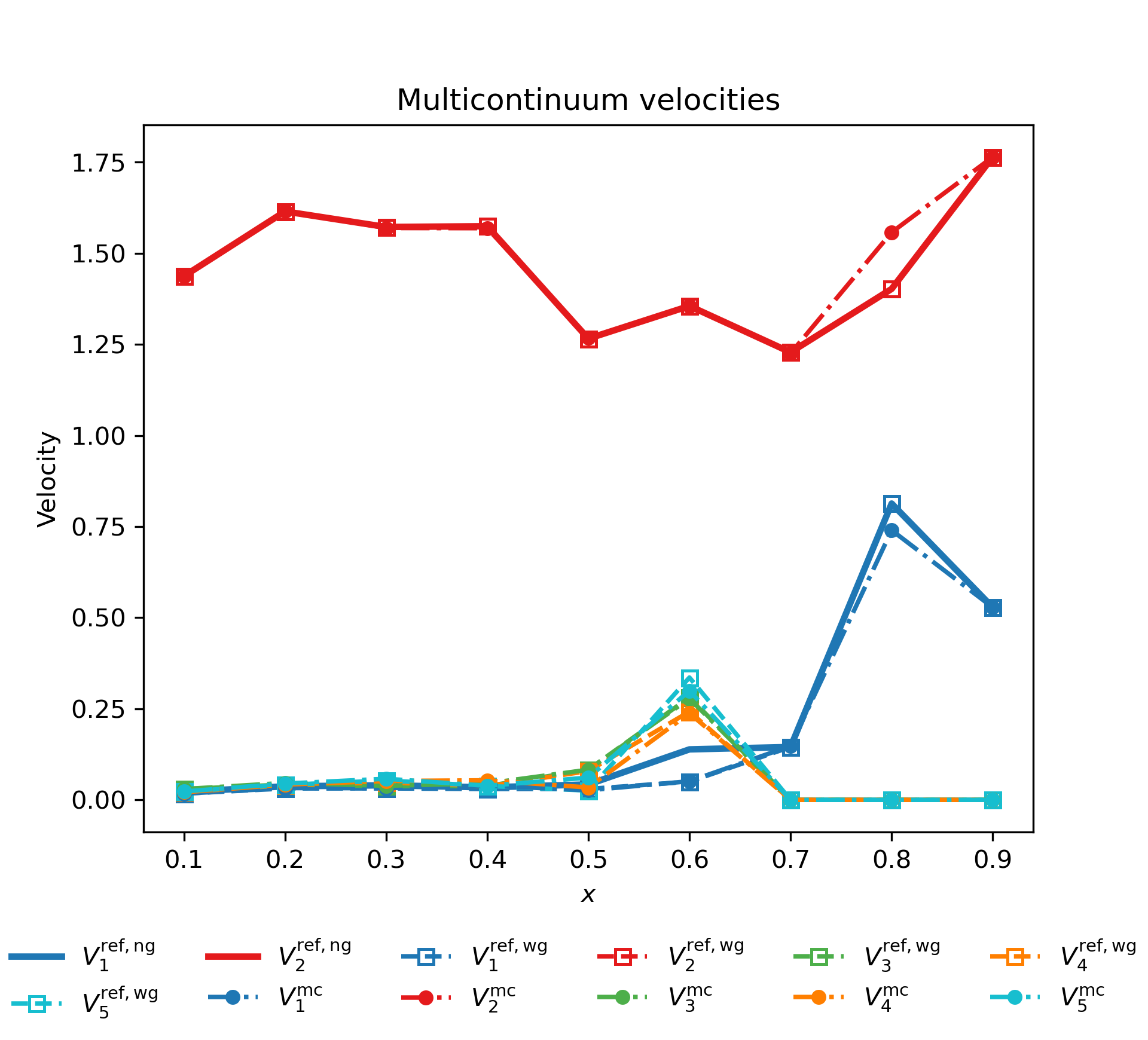}
\includegraphics[width=0.48\textwidth]{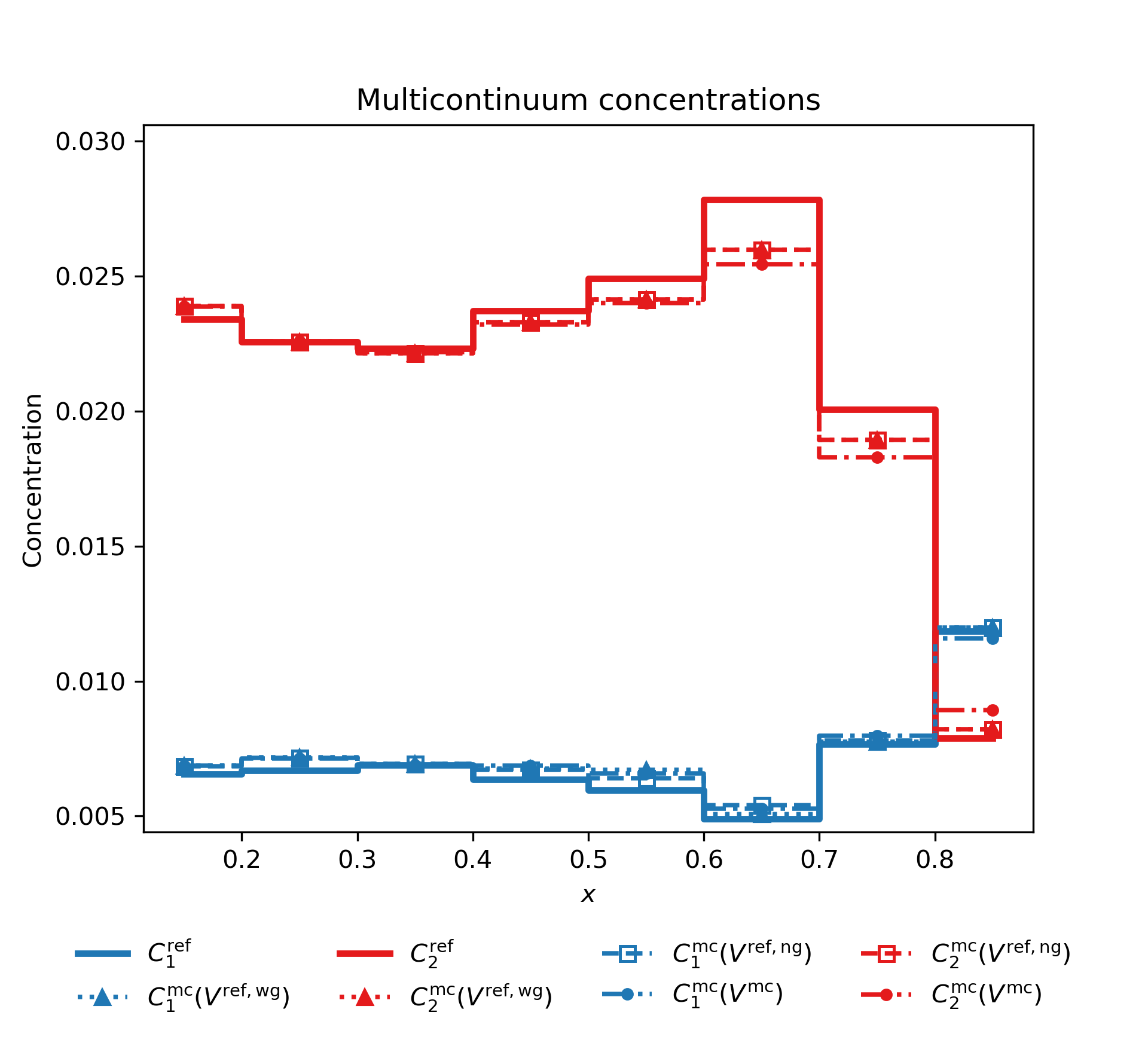}
\includegraphics[width=0.48\textwidth]{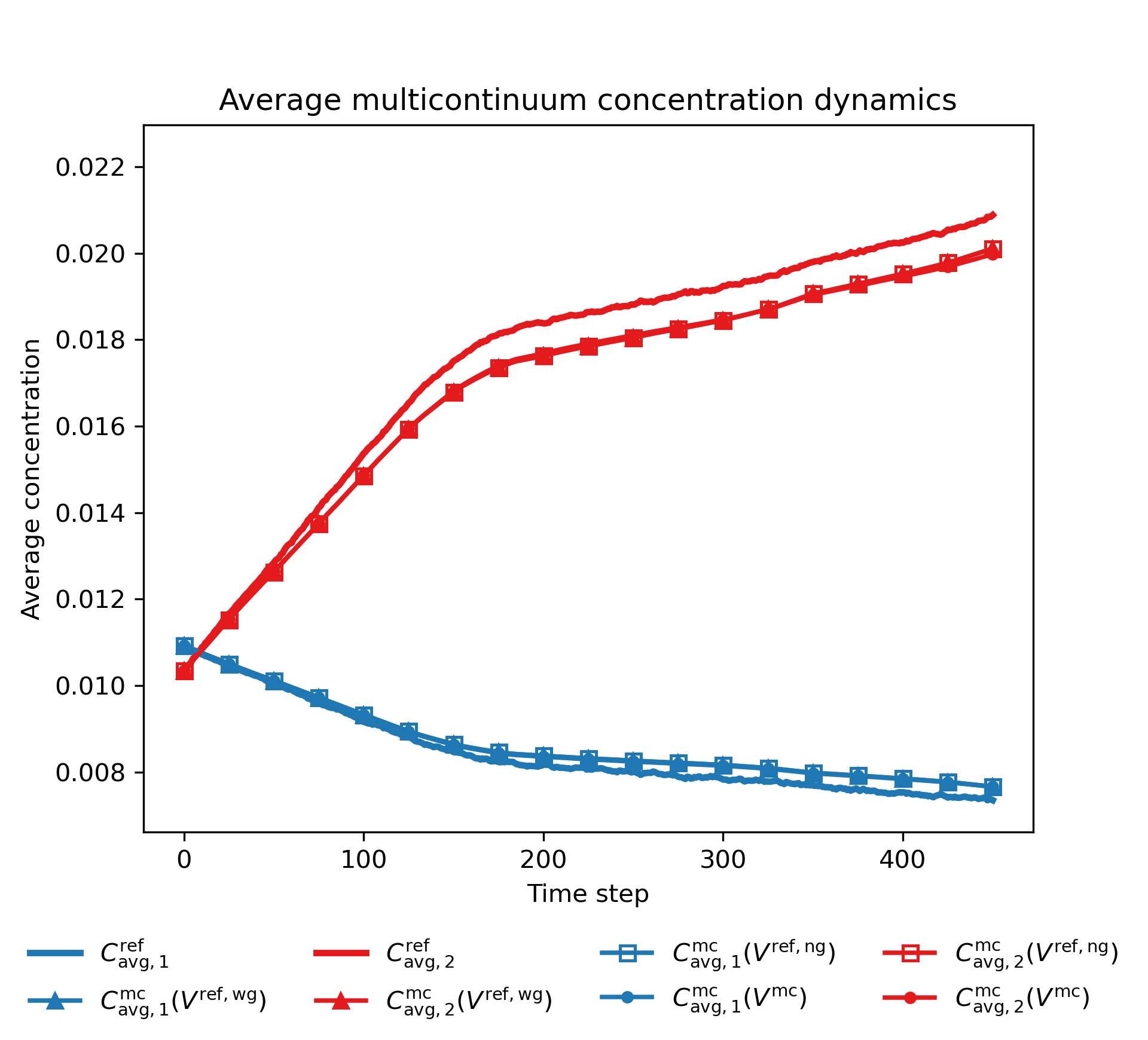}
\caption{Multicontinuum velocities, concentrations, and average concentration dynamics for Test B.}
\label{fig:two_phase_test_b_plots}
\end{figure}

In Figure \ref{fig:two_phase_test_b_plots}, we show the multicontinuum velocities, concentrations, and average concentration dynamics. The corresponding curves agree well in all plots, indicating the high accuracy of the MC-GMsFEM. Compared to the previous case, the curves are steeper due to the influence of the larger target domain.

Tables \ref{tabs:two_phase_test_b_V_errors}, \ref{tabs:two_phase_test_b_C_errors}, and \ref{tabs:two_phase_test_b_C_tracer_errors} present the relative errors for the multicontinuum velocities, concentrations, and average concentration dynamics, respectively. Since the structure of these comparisons follows that of Test A, we focus only on the main observations. All errors remain small, including those for the ganglia velocities. The concentration errors indicate that the effect of ganglia separation is negligible, while the last rows of the concentration and average-concentration tables show that the MC-GMsFEM approximation error remains low. These results confirm that the proposed approach accurately captures both the multicontinuum velocities and the corresponding concentration dynamics in a larger domain.
\begin{table}[hbt!]
\caption{Relative errors of multicontinuum velocities at the final time for Test B.}
\label{tabs:two_phase_test_b_V_errors}
\centering
\begin{tabular}{lll}
\hline
 Velocity error $e_V^{(i)}$                                                                                   & $e_V^{(1)}$   & $e_V^{(2)}$   \\
\hline
 $\|V_i^{\mathrm{ref, wg}} - V_i^{\mathrm{ref, ng}}\|_2 / \|V_i^{\mathrm{ref, ng}}\|_2 \times 100\%$ & 9.10 \%       & 0.00 \%       \\
 $\|V_i^{\mathrm{mc}} - V_i^{\mathrm{ref, ng}}\|_2 / \|V_i^{\mathrm{ref, ng}}\|_2 \times 100\%$      & 11.65 \%      & 3.51 \%       \\ \addlinespace
 $\|V_i^{\mathrm{mc}} - V_i^{\mathrm{ref, wg}}\|_2 / \|V_i^{\mathrm{ref, wg}}\|_2 \times 100\%$      & 7.52 \%       & 3.51 \%       \\
\hline
\end{tabular}

\vspace{0.4em}

\begin{tabular}{llll}
\hline
 Velocity error for ganglia $e_V^{(i)}$                                                                                   & $e_V^{(3)}$   & $e_V^{(4)}$   & $e_V^{(5)}$   \\
\hline
 $\|V_i^{\mathrm{mc}} - V_i^{\mathrm{ref, wg}}\|_2 / \|V_i^{\mathrm{ref, wg}}\|_2 \times 100\%$      & 3.06 \%       & 18.18 \%      & 14.77 \%      \\
\hline
\end{tabular}
\end{table}

\begin{table}[hbt!]
\caption{Relative errors of multicontinuum concentrations at the final time for Test B.}
\label{tabs:two_phase_test_b_C_errors}
\centering
\begin{tabular}{lll}
\hline
 Concentration error $e_C^{(i)}$                                                                                                                                      & $e_C^{(1)}$   & $e_C^{(2)}$   \\
\hline
 $\|C_i^{\mathrm{mc}}(V^{\mathrm{ref, ng}}) - C_i^{\mathrm{ref}}\|_2 / \|C_i^{\mathrm{ref}}\|_2 \times 100\%$                                           & 4.58 \%       & 3.85 \%       \\
 $\|C_i^{\mathrm{mc}}(V^{\mathrm{ref, wg}}) - C_i^{\mathrm{ref}}\|_2 / \|C_i^{\mathrm{ref}}\|_2 \times 100\%$                                           & 4.99 \%       & 3.85 \%       \\
 $\|C_i^{\mathrm{mc}}(V^{\mathrm{mc}}) - C_i^{\mathrm{ref}}\|_2 / \|C_i^{\mathrm{ref}}\|_2 \times 100\%$                                                & 5.34 \%       & 5.29 \%       \\ \addlinespace
 $\|C_i^{\mathrm{mc}}(V^{\mathrm{ref, wg}}) - C_i^{\mathrm{mc}}(V^{\mathrm{ref, ng}})\|_2 / \|C_i^{\mathrm{mc}}(V^{\mathrm{ref, ng}})\|_2 \times 100\%$ & 2.07 \%       & 0.00 \%       \\
 $\|C_i^{\mathrm{mc}}(V^{\mathrm{mc}}) - C_i^{\mathrm{mc}}(V^{\mathrm{ref, ng}})\|_2 / \|C_i^{\mathrm{mc}}(V^{\mathrm{ref, ng}})\|_2 \times 100\%$      & 2.42 \%       & 1.77 \%       \\ \addlinespace
 $\|C_i^{\mathrm{mc}}(V^{\mathrm{mc}}) - C_i^{\mathrm{mc}}(V^{\mathrm{ref, wg}})\|_2 / \|C_i^{\mathrm{mc}}(V^{\mathrm{ref, wg}})\|_2 \times 100\%$      & 2.41 \%       & 1.77 \%       \\
\hline
\end{tabular}
\end{table}

\begin{table}[hbt!]
\caption{Relative errors of average multicontinuum concentration dynamics for Test B.}
\label{tabs:two_phase_test_b_C_tracer_errors}
\centering
\begin{tabular}{lll}
\hline
 Average concentration dynamics error $e_{C, \mathrm{avg}}^{(i)}$                                                                                                                                                                      & $e_{C, \mathrm{avg}}^{(1)}$   & $e_{C, \mathrm{avg}}^{(2)}$   \\
\hline
 $\|C_{\mathrm{avg}, i}^{\mathrm{mc}}(V^{\mathrm{ref, ng}}) - C_{\mathrm{avg}, i}^{\mathrm{ref}}\|_2 / \|C_{\mathrm{avg}, i}^{\mathrm{ref}}\|_2 \times 100\%$                                           & 2.91 \%                       & 3.73 \%                       \\
 $\|C_{\mathrm{avg}, i}^{\mathrm{mc}}(V^{\mathrm{ref, wg}}) - C_{\mathrm{avg}, i}^{\mathrm{ref}}\|_2 / \|C_{\mathrm{avg}, i}^{\mathrm{ref}}\|_2 \times 100\%$                                           & 2.89 \%                       & 3.73 \%                       \\
 $\|C_{\mathrm{avg}, i}^{\mathrm{mc}}(V^{\mathrm{mc}}) - C_{\mathrm{avg}, i}^{\mathrm{ref}}\|_2 / \|C_{\mathrm{avg}, i}^{\mathrm{ref}}\|_2 \times 100\%$                                                & 2.75 \%                       & 3.73 \%                       \\ \addlinespace
 $\|C_{\mathrm{avg}, i}^{\mathrm{mc}}(V^{\mathrm{ref, wg}}) - C_{\mathrm{avg}, i}^{\mathrm{mc}}(V^{\mathrm{ref, ng}})\|_2 / \|C_{\mathrm{avg}, i}^{\mathrm{mc}}(V^{\mathrm{ref, ng}})\|_2 \times 100\%$ & 0.02 \%                       & 0.00 \%                       \\
 $\|C_{\mathrm{avg}, i}^{\mathrm{mc}}(V^{\mathrm{mc}}) - C_{\mathrm{avg}, i}^{\mathrm{mc}}(V^{\mathrm{ref, ng}})\|_2 / \|C_{\mathrm{avg}, i}^{\mathrm{mc}}(V^{\mathrm{ref, ng}})\|_2 \times 100\%$      & 0.17 \%                       & 0.26 \%                       \\ \addlinespace
 $\|C_{\mathrm{avg}, i}^{\mathrm{mc}}(V^{\mathrm{mc}}) - C_{\mathrm{avg}, i}^{\mathrm{mc}}(V^{\mathrm{ref, wg}})\|_2 / \|C_{\mathrm{avg}, i}^{\mathrm{mc}}(V^{\mathrm{ref, wg}})\|_2 \times 100\%$      & 0.16 \%                       & 0.26 \%                       \\
\hline
\end{tabular}
\end{table}

\FloatBarrier
\subsubsection{Test C}

\begin{figure}[hbt!]
\centering
\includegraphics[width=0.48\textwidth]{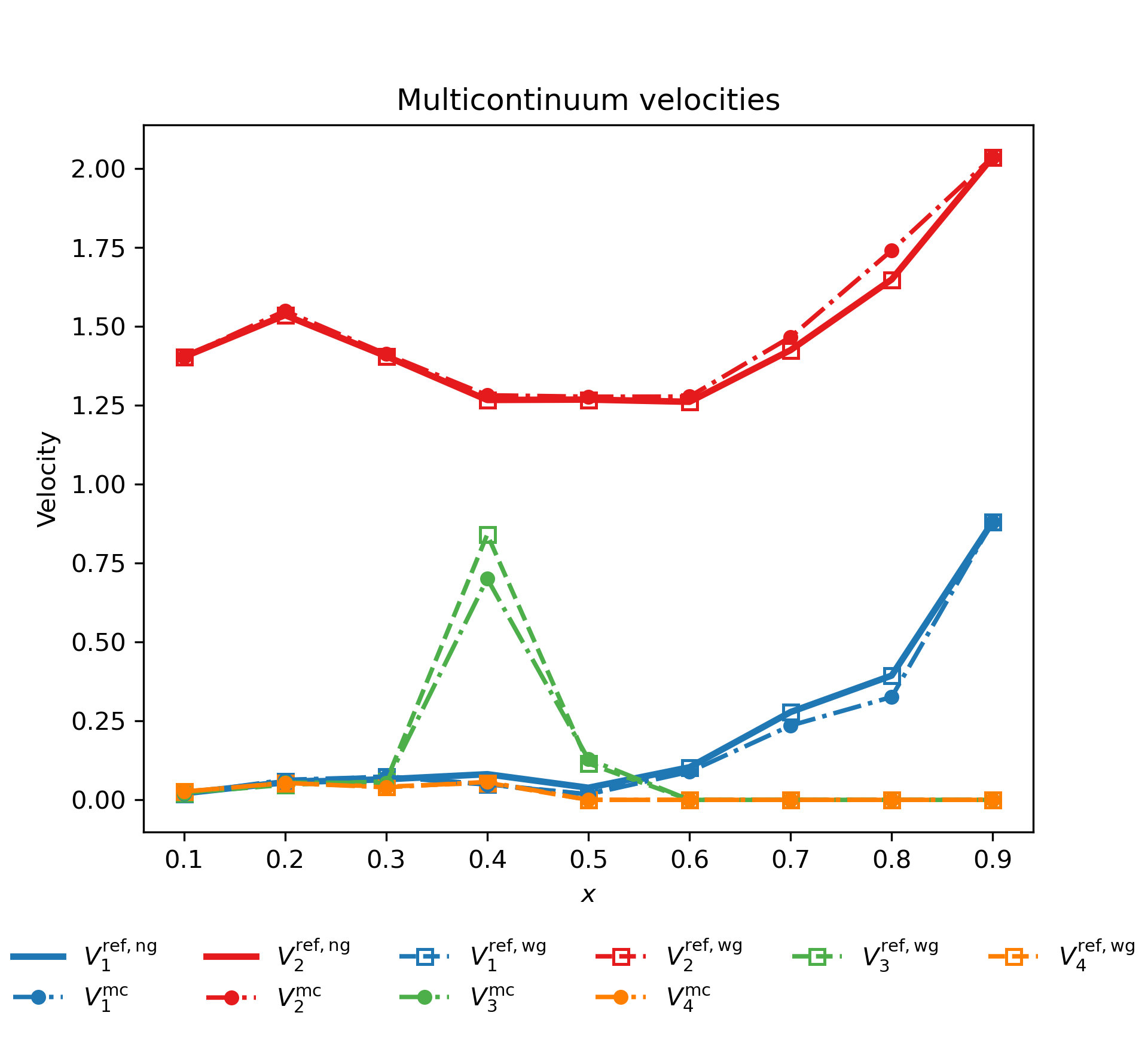}
\includegraphics[width=0.48\textwidth]{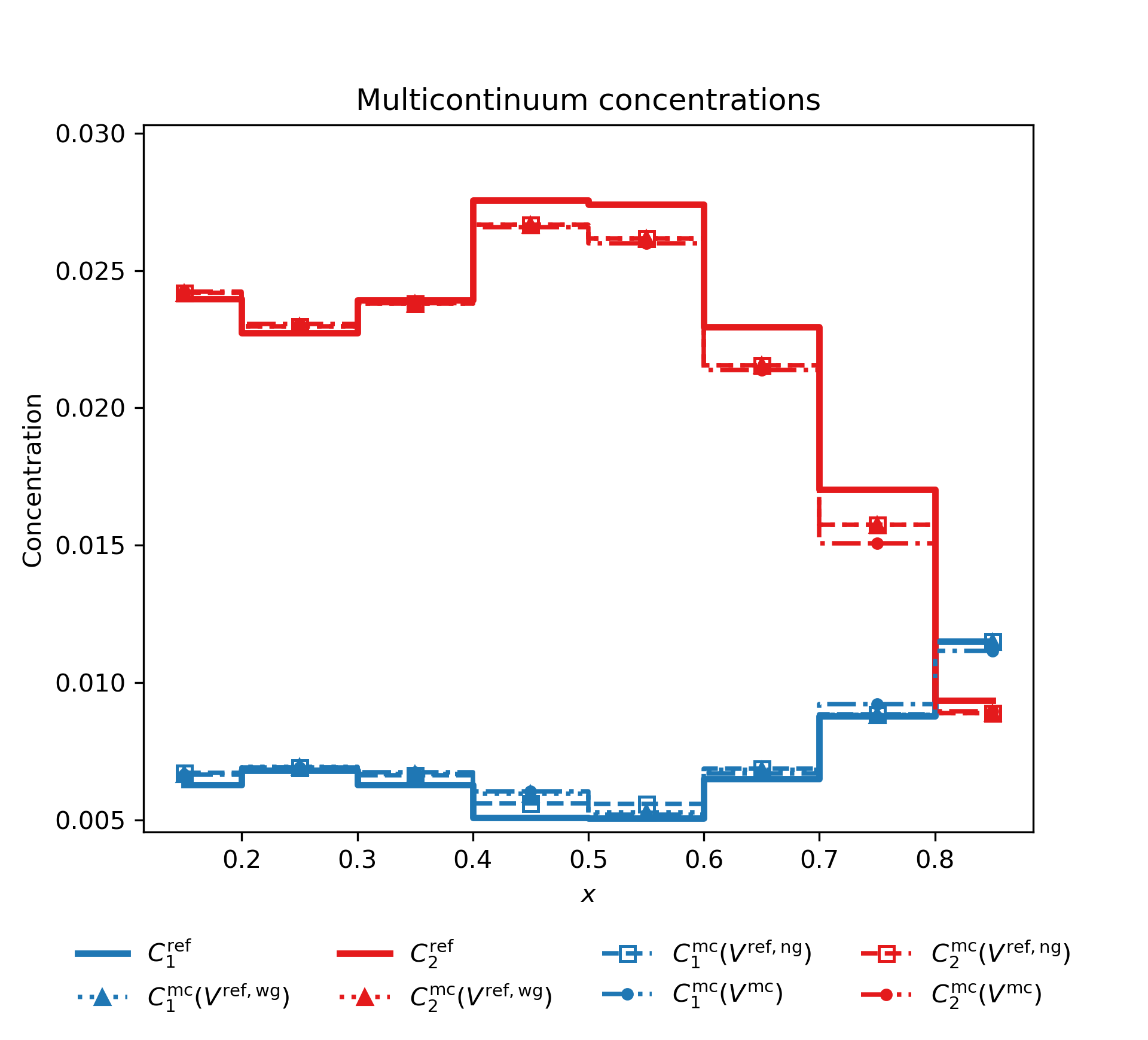}
\includegraphics[width=0.48\textwidth]{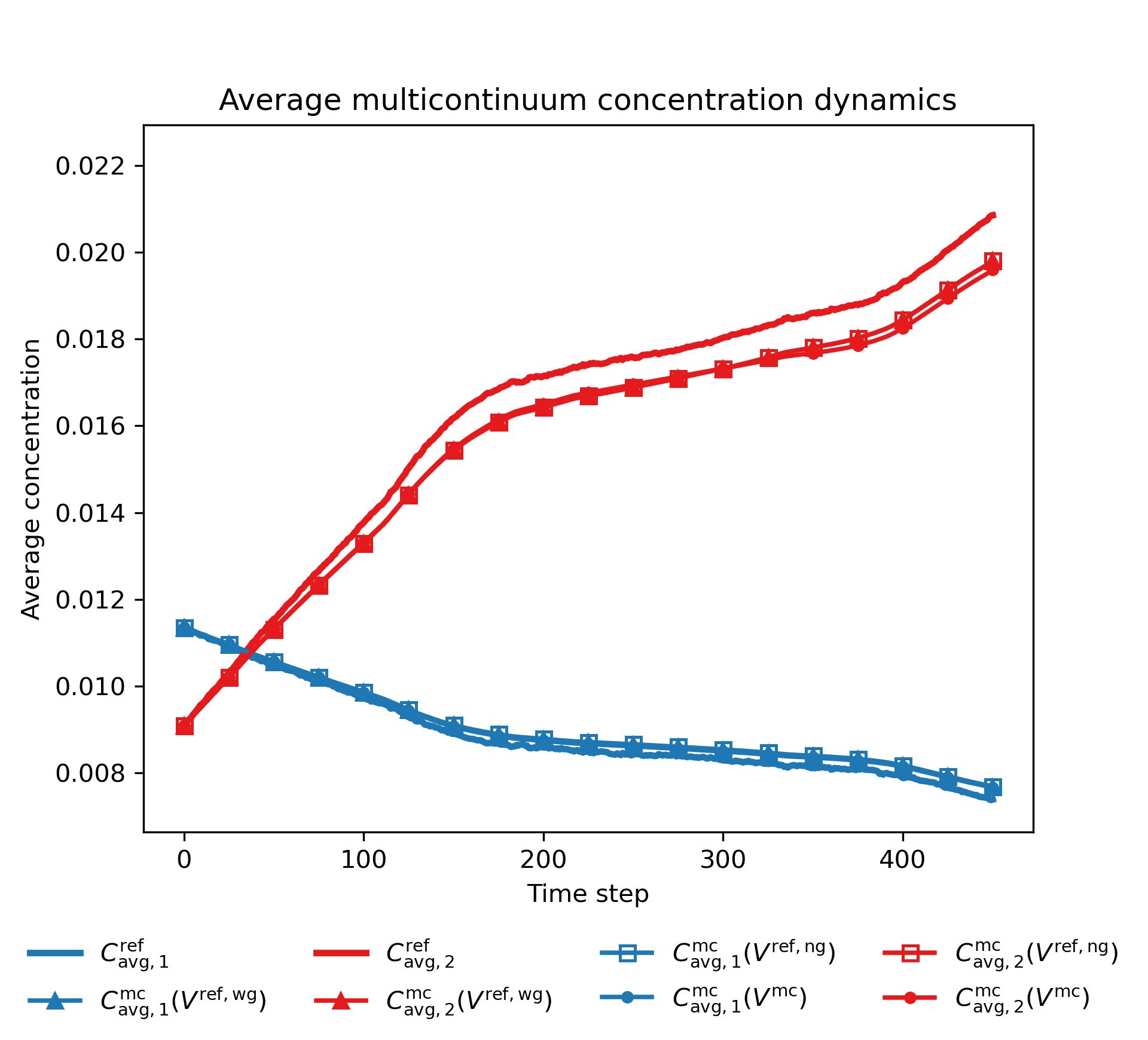}
\caption{Multicontinuum velocities, concentrations, and average concentration dynamics for Test C.}
\label{fig:two_phase_test_c_plots}
\end{figure}

In this case, we also consider the global MC-GMsFEM approach in Target domain 2. The fine-grid concentrations are shown in the third row of Figure \ref{fig:two_phase_fine_results}. Figure \ref{fig:two_phase_test_c_plots} presents the multicontinuum velocities, concentrations, and average concentration dynamics. We observe good agreement between the corresponding solution curves.

Tables \ref{tabs:two_phase_test_c_V_errors}, \ref{tabs:two_phase_test_c_C_errors}, and \ref{tabs:two_phase_test_c_C_tracer_errors} present the relative errors of the multicontinuum velocities, concentrations, and average concentration dynamics, respectively. All errors are small, demonstrating the high accuracy of the proposed MC-GMsFEM.

\begin{table}[hbt!]
\caption{Relative errors of multicontinuum velocities at the final time for Test C.}
\label{tabs:two_phase_test_c_V_errors}
\centering
\begin{tabular}{lll}
\hline
 Velocity error $e_V^{(i)}$                                                                                   & $e_V^{(1)}$   & $e_V^{(2)}$   \\
\hline
 $\|V_i^{\mathrm{ref, wg}} - V_i^{\mathrm{ref, ng}}\|_2 / \|V_i^{\mathrm{ref, ng}}\|_2 \times 100\%$ & 3.85 \%       & 0.00 \%       \\
 $\|V_i^{\mathrm{mc}} - V_i^{\mathrm{ref, ng}}\|_2 / \|V_i^{\mathrm{ref, ng}}\|_2 \times 100\%$      & 8.95 \%       & 2.35 \%       \\ \addlinespace
 $\|V_i^{\mathrm{mc}} - V_i^{\mathrm{ref, wg}}\|_2 / \|V_i^{\mathrm{ref, wg}}\|_2 \times 100\%$      & 8.04 \%       & 2.35 \%       \\
\hline
\end{tabular}

\vspace{0.4em}

\begin{tabular}{lll}
\hline
 Velocity error for ganglia $e_V^{(i)}$                                                                                   & $e_V^{(3)}$   & $e_V^{(4)}$   \\
\hline
 $\|V_i^{\mathrm{mc}} - V_i^{\mathrm{ref, wg}}\|_2 / \|V_i^{\mathrm{ref, wg}}\|_2 \times 100\%$      & 16.53 \%      & 3.75 \%       \\
\hline
\end{tabular}
\end{table}

\begin{table}[hbt!]
\caption{Relative errors of multicontinuum concentrations at the final time for Test C.}
\label{tabs:two_phase_test_c_C_errors}
\centering
\begin{tabular}{lll}
\hline
 Concentration error $e_C^{(i)}$                                                                                                                                      & $e_C^{(1)}$   & $e_C^{(2)}$   \\
\hline
 $\|C_i^{\mathrm{mc}}(V^{\mathrm{ref, ng}}) - C_i^{\mathrm{ref}}\|_2 / \|C_i^{\mathrm{ref}}\|_2 \times 100\%$                                           & 4.85 \%       & 3.91 \%       \\
 $\|C_i^{\mathrm{mc}}(V^{\mathrm{ref, wg}}) - C_i^{\mathrm{ref}}\|_2 / \|C_i^{\mathrm{ref}}\|_2 \times 100\%$                                           & 5.35 \%       & 3.91 \%       \\
 $\|C_i^{\mathrm{mc}}(V^{\mathrm{mc}}) - C_i^{\mathrm{ref}}\|_2 / \|C_i^{\mathrm{ref}}\|_2 \times 100\%$                                                & 6.19 \%       & 4.81 \%       \\ \addlinespace
 $\|C_i^{\mathrm{mc}}(V^{\mathrm{ref, wg}}) - C_i^{\mathrm{mc}}(V^{\mathrm{ref, ng}})\|_2 / \|C_i^{\mathrm{mc}}(V^{\mathrm{ref, ng}})\|_2 \times 100\%$ & 2.21 \%       & 0.00 \%       \\
 $\|C_i^{\mathrm{mc}}(V^{\mathrm{mc}}) - C_i^{\mathrm{mc}}(V^{\mathrm{ref, ng}})\|_2 / \|C_i^{\mathrm{mc}}(V^{\mathrm{ref, ng}})\|_2 \times 100\%$      & 3.77 \%       & 1.15 \%       \\ \addlinespace
 $\|C_i^{\mathrm{mc}}(V^{\mathrm{mc}}) - C_i^{\mathrm{mc}}(V^{\mathrm{ref, wg}})\|_2 / \|C_i^{\mathrm{mc}}(V^{\mathrm{ref, wg}})\|_2 \times 100\%$      & 2.53 \%       & 1.15 \%       \\
\hline
\end{tabular}
\end{table}

\begin{table}[hbt!]
\caption{Relative errors of average multicontinuum concentration dynamics for Test C.}
\label{tabs:two_phase_test_c_C_tracer_errors}
\centering
\begin{tabular}{lll}
\hline
 Average concentration dynamics error $e_{C, \mathrm{avg}}^{(i)}$                                                                                                                                                                      & $e_{C, \mathrm{avg}}^{(1)}$   & $e_{C, \mathrm{avg}}^{(2)}$   \\
\hline
 $\|C_{\mathrm{avg}, i}^{\mathrm{mc}}(V^{\mathrm{ref, ng}}) - C_{\mathrm{avg}, i}^{\mathrm{ref}}\|_2 / \|C_{\mathrm{avg}, i}^{\mathrm{ref}}\|_2 \times 100\%$                                           & 2.31 \%                       & 4.06 \%                       \\
 $\|C_{\mathrm{avg}, i}^{\mathrm{mc}}(V^{\mathrm{ref, wg}}) - C_{\mathrm{avg}, i}^{\mathrm{ref}}\|_2 / \|C_{\mathrm{avg}, i}^{\mathrm{ref}}\|_2 \times 100\%$                                           & 2.23 \%                       & 4.06 \%                       \\
 $\|C_{\mathrm{avg}, i}^{\mathrm{mc}}(V^{\mathrm{mc}}) - C_{\mathrm{avg}, i}^{\mathrm{ref}}\|_2 / \|C_{\mathrm{avg}, i}^{\mathrm{ref}}\|_2 \times 100\%$                                                & 1.92 \%                       & 4.31 \%                       \\ \addlinespace
 $\|C_{\mathrm{avg}, i}^{\mathrm{mc}}(V^{\mathrm{ref, wg}}) - C_{\mathrm{avg}, i}^{\mathrm{mc}}(V^{\mathrm{ref, ng}})\|_2 / \|C_{\mathrm{avg}, i}^{\mathrm{mc}}(V^{\mathrm{ref, ng}})\|_2 \times 100\%$ & 0.08 \%                       & 0.00 \%                       \\
 $\|C_{\mathrm{avg}, i}^{\mathrm{mc}}(V^{\mathrm{mc}}) - C_{\mathrm{avg}, i}^{\mathrm{mc}}(V^{\mathrm{ref, ng}})\|_2 / \|C_{\mathrm{avg}, i}^{\mathrm{mc}}(V^{\mathrm{ref, ng}})\|_2 \times 100\%$      & 0.40 \%                       & 0.60 \%                       \\ \addlinespace
 $\|C_{\mathrm{avg}, i}^{\mathrm{mc}}(V^{\mathrm{mc}}) - C_{\mathrm{avg}, i}^{\mathrm{mc}}(V^{\mathrm{ref, wg}})\|_2 / \|C_{\mathrm{avg}, i}^{\mathrm{mc}}(V^{\mathrm{ref, wg}})\|_2 \times 100\%$      & 0.32 \%                       & 0.60 \%                       \\
\hline
\end{tabular}
\end{table}

\FloatBarrier
\subsubsection{Test D}

\begin{figure}[hbt!]
\centering
\includegraphics[width=0.48\textwidth]{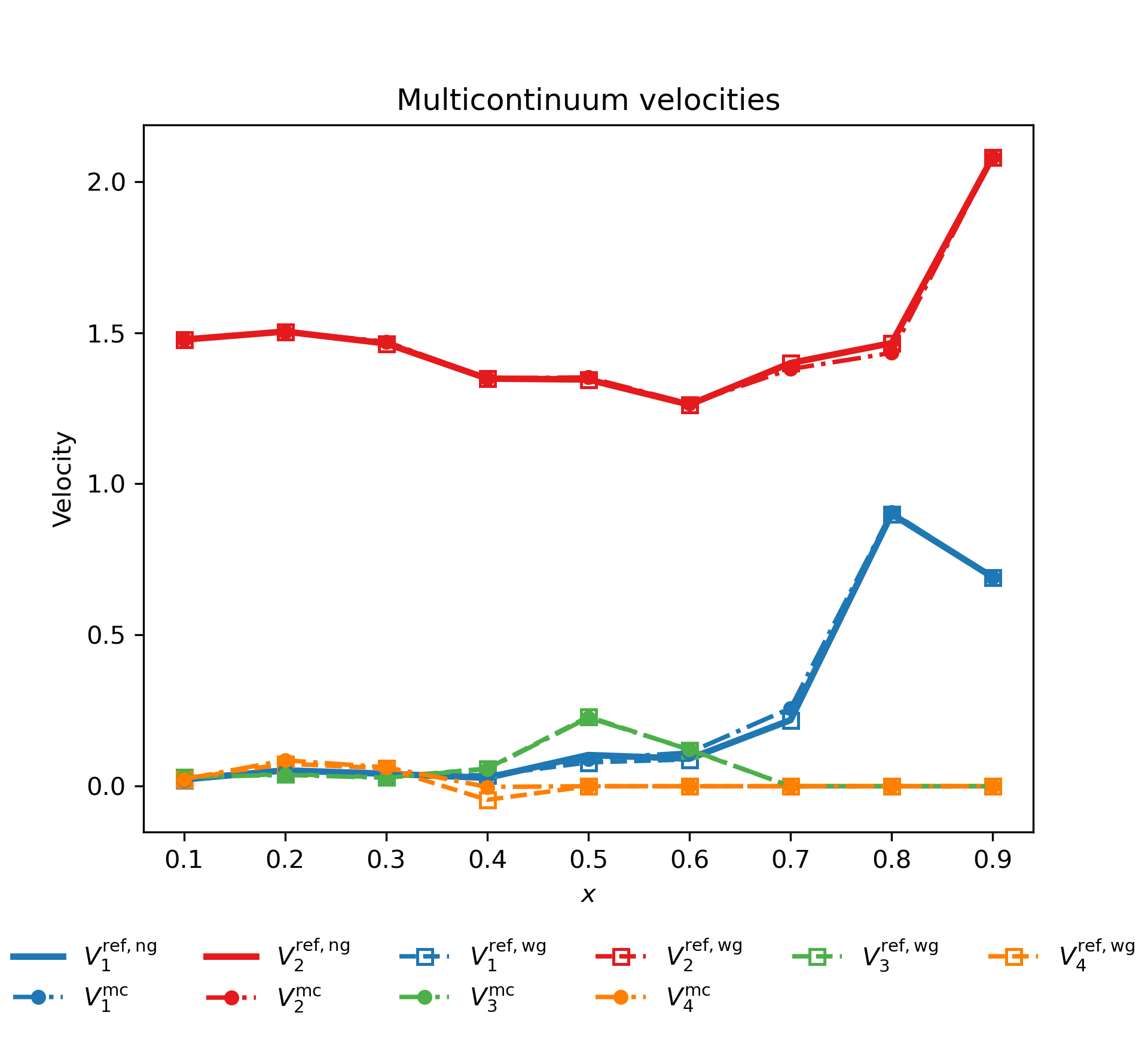}
\includegraphics[width=0.48\textwidth]{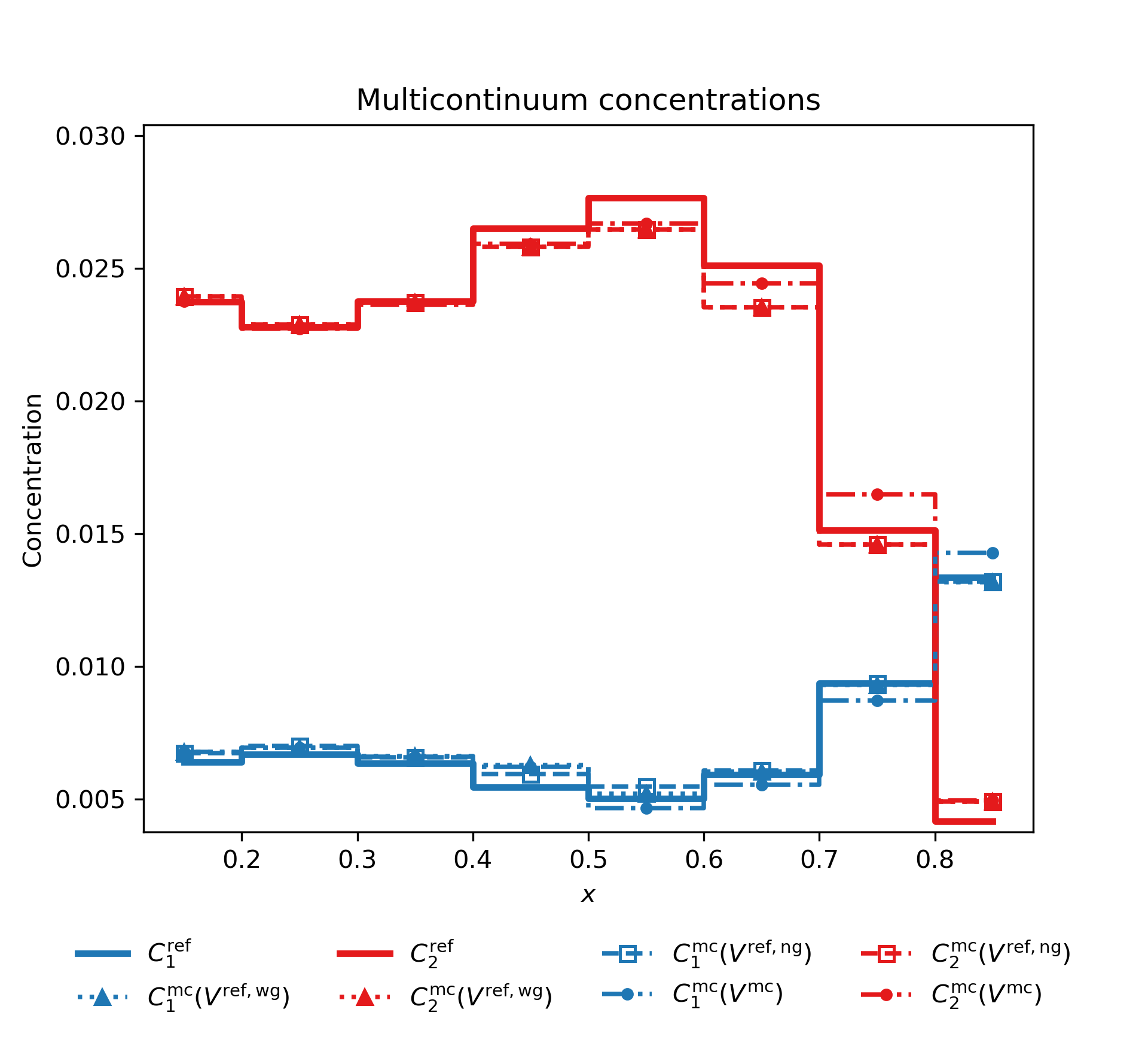}
\includegraphics[width=0.48\textwidth]{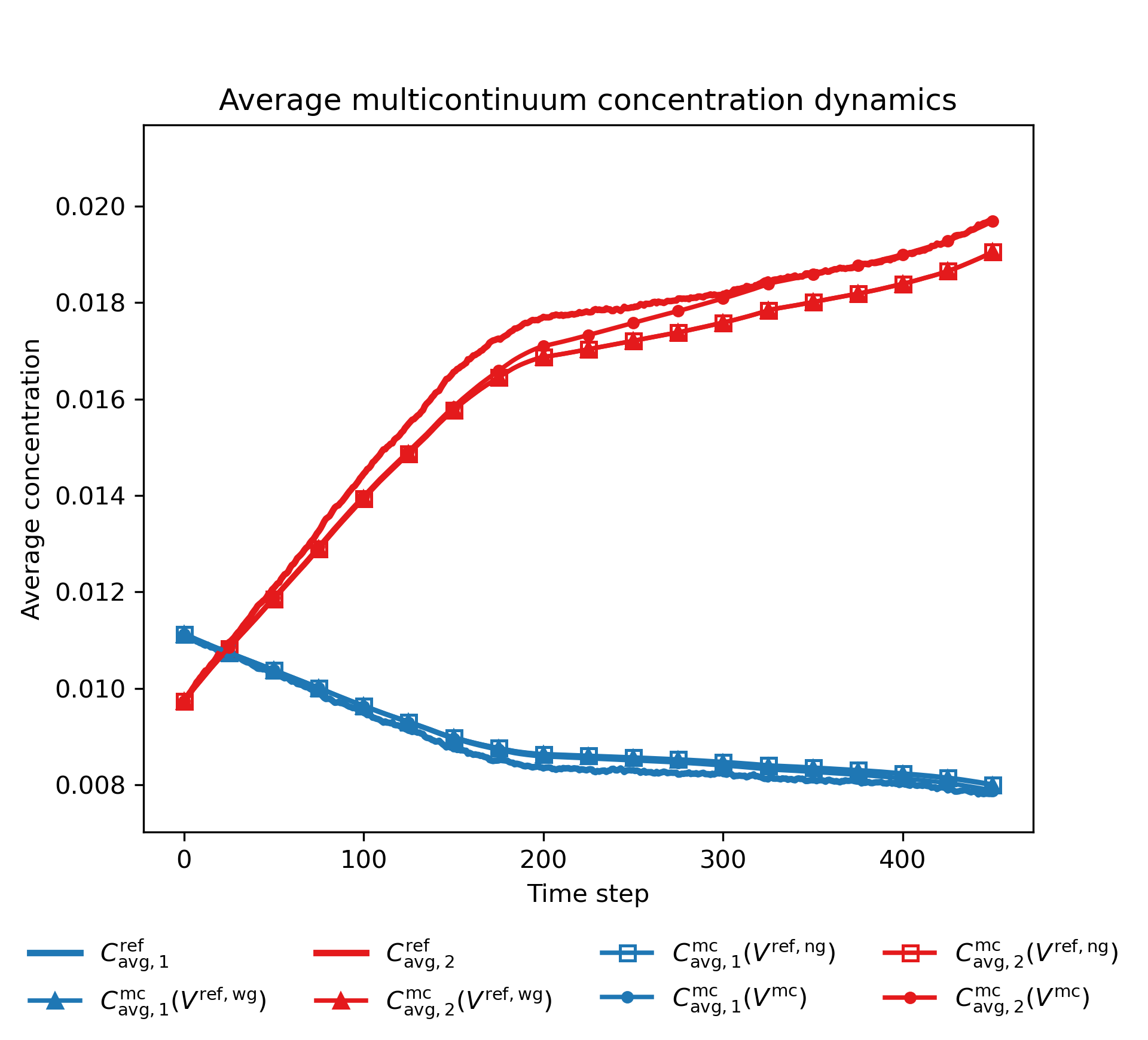}
\caption{Multicontinuum velocities, concentrations, and average concentration dynamics for Test D.}
\label{fig:two_phase_test_d_plots}
\end{figure}

In this test, we consider the local MC-GMsFEM approach in Target domain 2. The fine-grid concentrations are depicted in the fourth row of Figure \ref{fig:two_phase_fine_results}. Figure \ref{fig:two_phase_test_d_plots} presents the multicontinuum velocities, concentrations, and average concentration dynamics. One can see that the corresponding solution curves are close to each other.

Tables \ref{tabs:two_phase_test_d_V_errors}, \ref{tabs:two_phase_test_d_C_errors}, and \ref{tabs:two_phase_test_d_C_tracer_errors} present the relative errors of the multicontinuum velocities, concentrations, and average concentration dynamics. The velocity errors for the two main fluid phases and the first ganglia continuum are small. The relative error for the second ganglia continuum velocity is larger because the corresponding reference velocity is very close to zero (see Figure \ref{fig:two_phase_test_d_plots}); therefore, the relative error does not accurately reflect the approximation quality in this case. The concentration and average concentration dynamics errors are minor.

\begin{table}[hbt!]
\caption{Relative errors of multicontinuum velocities at the final time for Test D.}
\label{tabs:two_phase_test_d_V_errors}
\centering
\begin{tabular}{lll}
\hline
 Velocity error $e_V^{(i)}$                                                                                   & $e_V^{(1)}$   & $e_V^{(2)}$   \\
\hline
 $\|V_i^{\mathrm{ref, wg}} - V_i^{\mathrm{ref, ng}}\|_2 / \|V_i^{\mathrm{ref, ng}}\|_2 \times 100\%$ & 2.23 \%       & 0.00 \%       \\
 $\|V_i^{\mathrm{mc}} - V_i^{\mathrm{ref, ng}}\|_2 / \|V_i^{\mathrm{ref, ng}}\|_2 \times 100\%$      & 4.05 \%       & 0.87 \%       \\ \addlinespace
 $\|V_i^{\mathrm{mc}} - V_i^{\mathrm{ref, wg}}\|_2 / \|V_i^{\mathrm{ref, wg}}\|_2 \times 100\%$      & 4.07 \%       & 0.87 \%       \\
\hline
\end{tabular}

\vspace{0.4em}

\begin{tabular}{lll}
\hline
 Velocity error for ganglia $e_V^{(i)}$                                                                                   &    $e_V^{(3)}$   & $e_V^{(4)}$   \\
\hline
 $\|V_i^{\mathrm{mc}} - V_i^{\mathrm{ref, wg}}\|_2 / \|V_i^{\mathrm{ref, wg}}\|_2 \times 100\%$      & 1.90 \%       & 41.13 \%      \\
\hline
\end{tabular}
\end{table}

\begin{table}[hbt!]
\caption{Relative errors of multicontinuum concentrations at the final time for Test D.}
\label{tabs:two_phase_test_d_C_errors}
\centering
\begin{tabular}{lll}
\hline
 Concentration error $e_C^{(i)}$                                                                                                                                      & $e_C^{(1)}$   & $e_C^{(2)}$   \\
\hline
 $\|C_i^{\mathrm{mc}}(V^{\mathrm{ref, ng}}) - C_i^{\mathrm{ref}}\|_2 / \|C_i^{\mathrm{ref}}\|_2 \times 100\%$                                           & 3.96 \%       & 3.63 \%       \\
 $\|C_i^{\mathrm{mc}}(V^{\mathrm{ref, wg}}) - C_i^{\mathrm{ref}}\|_2 / \|C_i^{\mathrm{ref}}\|_2 \times 100\%$                                           & 4.68 \%       & 3.63 \%       \\
 $\|C_i^{\mathrm{mc}}(V^{\mathrm{mc}}) - C_i^{\mathrm{ref}}\|_2 / \|C_i^{\mathrm{ref}}\|_2 \times 100\%$                                                & 7.02 \%       & 3.21 \%       \\ \addlinespace
 $\|C_i^{\mathrm{mc}}(V^{\mathrm{ref, wg}}) - C_i^{\mathrm{mc}}(V^{\mathrm{ref, ng}})\|_2 / \|C_i^{\mathrm{mc}}(V^{\mathrm{ref, ng}})\|_2 \times 100\%$ & 2.00 \%       & 0.00 \%       \\
 $\|C_i^{\mathrm{mc}}(V^{\mathrm{mc}}) - C_i^{\mathrm{mc}}(V^{\mathrm{ref, ng}})\|_2 / \|C_i^{\mathrm{mc}}(V^{\mathrm{ref, ng}})\|_2 \times 100\%$      & 7.18 \%       & 3.44 \%       \\ \addlinespace
 $\|C_i^{\mathrm{mc}}(V^{\mathrm{mc}}) - C_i^{\mathrm{mc}}(V^{\mathrm{ref, wg}})\|_2 / \|C_i^{\mathrm{mc}}(V^{\mathrm{ref, wg}})\|_2 \times 100\%$      & 6.45 \%       & 3.44 \%       \\
\hline
\end{tabular}
\end{table}

\begin{table}[hbt!]
\caption{Relative errors of average multicontinuum concentration dynamics for Test D.}
\label{tabs:two_phase_test_d_C_tracer_errors}
\centering
\begin{tabular}{lll}
\hline
 Average concentration dynamics error $e_{C, \mathrm{avg}}^{(i)}$                                                                                                                                                                      & $e_{C, \mathrm{avg}}^{(1)}$   & $e_{C, \mathrm{avg}}^{(2)}$   \\
\hline
 $\|C_{\mathrm{avg}, i}^{\mathrm{mc}}(V^{\mathrm{ref, ng}}) - C_{\mathrm{avg}, i}^{\mathrm{ref}}\|_2 / \|C_{\mathrm{avg}, i}^{\mathrm{ref}}\|_2 \times 100\%$                                           & 2.55 \%                       & 3.62 \%                       \\
 $\|C_{\mathrm{avg}, i}^{\mathrm{mc}}(V^{\mathrm{ref, wg}}) - C_{\mathrm{avg}, i}^{\mathrm{ref}}\|_2 / \|C_{\mathrm{avg}, i}^{\mathrm{ref}}\|_2 \times 100\%$                                           & 2.45 \%                       & 3.62 \%                       \\
 $\|C_{\mathrm{avg}, i}^{\mathrm{mc}}(V^{\mathrm{mc}}) - C_{\mathrm{avg}, i}^{\mathrm{ref}}\|_2 / \|C_{\mathrm{avg}, i}^{\mathrm{ref}}\|_2 \times 100\%$                                                & 1.92 \%                       & 2.15 \%                       \\ \addlinespace
 $\|C_{\mathrm{avg}, i}^{\mathrm{mc}}(V^{\mathrm{ref, wg}}) - C_{\mathrm{avg}, i}^{\mathrm{mc}}(V^{\mathrm{ref, ng}})\|_2 / \|C_{\mathrm{avg}, i}^{\mathrm{mc}}(V^{\mathrm{ref, ng}})\|_2 \times 100\%$ & 0.12 \%                       & 0.00 \%                       \\
 $\|C_{\mathrm{avg}, i}^{\mathrm{mc}}(V^{\mathrm{mc}}) - C_{\mathrm{avg}, i}^{\mathrm{mc}}(V^{\mathrm{ref, ng}})\|_2 / \|C_{\mathrm{avg}, i}^{\mathrm{mc}}(V^{\mathrm{ref, ng}})\|_2 \times 100\%$      & 0.73 \%                       & 2.39 \%                       \\ \addlinespace
 $\|C_{\mathrm{avg}, i}^{\mathrm{mc}}(V^{\mathrm{mc}}) - C_{\mathrm{avg}, i}^{\mathrm{mc}}(V^{\mathrm{ref, wg}})\|_2 / \|C_{\mathrm{avg}, i}^{\mathrm{mc}}(V^{\mathrm{ref, wg}})\|_2 \times 100\%$      & 0.61 \%                       & 2.39 \%                       \\
\hline
\end{tabular}
\end{table}

\FloatBarrier
\subsubsection{Test E}

\begin{figure}[hbt!]
\centering
\includegraphics[width=0.48\textwidth]{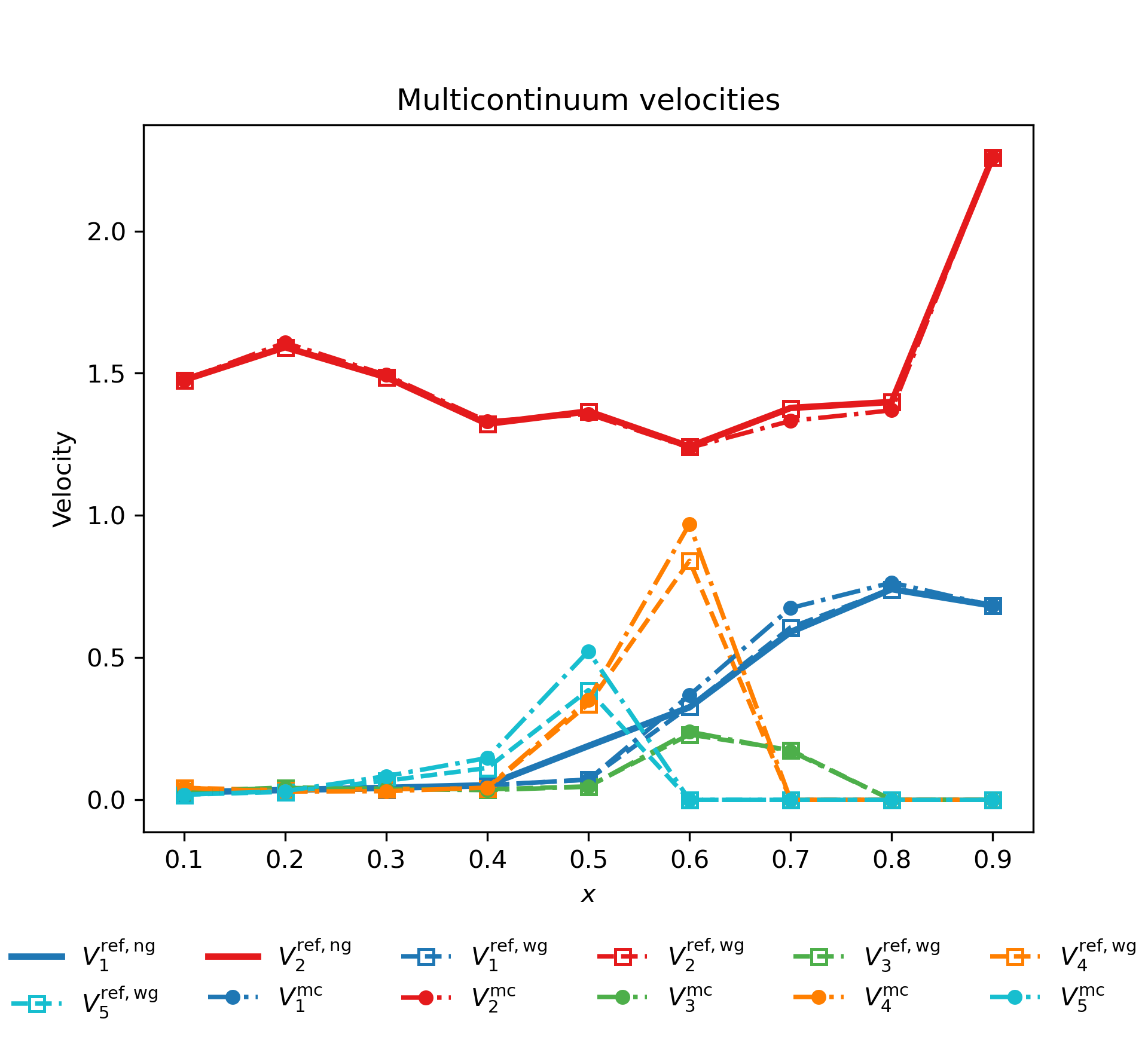}
\includegraphics[width=0.48\textwidth]{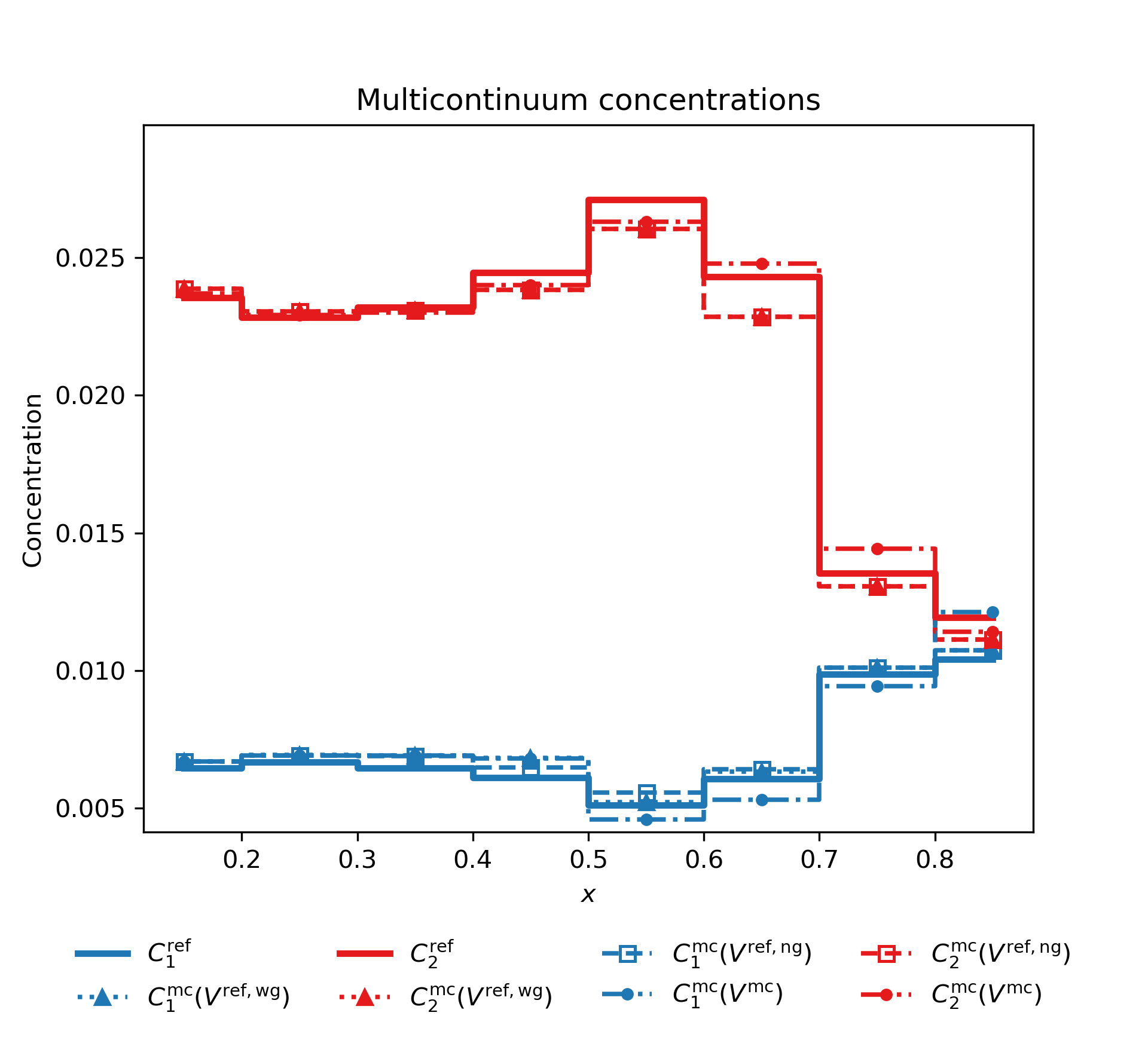}
\includegraphics[width=0.48\textwidth]{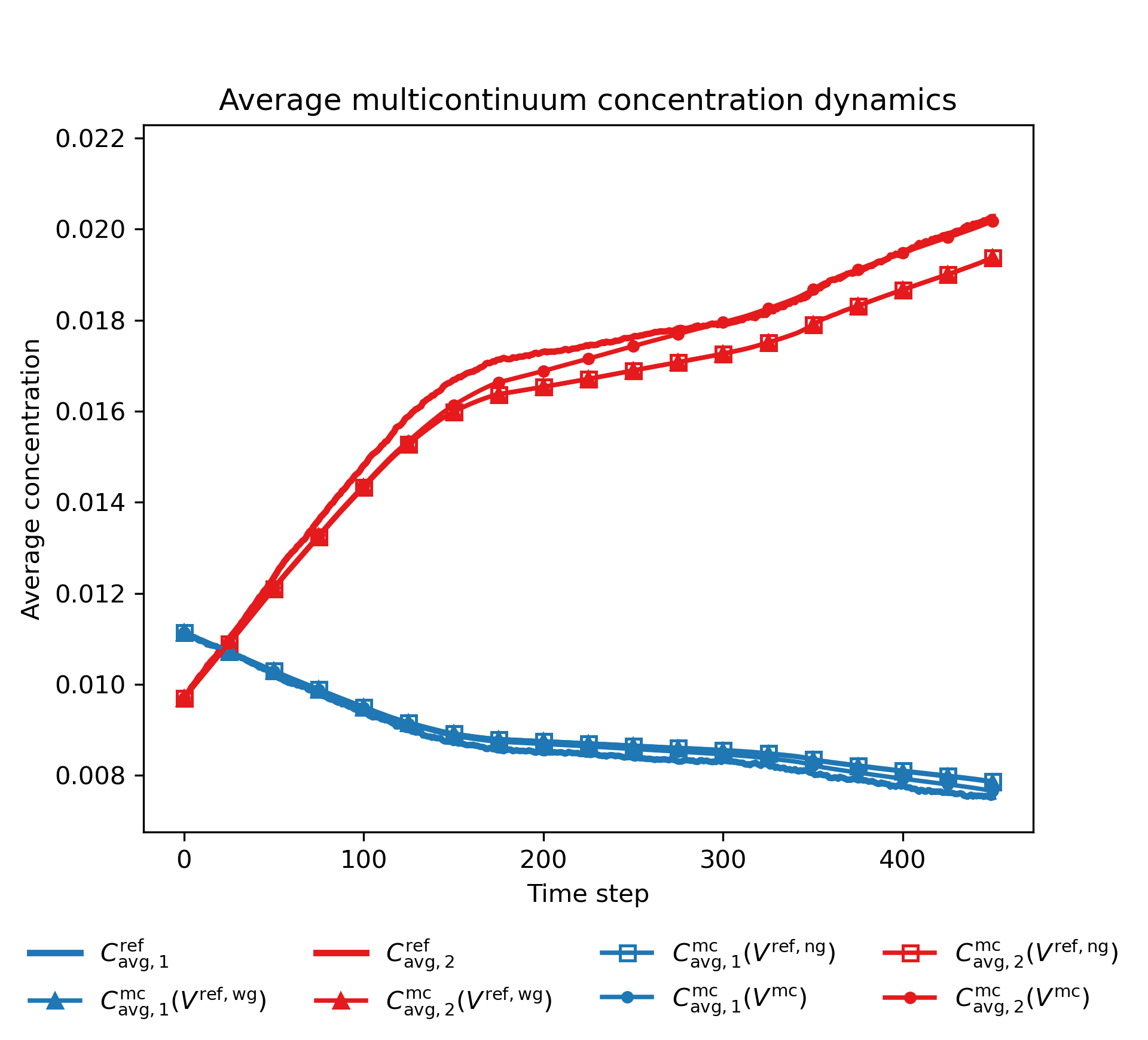}
\caption{Multicontinuum velocities, concentrations, and average concentration dynamics for Test E.}
\label{fig:two_phase_test_e_plots}
\end{figure}

This test also considers the local approach in Target domain 2. The fine-grid concentrations are depicted in the last row of Figure \ref{fig:two_phase_fine_results}. Figure \ref{fig:two_phase_test_e_plots} presents the multicontinuum velocities, concentrations, and average concentration dynamics. One can observe good agreement between the corresponding solution curves, indicating the high accuracy of the MC-GMsFEM approach.

Tables \ref{tabs:two_phase_test_e_V_errors}, \ref{tabs:two_phase_test_e_C_errors}, and \ref{tabs:two_phase_test_e_C_tracer_errors} show the relative errors of the multicontinuum velocities, concentrations, and average concentration dynamics, respectively. The velocity errors for the fluid phases are small, and the ganglia velocity errors are also reasonable. The multicontinuum concentration and average concentration dynamics errors are minor. Thus, the local MC-GMsFEM approach can also provide accurate solutions.

\begin{table}[hbt!]
\caption{Relative errors of multicontinuum velocities at the final time for Test E.}
\label{tabs:two_phase_test_e_V_errors}
\centering
\begin{tabular}{lll}
\hline
 Velocity error $e_V^{(i)}$                                                                                   & $e_V^{(1)}$   & $e_V^{(2)}$   \\
\hline
 $\|V_i^{\mathrm{ref, wg}} - V_i^{\mathrm{ref, ng}}\|_2 / \|V_i^{\mathrm{ref, ng}}\|_2 \times 100\%$ & 9.72 \%       & 0.00 \%       \\
 $\|V_i^{\mathrm{mc}} - V_i^{\mathrm{ref, ng}}\|_2 / \|V_i^{\mathrm{ref, ng}}\|_2 \times 100\%$      & 12.48 \%      & 1.28 \%       \\ \addlinespace
 $\|V_i^{\mathrm{mc}} - V_i^{\mathrm{ref, wg}}\|_2 / \|V_i^{\mathrm{ref, wg}}\|_2 \times 100\%$      & 6.72 \%       & 1.28 \%       \\
\hline
\end{tabular}

\vspace{0.4em}

\begin{tabular}{llll}
\hline
 Velocity error for ganglia $e_V^{(i)}$                                                                                   &  $e_V^{(3)}$   & $e_V^{(4)}$   & $e_V^{(5)}$   \\
\hline
 $\|V_i^{\mathrm{mc}} - V_i^{\mathrm{ref, wg}}\|_2 / \|V_i^{\mathrm{ref, wg}}\|_2 \times 100\%$      &  3.70 \%       & 14.31 \%      & 34.76 \%      \\
\hline
\end{tabular}
\end{table}

\begin{table}[hbt!]
\caption{Relative errors of multicontinuum concentrations at the final time for Test E.}
\label{tabs:two_phase_test_e_C_errors}
\centering
\begin{tabular}{lll}
\hline
 Concentration error $e_C^{(i)}$                                                                                                                                      & $e_C^{(1)}$   & $e_C^{(2)}$   \\
\hline
 $\|C_i^{\mathrm{mc}}(V^{\mathrm{ref, ng}}) - C_i^{\mathrm{ref}}\|_2 / \|C_i^{\mathrm{ref}}\|_2 \times 100\%$                                           & 4.60 \%       & 3.49 \%       \\
 $\|C_i^{\mathrm{mc}}(V^{\mathrm{ref, wg}}) - C_i^{\mathrm{ref}}\|_2 / \|C_i^{\mathrm{ref}}\|_2 \times 100\%$                                           & 5.01 \%       & 3.49 \%       \\
 $\|C_i^{\mathrm{mc}}(V^{\mathrm{mc}}) - C_i^{\mathrm{ref}}\|_2 / \|C_i^{\mathrm{ref}}\|_2 \times 100\%$                                                & 10.52 \%      & 2.40 \%       \\ \addlinespace
 $\|C_i^{\mathrm{mc}}(V^{\mathrm{ref, wg}}) - C_i^{\mathrm{mc}}(V^{\mathrm{ref, ng}})\|_2 / \|C_i^{\mathrm{mc}}(V^{\mathrm{ref, ng}})\|_2 \times 100\%$ & 2.24 \%       & 0.00 \%       \\
 $\|C_i^{\mathrm{mc}}(V^{\mathrm{mc}}) - C_i^{\mathrm{mc}}(V^{\mathrm{ref, ng}})\|_2 / \|C_i^{\mathrm{mc}}(V^{\mathrm{ref, ng}})\|_2 \times 100\%$      & 9.94 \%       & 3.99 \%       \\ \addlinespace
 $\|C_i^{\mathrm{mc}}(V^{\mathrm{mc}}) - C_i^{\mathrm{mc}}(V^{\mathrm{ref, wg}})\|_2 / \|C_i^{\mathrm{mc}}(V^{\mathrm{ref, wg}})\|_2 \times 100\%$      & 9.03 \%       & 3.99 \%       \\
\hline
\end{tabular}
\end{table}

\begin{table}[hbt!]
\caption{Relative errors of average multicontinuum concentration dynamics for Test E.}
\label{tabs:two_phase_test_e_C_tracer_errors}
\centering
\begin{tabular}{lll}
\hline
 Average concentration dynamics error $e_{C, \mathrm{avg}}^{(i)}$                                                                                                                                                                      & $e_{C, \mathrm{avg}}^{(1)}$   & $e_{C, \mathrm{avg}}^{(2)}$   \\
\hline
 $\|C_{\mathrm{avg}, i}^{\mathrm{mc}}(V^{\mathrm{ref, ng}}) - C_{\mathrm{avg}, i}^{\mathrm{ref}}\|_2 / \|C_{\mathrm{avg}, i}^{\mathrm{ref}}\|_2 \times 100\%$                                           & 2.80 \%                       & 3.92 \%                       \\
 $\|C_{\mathrm{avg}, i}^{\mathrm{mc}}(V^{\mathrm{ref, wg}}) - C_{\mathrm{avg}, i}^{\mathrm{ref}}\|_2 / \|C_{\mathrm{avg}, i}^{\mathrm{ref}}\|_2 \times 100\%$                                           & 2.72 \%                       & 3.92 \%                       \\
 $\|C_{\mathrm{avg}, i}^{\mathrm{mc}}(V^{\mathrm{mc}}) - C_{\mathrm{avg}, i}^{\mathrm{ref}}\|_2 / \|C_{\mathrm{avg}, i}^{\mathrm{ref}}\|_2 \times 100\%$                                                & 1.75 \%                       & 1.70 \%                       \\ \addlinespace
 $\|C_{\mathrm{avg}, i}^{\mathrm{mc}}(V^{\mathrm{ref, wg}}) - C_{\mathrm{avg}, i}^{\mathrm{mc}}(V^{\mathrm{ref, ng}})\|_2 / \|C_{\mathrm{avg}, i}^{\mathrm{mc}}(V^{\mathrm{ref, ng}})\|_2 \times 100\%$ & 0.10 \%                       & 0.00 \%                       \\
 $\|C_{\mathrm{avg}, i}^{\mathrm{mc}}(V^{\mathrm{mc}}) - C_{\mathrm{avg}, i}^{\mathrm{mc}}(V^{\mathrm{ref, ng}})\|_2 / \|C_{\mathrm{avg}, i}^{\mathrm{mc}}(V^{\mathrm{ref, ng}})\|_2 \times 100\%$      & 1.11 \%                       & 3.25 \%                       \\ \addlinespace
 $\|C_{\mathrm{avg}, i}^{\mathrm{mc}}(V^{\mathrm{mc}}) - C_{\mathrm{avg}, i}^{\mathrm{mc}}(V^{\mathrm{ref, wg}})\|_2 / \|C_{\mathrm{avg}, i}^{\mathrm{mc}}(V^{\mathrm{ref, wg}})\|_2 \times 100\%$      & 1.01 \%                       & 3.25 \%                       \\
\hline
\end{tabular}
\end{table}

\FloatBarrier
\subsubsection{Influence of the velocity gradient terms}

In the final numerical experiment, we study the influence of the velocity gradient terms (the second term in the expansion \eqref{eq:vel_exp}) in the macroscopic flow model. These terms arise from the gradient components of the multiscale basis functions. Without them, the macroscopic flow model reduces to the multicontinuum (multiphase) Darcy model, and when using them, we obtain Brinkman-like models. We consider Test C, which uses the global auxiliary functions in Target domain 2. We compute the multicontinuum velocities at the final time with and without the gradient terms to assess their influence. To compare the obtained solutions, we use the following relative error
\begin{equation*}
e_V^{(i)} = \frac{\|V_i^{\mathrm{mc}} - V_i^{\mathrm{ref, wg}}\|_2}{\|V_i^{\mathrm{ref, wg}}\|_2} \times 100\%,
\end{equation*}
where $V_i^{\mathrm{mc}}$ denotes the multicontinuum velocity obtained with or without the gradient terms.

Table \ref{tabs:two_phase_Brinkman_Darcy_V_errors} presents the relative errors of the obtained multicontinuum velocities. Here, $e_V^{(1)}$ and $e_V^{(2)}$ correspond to the two main fluid phases, while $e_V^{(3)}$ and $e_V^{(4)}$ are the ganglia velocity errors. Without gradient terms, the errors for the two main fluid phases are approximately twice as large, although they remain at a reasonable level. More importantly, ignoring the gradient terms leads to significantly larger errors for the ganglia velocities. These results suggest that gradient effects are important for accurately capturing trapped regions.
\begin{table}[hbt!]
\caption{Relative $L^2$ errors of multicontinuum velocities obtained with and without gradient terms.}
\label{tabs:two_phase_Brinkman_Darcy_V_errors}
\centering
\begin{tabular}{lllll}
\hline
 Multiscale basis functions              & $e_V^{(1)}$   & $e_V^{(2)}$   & $e_V^{(3)}$   & $e_V^{(4)}$   \\
\hline
 $N_i^v \phi_0^j + M_i^v \cdot \nabla \phi_0^j$      & 8.04 \%       & 2.35 \%       & 16.53 \%      & 3.75 \%       \\
 $N_i^v \phi_0^j$                                 & 18.92 \%       & 4.32 \%       & 96.13 \%      & 21.62 \%       \\
\hline
\end{tabular}
\end{table}

\FloatBarrier
\section{Conclusions}

In this paper, we introduced the multicontinuum Generalized Multiscale Finite Element Method (MC-GMsFEM), a unified framework that combines the construction of conforming multiscale basis functions with the derivation of physically meaningful macroscopic equations. The proposed approach establishes a direct connection between multiscale finite element approximations and multicontinuum homogenization through the representation. This formulation provides a systematic mechanism for deriving coarse-scale governing equations while preserving important fine-scale information.

The proposed framework extends existing multiscale methods in several important directions. First, it enables the identification and incorporation of multiple interacting continua within a consistent finite element setting. Second, it provides flexibility for heterogeneous media where the number and structure of continua may vary spatially. Third, the resulting formulation remains compatible with standard numerical implementations and can therefore be integrated into existing multiscale simulation methodologies.

We demonstrated the applicability of the method for elliptic problems in heterogeneous and perforated domains. We mainly applied the framework to two-phase immiscible flow problems and showed that trapped and mobile fluid regions can naturally be represented as distinct continua. In such settings, classical coarse-scale Darcy models may fail to accurately capture fine-scale flow mechanisms, particularly in the presence of long disconnected or weakly connected trapped regions. The proposed multicontinuum formulation significantly improves the accuracy of coarse-scale approximations by incorporating additional continuum variables that describe these effects. In two-phase flow modeling, we use fluid regions
as continua, which introduces additional errors as the more accurate continua
regions is defined by the velocity.

Though the proposed methods have advantages over previously proposed multicontinuum homogenization, there are some disadvantages. The main disadvantage is that the proposed method does not 
explicitly define the gradient term and the coefficients in macroscopic equations, which can additionally be derived. The main advantage is that the proposed method can easily be adapted for varying number of continua and general heterogeneities. 

The numerical results demonstrate that MC-GMsFEM provides accurate coarse-scale solutions while simultaneously yielding interpretable macroscopic models. The framework therefore offers both computational efficiency and a physically consistent upscaling methodology for complex multiscale systems.

\bibliographystyle{unsrt}
\bibliography{lit,lit1}

\appendix

\section{Ganglia detection algorithm}\label{sec:ganglia_detection_algorithm}

\begin{algorithm}[hbt!]
\caption{Ganglion extraction from piecewise-constant indicator fields}
\label{alg:ganglia-main}
\DontPrintSemicolon

\AlgInput{Fields $\Psi=(\psi_1,\ldots,\psi_J)$ on $\mathcal{T}_h$;
avoid lines $\mathcal{A}$; minimum length $\ell_{\min}$; weight $\omega$.}
\AlgOutput{Fields
$\widehat{\Psi}
=(\widehat{\psi}_1,\ldots,\widehat{\psi}_J,
\widehat{\psi}^{\mathrm g}_1,\ldots,
\widehat{\psi}^{\mathrm g}_{q_{\mathrm g}})$.}

Mark facets by $b(F)=1$ if $F$ lies on a line in $\mathcal{A}$, and
$b(F)=0$ otherwise.\;
For each element $K$, let $x_K$ be its barycenter and set
$\tau_K=\mathbf{1}_{\{\psi_1(K)>1/2\}}$, $\rho_K=-1$.\;
$G\leftarrow\emptyset$, $S\leftarrow\emptyset$.\;

\ForEach{$K\in\mathcal{T}_h$ with $\rho_K<0$ and $\tau_K=1$}{
    $(\mathcal{R},\sigma)\leftarrow\mathrm{Flood}(K)$\;
    \uIf{$\sigma=g$}{
        append $\mathcal{R}$ to $G$\;
    }
    \ElseIf{$\sigma=s$}{
        append $\mathcal{R}$ to $S$\;
    }
}

\If{$G\neq\emptyset$}{
    Let $G^{0}$ be a fixed copy of $G$.\;
    \ForEach{$\mathcal{I}\in S$}{
        $j^\star\leftarrow
        \arg\min_j\{d(\mathcal{I},G^{0}_j)
        +\omega\,\gamma_y(\mathcal{I},G^{0}_j)\}$\;
        $G_{j^\star}\leftarrow G_{j^\star}\cup\mathcal{I}$\;
    }
}

Set $\widehat{\psi}_i\leftarrow\psi_i$ for $i=1,\ldots,J$ and
$q_{\mathrm g}\leftarrow|G|$.\;
\For{$j=1,\ldots,q_{\mathrm g}$}{
    Set $\widehat{\psi}^{\mathrm g}_j(K)=0$ for all $K\in\mathcal{T}_h$.\;
    \ForEach{$K\in G_j$}{
        $\widehat{\psi}_1(K)\leftarrow0$;
        $\widehat{\psi}^{\mathrm g}_j(K)\leftarrow1$.\;
    }
}

\AlgReturn{$
\widehat{\Psi}
=(\widehat{\psi}_1,\ldots,\widehat{\psi}_J,
\widehat{\psi}^{\mathrm g}_1,\ldots,
\widehat{\psi}^{\mathrm g}_{q_{\mathrm g}})
$}\;
\end{algorithm}

\begin{algorithm}[hbt!]
\caption{$\mathrm{Flood}(K_0)$}
\label{alg:ganglia-flood}
\DontPrintSemicolon

\AlgInput{Start element $K_0$; arrays $\tau_K$, $\rho_K$, and facet marker $b(F)$.}
\AlgOutput{Region $\mathcal{R}$ and label
$\sigma\in\{g,s,d\}$.}

$Q\leftarrow[(K_0,1)]$;
$\rho_{K_0}\leftarrow1$;
$\mathcal{R}\leftarrow\emptyset$;
$\mathrm{admissible}\leftarrow\mathrm{true}$.\;

\While{$Q\neq\emptyset$}{
    Remove the first pair $(K,r)$ from $Q$.\;

    \If{$\tau_K=1$}{
        $r\leftarrow1$;
        $\mathcal{R}\leftarrow\mathcal{R}\cup\{K\}$.\;
        \If{there is $F\in\mathcal{F}(K)$ with $b(F)=1$}{
            $\mathrm{admissible}\leftarrow\mathrm{false}$.\;
        }
    }

    \ForEach{$F\in\mathcal{F}(K)$}{
        \ForEach{$L\in\operatorname{Adj}(F)\setminus\{K\}$}{
            \uIf{$\tau_L=1$}{
                $r'\leftarrow1$.\;
            }
            \uElseIf{$r>0$}{
                $r'\leftarrow r-1$.\;
            }
            \Else{
                \AlgContinue\;
            }

            \If{$r'>\rho_L$}{
                $\rho_L\leftarrow r'$;
                append $(L,r')$ to $Q$.\;
            }
        }
    }
}

\If{$\mathrm{admissible}=\mathrm{false}$ or $\mathcal{R}=\emptyset$}{
    \AlgReturn{$(\emptyset,d)$}\;
}
\uIf{$\ell_x(\mathcal{R})>\ell_{\min}$}{
    \AlgReturn{$(\mathcal{R},g)$}\;
}
\Else{
    \AlgReturn{$(\mathcal{R},s)$}\;
}
\end{algorithm}

We assume that all indicator fields $\psi_i$ are piecewise constant. In Algorithm~\ref{alg:ganglia-main}, $G$ stores the detected ganglia, while $S$ stores small admissible islands before they are attached to the nearest ganglion. We treat the left and right boundaries as avoid lines $\mathcal{A}$. When attaching a small island to the ganglia, we give preference to the ganglion located at approximately the same height in the slab domain. This preference is represented by the weight $\omega$, which is set to $0.25$ in our numerical simulations. In Algorithm~\ref{alg:ganglia-flood}, $Q$ denotes the breadth-first search queue. The label $\sigma=g$ denotes a ganglion, $\sigma=s$ denotes a small island, and $\sigma=d$ denotes a discarded region. Note that the minimal horizontal length of a ganglion $l_{\mathrm{min}}$ is set to be greater than two coarse-grid blocks.

Here $\mathcal{F}(K)$ is the set of facets of element $K$, and $\operatorname{Adj}(F)$ is the set of elements adjacent to facet $F$. Each queue entry $(K,r)$ consists of an element $K$ and the remaining gap budget $r\in\{0,1\}$. The value $r=1$ means that one non-target element may still be crossed. Whenever a target element is reached, the budget is reset to one. Thus the algorithm connects target regions across at most one consecutive non-target element.

For measuring lengths and distances, we use the following notations. First, let
\[
x_K=(x_{K,1},\ldots,x_{K,d})
\]
denote the barycenter of element $K$. 

Then, for a region $\mathcal{R}$, its horizontal length is
\[
\ell_x(\mathcal{R})
=
\max_{K\in\mathcal{R}} x_{K,1}
-
\min_{K\in\mathcal{R}} x_{K,1}.
\]
For two regions $\mathcal{R}$ and $\mathcal{U}$, we define the distance
\[
d(\mathcal{R},\mathcal{U})
=
\min_{K\in\mathcal{R},\,L\in\mathcal{U}}
\|x_K-x_L\|_2 .
\]
The vertical interval gap is
\[
\gamma_y(\mathcal{R},\mathcal{U})
=
\max\{0,\,
a_-(\mathcal{U})-a_+(\mathcal{R}),\,
a_-(\mathcal{R})-a_+(\mathcal{U})\},
\]
where
\[
a_-(\mathcal{R})=\min_{K\in\mathcal{R}} x_{K,2},
\qquad
a_+(\mathcal{R})=\max_{K\in\mathcal{R}} x_{K,2}.
\]

\section{The derivation of multicontinuum equations in perforated domains}
\label{sec:ell_perfor}
Let us consider elliptic problems in perforated domains
\[
\Omega_p=\Omega\setminus \mathcal{B},
\]
where $\mathcal{B}$ denotes the set of perforations. The derivation
is similar to the heterogeneous case and we omit some steps and consider a simplified model
problem
\begin{equation}\label{eq:ell_perf_problem_formulation}
\begin{split}
-\Delta u = f
\qquad &\text{in } \Omega_p,\\
u = 0, 
\qquad &\text{on } \partial \Omega_p.
\end{split}
\end{equation}
The weak formulation is 
\begin{equation}
a(u_H,v_H)=f(v_H),
\nonumber
\end{equation}
where
\begin{equation}
a(u_H,v_H)=\int_{\Omega_p} \nabla u_H \cdot \nabla v_H,
\qquad
f(v_H)=\int_{\Omega_p} f v_H.
\nonumber
\end{equation}
Using multicontinuum GMsFEM (see multiscale basis functions), we have
(see \eqref{eq:mcbasis} and \eqref{eq:mcrep})
\begin{equation}
u_H=\sum_i N_i U_i + \sum_i M_i \cdot \nabla U_i,
\ \
v_H=\sum_j N_j V_j + \sum_j M_j \cdot \nabla V_j.
\nonumber
\end{equation}
We refer to
\eqref{eq:Ni_ell_perf} for definition of $N_i$
used in numerical simulations and \eqref{eq:Mi_ell_perf} for the basis functions, which define $M_i$.
Ignoring second-order derivatives of macroscopic variables, we have 
\begin{equation}
\nabla u_H
\approx
\sum_i
\left(
\nabla N_i U_i
+
B_i \nabla U_i
\right),
\qquad
B_i:=N_i I+\nabla M_i.
\nonumber
\end{equation}
Substituting this into the weak formulation, we have 
\begin{equation}
\begin{split}
a(u_H,v_H)=
\sum_K \sum_{i,j}
\int_{K\cap\Omega_p}
\Big[
(\nabla N_i\cdot\nabla N_j)U_iV_j
+
(\nabla N_i\cdot B_j)U_i\nabla V_j
\\
+
(B_i\nabla U_i\cdot \nabla N_j)V_j
+
(B_i\nabla U_i)\cdot(B_j\nabla V_j)
\Big].
\nonumber
\end{split}
\end{equation}
Define coarse-scale effective coefficients:
\begin{equation}
\begin{split}
\alpha_{ij}
:=
\frac{1}{|K|}
\int_{K\cap\Omega_p}
\nabla N_i\cdot\nabla N_j,\quad
\beta_{ij}
:=
\frac{1}{|K|}
\int_{K\cap\Omega_p}
\nabla N_i\cdot B_j,\\
\gamma_{ij}
:=
\frac{1}{|K|}
\int_{K\cap\Omega_p}
B_i^T\nabla N_j,\quad
\theta_{ij}
:=
\frac{1}{|K|}
\int_{K\cap\Omega_p}
B_i^TB_j.
\nonumber
\end{split}
\end{equation}

The forcing term becomes
$f(v_H)
=
\sum_K
\int_{K\cap\Omega_p}
f
\left(
N_jV_j+M_j\cdot\nabla V_j
\right).$
We define
\begin{equation}
F_j
:=
\frac{1}{|K|}
\int_{K\cap\Omega_p}
fN_j,
\qquad
G_j
:=
\frac{1}{|K|}
\int_{K\cap\Omega_p}
fM_j.
\nonumber
\end{equation}

Assuming macroscopic variables are smooth over each coarse block
(as before)
\begin{equation}
\int_\Omega
\Big[
\alpha_{ij}U_iV_j
+
(\beta_{ij}U_i) \cdot \nabla V_j
+
\gamma_{ij} \cdot \nabla U_i V_j
+
\theta_{ji}\nabla U_i\cdot\nabla V_j
\Big]
=
\int_\Omega
\left(
F_jV_j+G_j\cdot\nabla V_j
\right).
\nonumber
\end{equation}

The  multicontinuum macroscopic system is 
\begin{equation}
\alpha_{ij}U_i
-\nabla\cdot
\left(
\theta_{ji}\nabla U_i
\right)
-
\nabla\cdot
\left(
 \beta_{ij}U_i
\right)
+
\gamma_{ij}\cdot\nabla U_i
=
F_j-\nabla\cdot G_j.
\nonumber
\end{equation}

Let us now consider the construction of the auxiliary functions and multiscale basis functions used in our numerical results. As in the heterogeneous case, we define the coarse neighborhood $\omega_l$ (see \eqref{eq:elliptic_omega_l}) and introduce the general notation $\omega_l^{*}$, which may denote $\omega_l$, an oversampled neighborhood $\omega_l^+$, or the global domain $\Omega$. To account for perforations, we additionally define $\omega_l^p = \omega_l \cap \Omega_p$ and $\omega_l^{*p} = \omega_l^{*} \cap \Omega_p$. The auxiliary functions $N_i^l$ are constructed by solving the following problems in $\omega_l^{*}$.
\begin{equation}
\label{eq:Ni_ell_perf}
\begin{split}
&\int_{\omega_l^{*p}} \nabla N_i^l \cdot \nabla v - \sum_{K_d \subset \omega_l^*} \sum_j \Lambda_{ij}^{ld} \int_{K_d \cap \omega_l^{*p}} \psi_j v = 0,\\
&\int_{K_d \cap \omega_l^{*p}} N_i^l \psi_j = \delta_{ij} \int_{K_d \cap \omega_l^{*p}} \psi_j, \quad \forall K_d \subset \omega_l^*,
\end{split}
\end{equation}
where we additionally set $N_i^l = 0$ on $\partial \omega_l^{*p} \setminus \partial \omega_l^{*}$ in accordance with \eqref{eq:ell_perf_problem_formulation}. Here, $\Lambda_{ij}^{ld}$ denote the Lagrange multipliers used to impose the block-average constraints.

Next, we introduce $z_i^l = M_i^l \cdot \nabla \phi_0^l$. Then, using the definition of $M_i^l$, we obtain the following local problems.
\begin{equation}
\label{eq:Mi_ell_perf}
\begin{split}
&\int_{\omega_l^p} \nabla (z_i^l + N_i^l \phi_0^l) \cdot \nabla v - \sum_{K_d \subset \omega_l} \sum_j \Theta_{ij}^{ld} \int_{K_d \cap \omega_l^p} \psi_j v = 0,\\
&\int_{K_d \cap \omega_l^p} z_i^l \psi_j = 0, \quad \forall K_d \subset \omega_l,\\
&z_i^l|_{\partial K_d \cap \overline{\omega}_l^p} = 0, \quad \forall K_d \subset \omega_l,
\end{split}
\end{equation}
where $z_i^l = 0$ on $\partial \omega_l^p \setminus \partial \omega_l$. Here, $\Theta_{ij}^{ld}$ denote the Lagrange multipliers used to impose the zero block-average constraints.

Finally, as in the heterogeneous case, the multiscale basis functions are constructed using \eqref{eq:phi_expansion_with_z}.

\end{document}